\documentclass{article}

\usepackage{amsmath, amsthm, amssymb}
\usepackage{arxiv}
\usepackage{algorithmicx}

\usepackage{graphicx}		
\usepackage{epstopdf}
\usepackage{algorithm}
\usepackage{algpseudocode}
\usepackage[caption=false]{subfig}

\usepackage{mathtools}
\usepackage{array}    
\usepackage{multicol}
\usepackage{multirow}
\usepackage{rotating}
\usepackage{hhline}
\usepackage{stackengine}

\usepackage[utf8]{inputenc} 
\usepackage[T1]{fontenc}    
\usepackage{hyperref}       
\usepackage{url}            
\usepackage{booktabs}       
\usepackage{amsfonts}       
\usepackage{nicefrac}       
\usepackage{microtype}      
\usepackage{lipsum}		
\usepackage{amsopn}

\ifpdf
  \DeclareGraphicsExtensions{.eps,.pdf,.png,.jpg}
\else
  \DeclareGraphicsExtensions{.eps}
\fi

\newtheorem{remark}{Remark}
\newtheorem{prop}{Proposition}
\newcommand{\figref}[1]{Figure~#1}

\newcommand{\sectref}[1]{Section~#1}
\newcommand{\subsectref}[1]{Subsection~#1}

\newcommand{\Equaref}[1]{Equation~#1}

\newcommand{\proref}[1]{Proposition~#1}
\newcommand{\remarkref}[1]{Remark~#1}

\title{A Discontinuous Galerkin Framework for Efficient Simulation of Cyclic Voltammetry Models}

\author{
  Assionvi Hove Kouevi\thanks{Schlumberger Cambridge Research
    .} \\
  Schlumberger Cambridge Research\\
  Cambridge, CB3 0EL, United Kingdom \\
   \And
Debora Campos de Faria \footnotemark[1]
     \\
  Schlumberger Cambridge Research\\
  Cambridge, CB3 0EL, United Kingdom \\
   \And 
   Can Evren Yarman \thanks{Study and Production schlumberger} \\
   Study and Production schlumberger\\
  1 Rue Henri Becquerel, 92140 Clamart, France \\
  \And     
 Gabriel J. Lord\thanks{Radboud University Nijmegen} \\
  Radboud University Nijmegen\\
  Houtlaan 4, 6525 XZ Nijmegen, Netherlands \\
}

\begin{document}
\maketitle

\begin{abstract}


We develop a numerical framework to efficiently simulate the cyclic voltammetry models based on discontinuous Galerkin (DG) spatial discretization combined with linearly implicit Euler time integration. The method is applied to two electrochemical mechanisms: a heterogeneous electron-transfer (ET) reaction and an electron-transfer reaction coupled with a catalytic chemical step (EC'). The governing time-dependent reaction-diffusion equations are expressed in dimensionless form and discretized in space using locally defined orthonormal Legendre polynomial basis functions. This discretization leads to a system of ordinary differential equations that is subsequently advanced in time using the linearly implicit Euler method. The resulting solver, denoted $\text{DG}_{\text{Impl}}$, is evaluated through numerical experiments and compared with MATLAB's $pdepe$ solver in terms of accuracy and computational efficiency. The numerical results demonstrate the convergence of $\text{DG}_{\text{Impl}}$ with respect to spatial and temporal discretizations and show that the method accurately reproduces the voltammetric responses of both ET and EC' mechanisms. For a prescribed relative error, $\text{DG}_{\text{Impl}}$ requires less computational time than $pdepe$ in the numerical experiments considered. These results demonstrate the effectiveness of the proposed DG formulation for the forward simulation of cyclic voltammetry models and provide a numerical foundation for subsequent adjoint-based parameter-estimation and inverse-problem formulations.

\end{abstract}

\section{Introduction}

Cyclic voltammetry (CV) is a widely used electroanalytical technique for probing electron-transfer reactions and coupled chemical processes at electrode–electrolyte interfaces. It has been frequently used by electrochemists for a variety of purposes, including electrochemical sensing and biosensing (e.g. glucose sensors \cite{ELEC1}, gas detectors \cite{ELEC2,ELEC3}, pH meters  \cite{ELEC4}, biosensors \cite{sankar2024biosensors}), energy storage and rechargeable batteries (e.g. solar cells \cite{ELEC8}, batteries \cite{tian2025batteries, ELEC9}, Supercapacitors \cite{mohamad2025supercapacitors}), technological applications (e.g. electroplating \cite{ELEC5}, electrochromic displays \cite{ELEC6}). CV is also an important tool for elucidating reaction mechanisms in electrolysis and electrosynthetic processes \cite{rafiee2024cv}, as well as for characterizing the activity, redox behavior and interfacial evolution of electrocatalysts \cite{zhang2026cv}. These applications are particularly relevant to electrochemical processes of broad technological interest, including hydrogen evolution \cite{yang2023her}, CO$_2$ reduction \cite{ghazi2026co2}, and electrochemical ammonia synthesis \cite{kou2025ammonia}. CV is further used in electrochemical imaging and related characterization applications \cite{ELEC10,lewis2024electrochemical}. These examples illustrate the continuing importance of CV as both a mechanistic and materials-characterization tool across electrochemistry. Comprehensive treatments of the theory and practice of voltammetry can be found, for example, in \cite{CV_ref3,ET}.

\begin{figure}[H]
\centering
\includegraphics[height=6cm,width=10cm]{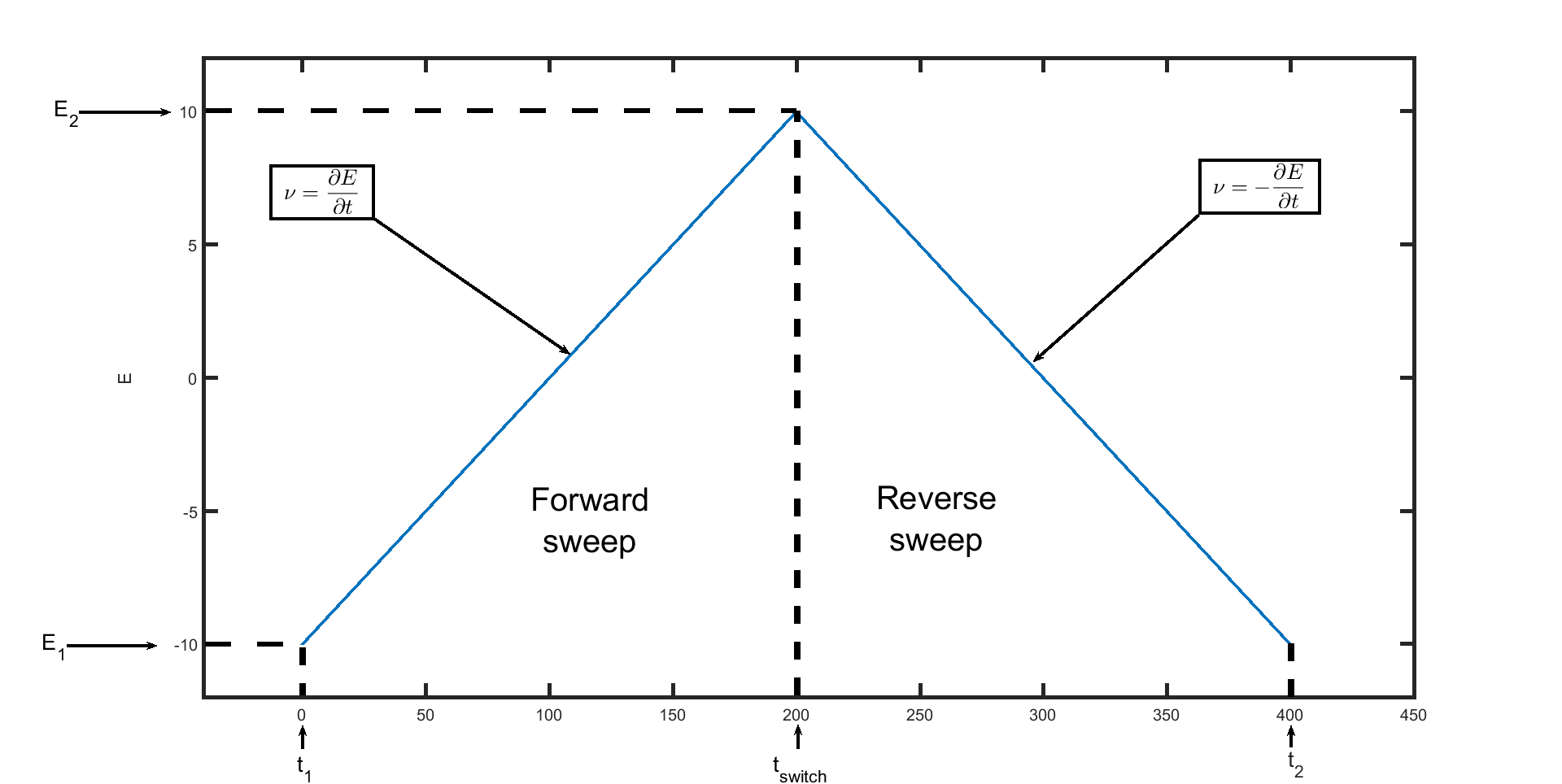}
\caption{\textbf{The waveform of the potential applied during a typical CV experiment}. In this case the initial potential, $E_1 = -10V$, and the vertex potential, $E_2 = 10V$, and the scan rate, $\nu = 0.1Vs^{-1}$.\label{Potential}}
\end{figure}

In its most basic form, CV involves applying a voltage to an electrode in contact with an electrolyte solution, and seeing how the output current changes.
Let us consider for example the heterogeneous  electron transfer (ET)  process at the electrode represented by  \Equaref{\eqref{ETO_reac}}.
\begin{equation}\label{ETO_reac}
Q \xrightleftharpoons[k_b]{k_f} Q^{+} +e^-,
\end{equation}
where $e^-$ denotes an electron; $Q$ and $Q^{+}$ are solution phase species; $k_b$ and $k_f$ are the heterogeneous rate constants given  by
\begin{align}
k_f = k_0\exp\left[ (1-\alpha)\frac{F}{RT}(E - E^0)\right], \quad
k_b = k_0\exp\left[ (-\alpha)\frac{F}{RT}(E - E^0)\right],  \label{kb}
\end{align}
called the  Butler–Volmer kinetics equations \cite{ET}. In \Equaref{\eqref{kb}}, $F$ is the Faraday constant, $T(K)$ is the Kelvin temperature, $R$ is the universal gas constant, $\alpha \in [0,1]$ is the charge transfer coefficient, $k_0(ms^{-1})$ is the standard electron transfer rate constants and  $E^0(V)$ is the formal potential of the half-reaction \Equaref{\eqref{ETO_reac}}. The heterogeneous rate constants $k_b$ and $k_f$  describe  how  the  flux  of  electrons  in  the  electrode solution interface depends on the applied potential.
In a typical CV experiment, the potential of a working electrode is cycled linearly between two potential values at which the oxidation (increasing the potential) and reduction (decreasing the potential) of an electroactive occurs. This is done in two steps, forward and reverse sweeps of the potential as illustrated in \figref{\ref{Potential}}.

\begin{figure}[H]
\centering
\includegraphics[height=6cm,width=10cm]{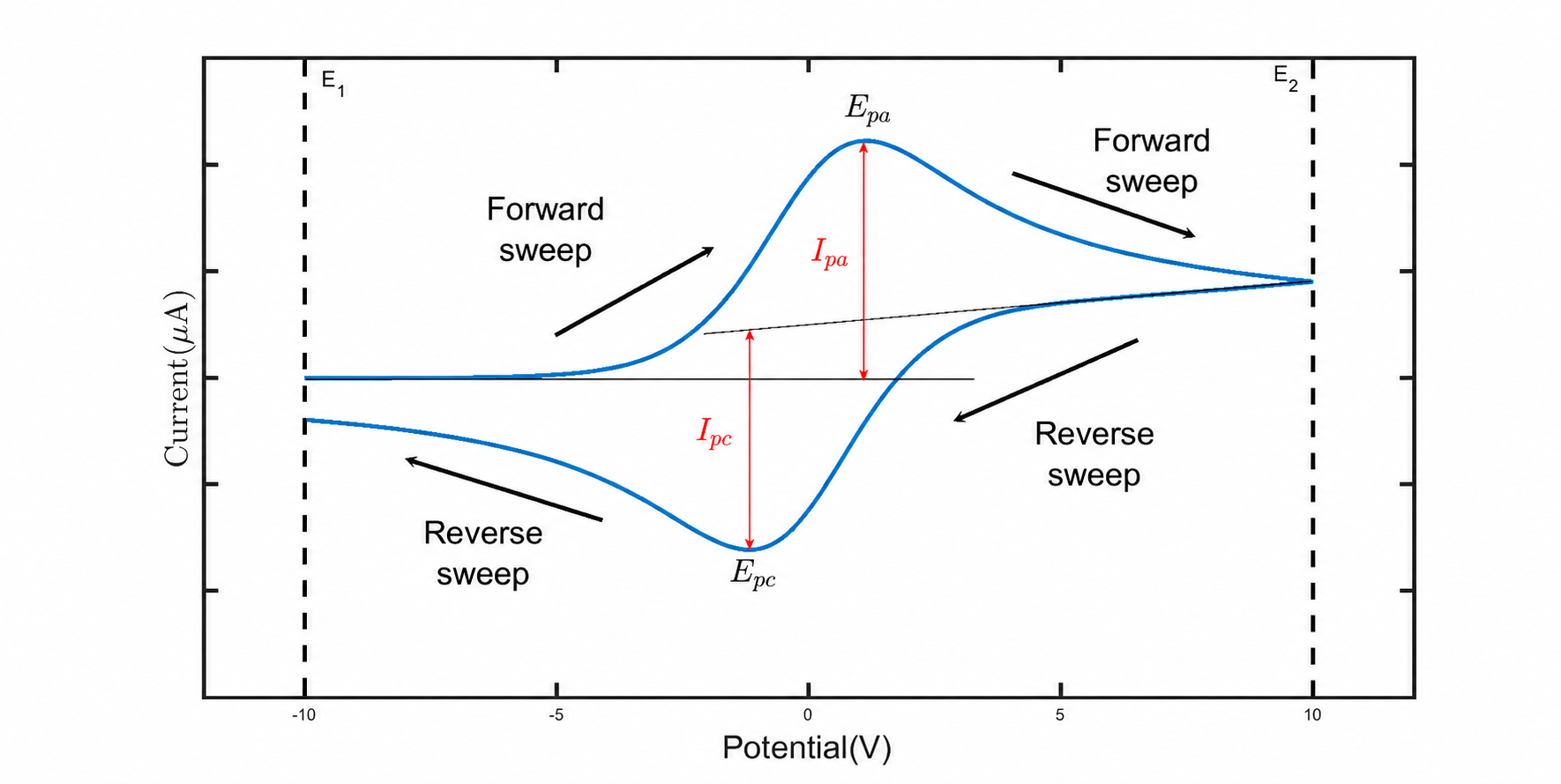}
\caption{\textbf{Representative cyclic voltammogram generated by the applied triangular potential waveform.} The cathodic and anodic peak currents, $I_{pc}$ and $I_{pa}$, occur at the corresponding peak potentials $E_{pc}$ and $E_{pa}$, respectively. The black arrows indicate the directions of the forward and reverse potential sweeps. \label{CVR}}
\end{figure}

On the forward sweep, the potential is given by  $ E = E_1 + \nu t$, for $t\geq 0$, where $\nu$ denotes the scan rate \cite{ET}. At the time $t = t_{switch}= (E_2 - E_1)/ \nu$, the potential reaches the vertex potential $E_2$ and change the direction. On the reverse sweep, the potential is given by $ E =  2E_2 - E_1 -\nu t$, for $t>t_{switch}$.
The process can then be repeated in a periodic, or cyclic manner. Throughout this process the current, $I$ which is proportional to the rate of electron transfer, is recorded. A representative plot of current against potential is presented in \figref{\ref{CVR}}. Here, $I_{pc}$ and  $I_{pa}$ denote the cathodic peak and anodic peak currents, occurring at the corresponding peak potentials $E_{pc}$ and  $E_{pa}$, respectively. Note that in \figref{\ref{CVR}}, as the potential is scanned in the positive direction, the current rises to a peak and then decays. The current depends on two steps in the overall process, the movement of electroactive material to the surface by diffusion and the heterogeneous electron transfer process happening at the electrode-solution interface.
This technique is extremely useful experimentally as the resulting peak shaped signal provides a direct fingerprint of the features of the reduction and oxidation processes. Modeling and analysis of the position, amplitude and shape of the peaks can give important information about the nature of the electrochemical process taking place and about the chemical species involved in the electrochemical reactions. The modelling of the CV experiment requires the definition of the electrical perturbation applied as well as the system under study, in term of mass transport, boundary conditions, heterogeneous electron transfer and any existing homogeneous chemical step. Therefore, the mathematical problem faced in CV is to resolve a system of partial differential equations (PDEs) for example by analytical, semi-analytical or numerical methods. This system of PDEs relate the concentration profiles of the species present in the chemical reaction to measured voltammograms. Numerical methods offer an accurate approximation to the true solution while capturing the key complexity of the process. 

Classical and contemporary electrochemical simulation packages have largely relied on finite difference (FD) formulations. DigiSim employs a fast fully implicit FD algorithm \cite{rudolph1991fifd, rudolph1994simulator, Digisim}, while DigiElch and MECSim use related FD approaches for the numerical simulation of electrochemical reaction transport problems \cite{bieniasz2024computer,digielch2026,kennedy2017mecsim}. More recent open-source tools, including Soft Potato and ElectroKitty, also employ FD based discretizations \cite{rodriguez2022softpotato,vodeb2025electrokitty}. In contrast, general purpose multiphysics environments such as COMSOL Multiphysics typically rely on finite element (FE) discretizations and have been applied to CV problems involving complex electrode geometries and coupled electrochemical transport models \cite{lavacchi2009comsol,mehandzhiyski2025digital,zhang2026comsol}. MATLAB's built-in solver $pdepe$ provides another convenient approach for simulating one-dimensional electrochemical models governed by parabolic PDEs, using a spatial FE discretization combined with the $ode15s$ time integrator \cite{PDEPEMATLAB}. Such tools provide flexible and well-established voltammetric simulators, while also motivating the investigation of alternative spatial and temporal discretizations offering greater local conservation properties, flexibility on nonuniform meshes and formulations that can be systematically adapted to different electrochemical mechanisms.

Alongside the development of simulation software, numerical formulations for CV have continued to evolve toward increasingly detailed physical models and more flexible spatial discretizations. For example, Levey et al.\ developed FE models that couple Butler Volmer kinetics with Nernst-Planck and Poisson equations, first for redox-active monolayers and then for CV at planar electrodes with explicit treatment of ion migration, ohmic potential drop, and the potential dependent electrical double layer \cite{levey2022finite,levey2023simulation}. Nonuniform and adaptive spatial discretizations have also been investigated for systems involving complex or porous electrode geometries. for instance, Tichter et al. used adaptive numerical resolution and nonuniform spatial grids to reduce the computational cost of simulations involving complex electrode geometries \cite{tichter2024rough}. In parallel, physics-informed neural networks (PINNs) have emerged as mesh-free alternatives for solving the PDEs governing voltammetric responses. Chen et al.\ demonstrated PINN-based simulations of one and two dimensional voltammetric problems and reported close agreement with FD as well as analytical reference solutions \cite{chen2022pinn}. A subsequent critical assessment considered a broader set of electroanalytical problems, including CV and chronoamperometry in one to three spatial dimensions, and highlighted both the potential of PINNs and the challenges associated with training, loss weighting, scaling, and domain decomposition for more demanding problems \cite{chen2022critical}. These developments demonstrate that both classical discretization-based solvers and emerging mesh-free approaches remain active areas of research in computational electrochemistry.

Within classical discretization methods, discontinuous Galerkin (DG) methods provide an attractive framework because of their local element formulation, local conservation properties, stability, high order approximation capabilities and flexibility on nonuniform and nonconforming meshes \cite{hesthaven2008nodal,DG1}. Originally developed for transport problems \cite{HDG1}, DG methods now constitute a mature framework for convection diffusion and reaction equations, while their development continues through formulations aimed at improving robustness and computational efficiency. Recent examples include reduced space DG methods for nonlinear convection diffusion reaction equations \cite{hou2024reduced}, hybridizable DG formulations for nonlinear diffusion problems \cite{musch2023hybridizable}, and stable DG schemes for convection dominated convection diffusion reaction equations \cite{boyana2024convergence}. These properties motivate the investigation of DG discretization for the time dependent reaction diffusion systems arising in CV.


In this work, we develop and assess a DG-based forward solver for two representative CV mechanisms: heterogeneous electron transfer (ET) and electron transfer coupled with a catalytic chemical reaction (EC'). The spatial discretization employs locally defined orthonormal Legendre polynomial bases on a nonuniform mesh and is combined with linearly implicit Euler time integration. Particular attention is given to validation against known analytical and physical properties, convergence under temporal and nonuniform spatial mesh refinement, and computational efficiency relative to MATLAB's $pdepe$ solver. Beyond efficient forward simulation, the proposed solver also provides a numerical foundation for subsequent parameter estimation and inverse-problem formulations in CV using adjoint method.

The paper is organized as follows. In \sectref{\ref{sec:main}}, we give a brief review of the two CV models considered. The DG method combined with linearly implicit Euler time integration is applied to the CV models in \sectref{\ref{sec:alg}}. Various experimental
results are produced in \sectref{\ref{sec:experiments1}} to investigate the efficiency of our solver, and the conclusions follow in \sectref{\ref{sec:conclusions}}.


\section{Mathematical Models for Cyclic Voltammetry}
\label{sec:main}
In this section, we describe two cyclic voltammetry mechanisms, namely the Heterogeneous Electron-Transfer (ET) Model and Electron Transfer with a Catalytic Chemical Step (EC') Model before deriving the dimensionless governing equations of the considered models.


\subsection{Heterogeneous Electron-Transfer Model}
\begin{figure}[H]
\centering
\includegraphics[height=5.5cm,width=14cm]{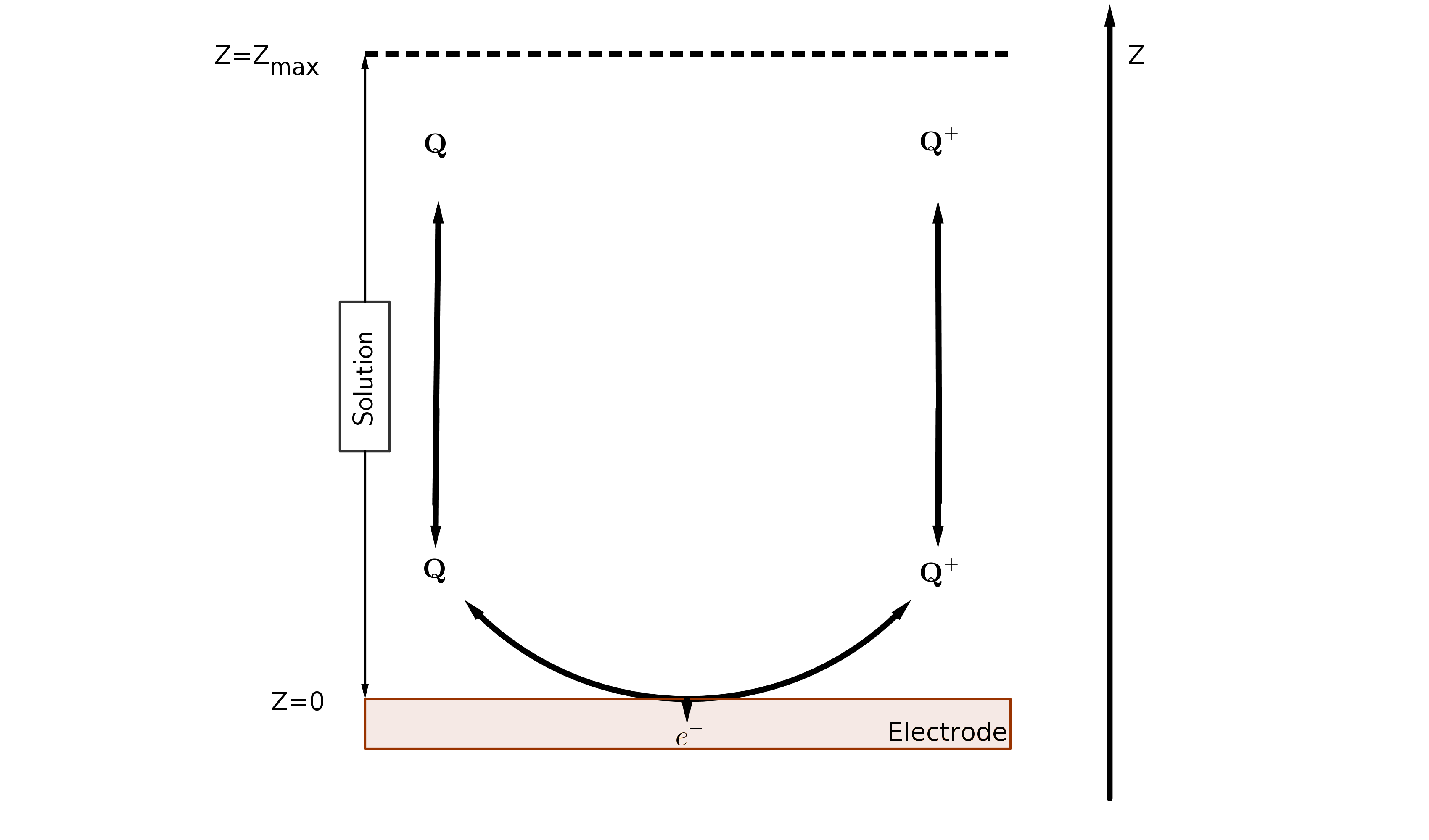}
\caption{\textbf{Schematic representation of the ET mechanism}
The electroactive species $Q$ and $Q^+$ undergo diffusion in the solution, while the reversible electron-transfer reaction takes place at the electrode surface $Z=0$. The computational domain extends from the electrode to the bulk-solution boundary at $Z=Z_{\max}$.
\label{ETO_mod}}
\end{figure}

Let us introduce the ET model, where the electroactive species ($Q$) undergo a heterogeneous charge transfer with the working electrode in an electrochemical cell. Schematically, the ET mechanism is presented in \figref{\ref{ETO_mod}}.

The mathematical model  to  describe  ET is derived from the Fick's Law for mass transfer \cite{Fick1, Fick2, Fick3}. Denoting the concentrations of species $Q$ and $Q^+$ by $C_{Q}$ and $C_{Q^{+}}$, respectively, the one-dimensional transport problem is governed by a coupled system of PDEs given by
\begin{subequations}
\begin{align}
\partial_tC_{Q} =& \partial_Z(D_{Q}\partial_ZC_{Q} ),\label{diffusion0}\\
\partial_tC_{Q^{+}} =& \partial_Z(D_{Q^{+}}\partial_ZC_{Q^{+}} ),\label{diffusion10}
\end{align}
\end{subequations}
subject to the boundary conditions
\begin{subequations}
\begin{align}
\partial_ZC_{i}\Bigr|_{Z = Z_{max}} =& 0\;\;\forall i\in \{Q,Q^{+}\},\label{diff10}\\
D_{Q}\partial_ZC_{Q}\Bigr|_{Z = 0} =& (k_fC_{Q} - k_bC_{Q^{+}})\Bigr|_{Z = 0},\label{diff20}\\
D_{Q^{+}}\partial_ZC_{Q^{+}}\Bigr|_{Z = 0}  =& -(k_fC_{Q} - k_bC_{Q^{+}})\Bigr|_{Z = 0},\label{diff30}
\end{align}
\end{subequations}
and the initial condition
\begin{equation}\label{diff40}
C_{i}\Bigr|_{t = 0} =C_i^0\;\;\forall i\in \{Q,Q^{+}\},
\end{equation}
where $D_{Q}$ and $D_{Q^{+}}$ are respectively the diffusion coefficients of $Q$ and $Q^{+}$. The current, $I$, is given by  the flux of electrons passing through the solution-electrode interface
\begin{equation*}
I = -FAD_{Q^+}\partial_Z C_{Q^{+}}\Bigr|_{Z = 0} =FAD_{Q}\partial_Z C_{Q}\Bigr|_{Z = 0},
\end{equation*}
where $A$ denotes the electrode surface area. A practical way to implement and numerically analyze the physical models is through dimensionless analysis.
The  dimensionless  concentration of the species $Q$ and $Q^{+}$  is  defined  as  the  ratio  between  the concentration of each  species ($C_i$) and a reference concentration, which we take to be the sum of the initial concentrations ($C_{Q}^0 + C_{Q^{+}}^0$). Then, we have
\begin{equation}\label{conc_dim}
\tilde{C}_i = \frac{C_i}{C_{Q}^0 + C_{Q^{+}}^0},\;\forall i\in \{Q,Q^{+}\}.
\end{equation}
For the simulations, we consider $C_{Q^{+}}^0 = 0$, then we have $\tilde{C}_{Q^{+}}^0 = 0$ and $\tilde{C}_{Q}^0 = 1$.
The dimensionless time, $\tilde{t}$, is defined by 
\begin{equation*}
\tilde{t} = \frac{t}{\tau}\;,\qquad \tau = \frac{RT}{F\nu}.
\end{equation*}
where $\tau$ represents  the  time  needed  to  increase or decrease  the  potential  by  one thermal voltage. Similarly, the space coordinates can be made dimensionless by using  the  diffusion  length $\delta$, defined by 
\begin{equation*}
\delta = \sqrt{D_{Q}\tau},
\end{equation*}
as  the  reference  parameter.  It  represents  the  distance  in which the electroactive species moves within the reference time scale $\tau$. Then,  the  dimensionless  coordinate, $z$, is  defined  by  the  ratio  of  the  dimensioned coordinate $Z$ and the diffusion length as
\begin{equation*}
z = \frac{Z}{\delta}.
\end{equation*}

The  diffusion  coefficients $D_{Q}$ and $D_{Q^{+}}$ are  normalized  by  the  diffusion  coefficient  of  the  electroactive species $Q$  to  calculate  a  dimensionless  diffusion  coefficient 
\begin{equation}\label{Diff_dim}
\tilde{D}_i = \frac{D_i}{D_{Q}},\;\forall i\in \{Q,Q^{+}\},
\end{equation}
for  each species involved in the chemical reaction shown in \Equaref{\eqref{ETO_reac}}.
Thus, the dimensionless diffusion coefficient $\tilde{D}_{Q}$ is always equal to one ($\tilde{D}_{Q} = 1$). 
The  dimensionless  electron  transfer  rate  constant, defined by \Equaref{\eqref{rate_dim}}, represents the  level  of electrochemical reversibility  of  the  system
\begin{equation}\label{rate_dim}
K_0 = k_0 \sqrt{\frac{\tau}{D_{Q}}}.
\end{equation}
The change in the applied potential $E$  is converted to the dimensionless overpotential $P$ by dividing the time coordinate  by time scale $\tau$
\begin{equation}\label{pot_dim}
P = \frac{F}{RT}(E-E^0).
\end{equation}
Therefore, the dimensionless time step and the dimensionless  potential step coincide during the simulations  of the dimensionless system. For the numerical simulation, we will assume that $E^0=0$.
Finally the dimensionless current is given by
\begin{equation}\label{cur_dim}
G =  \dfrac{I}{AC_{m}F\sqrt{\dfrac{F\nu D_{Q}}{RT}}}= \partial_z\tilde{C}_{Q}\Bigr|_{z = 0}.
\end{equation}
\begin{remark}\label{Aoki_approx} 
According to \cite{aoki84}, if the species $Q$ and $Q^+$ have the same diffusion coeficient (i.e. $\tilde{D}_{Q^+}=1$) then the dimensionless cathodic peak current, $G_{pc}$, can be approximated by
\begin{equation}\label{AokiEq}
G_{pc} \approx 0.446 + (0.247\sqrt{\alpha} - 0.223)(1-\tanh(0.63\log(K_0) + \frac{0.189}{1-\alpha} - 0.219)).
\end{equation}
The errors involved in the approximation \eqref{AokiEq} are $1.2\%$ for $\alpha$ values ranging from 0.3 to 0.8 and $3\%$ for $\alpha=0.2$. This shows that the reversibility of electron transfer depends on the parameters, namely $K_0$ and $\alpha$.
\end{remark}

Substituting  the  dimensionless  parameters  in Equations \eqref{diffusion0} to \eqref{diff40},   we  get  the dimensionless governing equations for the ET model given by
\begin{subequations}
\begin{align}
\partial_{\tilde{t}}\tilde{ C}_{Q} =& \partial_z(\tilde{ D}_{Q}\partial_z\tilde{ C}_{Q} )\label{diffusion},\\
\partial_{\tilde{t}}\tilde{ C}_{Q^{+}} =& \partial_z(\tilde{ D}_{Q^{+}}\partial_z\tilde{ C}_{Q^{+}} ),\label{diffusion1}
\end{align}
\label{eq:RDE-ETO}
\end{subequations}
subject to the following boundary conditions
\begin{subequations}
\begin{align}
\partial_z\tilde{ C}_{i}\Bigr|_{z = z_{max}} =& 0\;\;\forall i\in \{Q,Q^{+}\}\label{diff1},\\
\tilde{ D}_{Q}\partial_z\tilde{ C}_{Q}\Bigr|_{z = 0} =& (K_f\tilde{ C}_{Q} - K_b\tilde{ C}_{Q^{+}})\Bigr|_{z = 0}\label{diff2},\\
\tilde{ D}_{Q^{+}}\partial_z\tilde{ C}_{Q^{+}}\Bigr|_{z = 0}  =& -(K_f\tilde{ C}_{Q} - K_b\tilde{ C}_{Q^{+}})\Bigr|_{z = 0}\label{diff3},
\end{align}
\label{eq:BC-ETO}
\end{subequations}
and the following initial condition
\begin{equation}\label{diff4}
\tilde{ C}_{i}\Bigr|_{\tilde{t} = 0} =\tilde{C}_i^0\;\;\forall i\in \{Q,Q^{+}\}.
\end{equation}
The heterogeneous electron transfer rate constants, from Equation \eqref{kb}, can then be written in dimensionless form as follows
\begin{equation*}
K_f = K_0\exp[(1-\alpha)P],\qquad
K_b = K_0\exp[(-\alpha)P],
\end{equation*}
where the dimensionless potential, $P$, in term of the dimensionless time, is given by
\begin{align}
P = \left\lbrace
\begin{array}{ll}
P_1 + \tilde{t} , \; 0\leq \tilde{t}\leq \tilde{t}_\lambda\\
P_2 -  (\tilde{t} - \tilde{t}_\lambda) , \; \tilde{t}_\lambda\leq \tilde{t}\leq 2\tilde{t}_\lambda
\end{array} 
\right.,\qquad
\tilde{t}_\lambda = P_2-P_1,
\label{eq:P-ETO}
\end{align}
with $P_1$ and $P_2$ respectively the dimensionless initial and reverse potential, computed from \Equaref{\eqref{pot_dim}}.

\begin{prop}\label{pro_sum}
Let $\tilde{S}$ be the sum of the dimensionless concentration of the species $Q$ and $Q^+$ (i.e. $\tilde{S}= \tilde{ C}_{Q} + \tilde{ C}_{Q^{+}}$). Due to the mass conservation law, at any time $\tilde{t}$, we have
\begin{equation}\label{Eto_equa0}
R(\tilde{t}) := \frac{1}{z_{max}}\int_0^{z_{max}}\tilde{S}(\tilde{t}, z)\;dz = 1.
\end{equation}
Moreover, if the species $Q$ and $Q^+$ have the same diffusion coefficient, i.e.  $\tilde{D}_{Q^+} = 1$, then for any time $\tilde{t}$ and any position $z$, we have
\begin{equation}\label{Eto_equa1}
\tilde{S}(\tilde{t},z) = 1.
\end{equation}
\end{prop}

\begin{proof}[Proof of proposition 1:]
Owing to the mass conservation law, we have $R(\tilde{t}) = R(0)$ for all $\tilde{t}$. If $ \tilde{D}_{Q^+} = 1$, then according to the dimensionless governing  Equations \eqref{diffusion} to \eqref{diff4}, we have  
\begin{equation}\label{sum_eq}
\partial_{\tilde{t}}\tilde{S} = \partial_z(\partial_z\tilde{S} ), \quad
\partial_z\tilde{ S}\Bigr|_{z\in \{0,z_{max}\}} = 0,\;\;
\tilde{S}\Bigr|_{\tilde{t}=0} = \tilde{ C}_{m} =1,
\end{equation}
thus the result given by Equations \eqref{Eto_equa0} and \eqref{Eto_equa1}.
\end{proof}

\subsection{Electron Transfer with a Catalytic Chemical Step Model}
\begin{figure}[H]
\centering
\includegraphics[height=5cm,width=14cm]{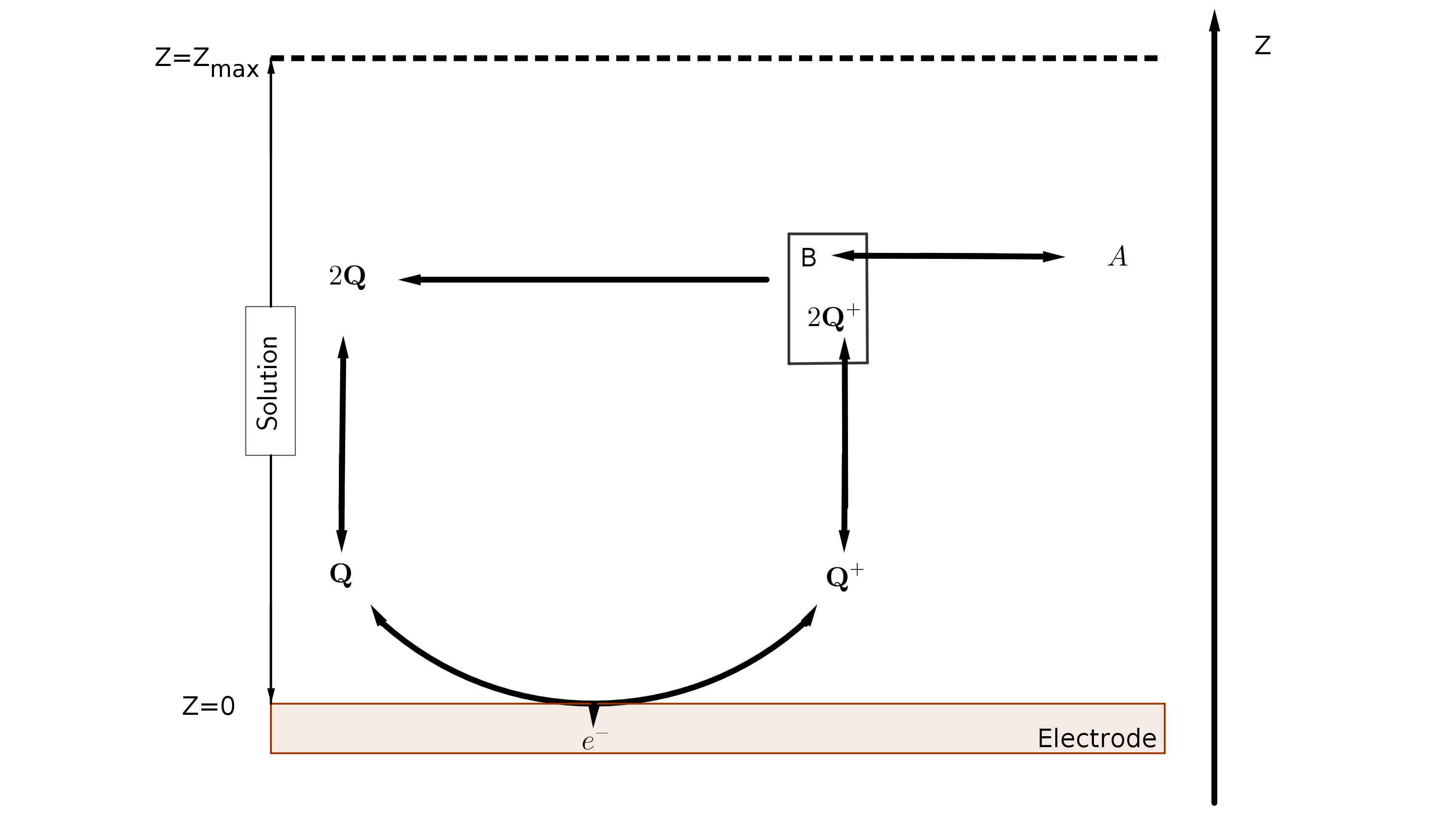}
\caption{\textbf{Schematic representation of the EC' mechanism.}
The species $Q$ and $Q^+$ undergo diffusion and interconversion through heterogeneous electron transfer at the electrode surface $Z=0$. The oxidized species $Q^+$ is additionally involved in a catalytic reaction with $B$, producing $Q$ and $A$ in the solution phase. The domain extends to the bulk-solution boundary at $Z=Z_{\max}$.
\label{EC_mod}}
\end{figure}

We now introduce an EC' model of four species defined by the following chemical reactions
\begin{equation}\label{EC}
Q \xrightleftharpoons[k_b]{k_f} Q^{+} +e^-, \quad B + 2Q^{+} \xrightharpoonup[]{k_r} 2Q  ,\quad A \xrightleftharpoons[k_{A}^-]{k_{A}} B .
\end{equation}
The EC' model can be interpreted as electron transfer followed by a catalytic reaction and an equilibrium reaction where $k_A$ is equal to $k_A^{-}$ and species $A$, $B$ only participate in the homogeneous reactions. Schematically, the EC'  mechanism is presented in \figref{\ref{EC_mod}}.
Denoting the concentrations of species $Q$, $Q^+$, $A$ and $B$ by $C_{Q}$, $C_{Q^{+}}$, $C_{A}$ and $C_{B}$, respectively, the one-dimensional EC' model is governed by a coupled system of PDEs given by
\begin{subequations}
\begin{align}
\partial_tC_{Q} =& \partial_Z(D_{Q}\partial_ZC_{Q} )+2k_rC_{Q^{+}}C_{B},\label{EC_eq1}\\
\partial_tC_{Q^{+}} =& \partial_Z(D_{Q^{+}}\partial_ZC_{Q^{+}} )- 2k_rC_{Q^{+}}C_{B},\\
\partial_tC_{B} =& \partial_Z(D_{B}\partial_ZC_{B} ) - k_rC_{Q^{+}}C_{B} + k_{A}(C_{A} -  C_{B}),\\
\partial_tC_{A} =& \partial_Z(D_{A}\partial_ZC_{A} ) - k_{A}(C_{A} -  C_{B}),\label{EC_eq4}
\end{align}
\end{subequations}
where $k_r$ is the electrocatalytic  rate constant, $k_{A}$ is the rate constant of species $A$ being converted to $B$, the diffusion coefficient of the species $A, B,$ $Q$ and $Q^{+}$ are respectively denoted by $D_{A}, D_{B}, D_{Q}$ and $D_{Q^{+}}$. The boundary conditions are given by Equations \eqref{diff10} to \eqref{diff30} for the species, $Q$ and $Q^{+}$, present in the ET model. Since species $A$ and $B$ only participate in the homogeneous reactions, they are subject to the  no  flux  boundary conditions,

\begin{equation}\label{Noflux}
\partial_ZC_{j}\Bigr|_{Z = Z_{max}} = \partial_ZC_{j}\Bigr|_{Z = 0} = 0,\qquad j\in \{A,B\}.
\end{equation}
Initial  conditions  for the species involved in the ET model and of the species $A$ are  given  by \Equaref{\eqref{diff40}}  and  
\begin{equation}\label{Ini00}
C_{A}\Bigr|_{t = 0} = C_{B}\Bigr|_{t = 0} = \frac{1}{2}C_{A^{total}}^0,
\end{equation}
where $C_{A^{total}}^0$ is the total concentration of the species $A$.
Note from the governing Equations \eqref{EC_eq1} to \eqref{EC_eq4} that the voltammogram of EC' model is equivalent to the voltammogram of the ET model when the electrocatalytic  rate constant $k_r = 0$ or when $C_{A^{total}}^0=0$. For consistency of the analysis and to make better comparisons to the ET model, the EC' model shown in Equations \eqref{EC_eq1} to \eqref{Ini00}, will also be examined in dimensionless  form. Therefore, in addition to the dimensionless parameters  of Equations \eqref{conc_dim} to \eqref{pot_dim}, the dimensionless electrocatalytic  rate constant, $K_r$, and the dimensionless rate of reaction, $K_A$, are  defined  as follows
\begin{equation}
K_r = k_r\tau C_m, \qquad K_A = k_A\tau.
\end{equation}

The dimensionless concentration profile of the species involved in the electron transfer followed by a catalytic reaction model \eqref{EC} follow the system of PDEs
\begin{subequations}
\begin{align}
\partial_{\tilde{t}}\tilde{ C}_{Q} =& \partial_z(\tilde{ D}_{Q}\partial_z\tilde{ C}_{Q} )+2K_r\tilde{ C}_{Q^{+}}\tilde{ C}_{B},\label{diffusion22}\\
\partial_{\tilde{t}}\tilde{ C}_{Q^{+}} =& \partial_z(\tilde{ D}_{Q^{+}}\partial_z\tilde{ C}_{Q^{+}} )- 2K_r\tilde{ C}_{Q^{+}}\tilde{ C}_{B},\label{diffusion32}\\
\partial_{\tilde{t}}\tilde{ C}_{B} =& \partial_z(\tilde{ D}_{B}\partial_z\tilde{ C}_{B} ) - K_r\tilde{ C}_{Q^{+}}\tilde{ C}_{B} + K_{A}(\tilde{ C}_{A} -  \tilde{ C}_{B}),\label{diffusion42}\\
\partial_{\tilde{t}}\tilde{ C}_{A} =& \partial_z(\tilde{ D}_{A}\partial_z\tilde{ C}_{A} ) - K_{A}(\tilde{ C}_{A} -  \tilde{ C}_{B}),\label{diffusion52}
\end{align}
\label{eq:RDE-EC'}
\end{subequations}
subject to the dimensionless boundary conditions in Equations \eqref{diff1} to \eqref{diff3} presented for the ET model with the additional dimensionless boundary conditions 
\begin{equation}\label{bisul12}
\partial_z\tilde{ C}_{j}\Bigr|_{z = z_{max}} = \partial_z\tilde{ C}_{j}\Bigr|_{z = 0} = 0,\qquad j\in \{A,B\}.
\end{equation}
We further impose the dimensionless initial conditions presented for ET model in \Equaref{\eqref{diff4}} with the additional dimensionless initial conditions given by
\begin{equation}\label{bisul22}
\tilde{ C}_{A}\Bigr|_{\tilde{t} = 0} = \tilde{C}_{B}\Bigr|_{\tilde{t} = 0} = \frac{1}{2}\tilde{ C}_{A^{total}}^0.
\end{equation}

In summary, to simulate the dimensionless voltammogram in the case of the EC' model (resp. ET model), we need to solve the system of equation shown in \Equaref{\eqref{eq:RDE-EC'}} (resp. \Equaref{\eqref{eq:RDE-ETO}}) subject to the boundaries conditions shown in \Equaref{\eqref{eq:BC-ETO}} and \Equaref{\eqref{bisul12}} (resp. \Equaref{\eqref{eq:BC-ETO}}) and initial conditions shown in \Equaref{\eqref{diff4}} and \Equaref{\eqref{bisul22}} (resp. \Equaref{\eqref{diff4}}). Here, the concentrations of the species $Q$ and $Q^+$ are also in agreement with \proref{\ref{pro_sum}}.

\section{DG and linearly implicit Euler discretization of the CV models}
\label{sec:alg}
Let us denote $\Omega$ be the one dimension spatial domain given by $ [0, z_{max}]$.  We then consider a non uniform partition of $\Omega$ into $n$ elements $I_i = [z_{i-1}, z_i]$ and denote $\mathcal{P}^{k_i}(I_i)$ the set of polynomials of degree $k_i$ on the interval $I_i$. The DG method approximates the solution on each element, $I_i,\;i=\{1,\cdots,n\}$ using polynomials, see \cite{DG1}. 
The step size $h_i$ given by $h_i=z_{i}-z_{i-1}$ for all $i$ in $\{1,\cdots,n\}$ is associated to each element $I_i$. 
We define the test space $V_h$ by 
\begin{equation}\label{fintSpac}
V_h = \{u \in L^2(\Omega)\text{ such that }\;\; u\Bigr|_{I_i} \in \mathcal{P}^{k_i}(I_i),\forall\; i=\{1,\cdots,n\}\}.
\end{equation}
For a given function $u\in V_h$, let us introduce the notation
\begin{equation}
u(z_i^+) = \underset{\epsilon\underset{\epsilon>0}{\rightarrow} 0}{\lim}u(z_i + \epsilon),\qquad u(z_i^-) = \underset{\epsilon\underset{\epsilon>0}{\rightarrow} 0}{\lim}u(z_i - \epsilon).
\end{equation}
At an interior node $z_i$ (i.e. $z_i,\;\;i=\{1,\cdots,n-1\}$), the jump function and the average function of $u$ are respectively defined by
\begin{equation}
[u(z_i)] = u(z_i^-) - u(z_i^+),\qquad \{u(z_i)\} = \frac{1}{2} (u(z_i^-) + u(z_i^+)).
\end{equation}
By convention, the definition of jump and average is extended to the external nodes. At the external nodes $z_0$ and $z_n$, the jump function is defined by
\begin{equation}
[u(z_0)] = -u(z_0^+),\qquad [u(z_n)] = u(z_n^-),
\end{equation}
and the average function is defined by
\begin{equation}
\{u(z_0)\} = u(z_0^+),\qquad \{u(z_n)\} = u(z_n^-).
\end{equation}

For simplicity, we introduce the following notations
$$D = \tilde{ D}_{Q},\;\;D^+ = \tilde{ D}_{Q^{+}},\;\;  C = \tilde{ C}_{Q},\;\;C^+ = \tilde{ C}_{Q^{+}} $$
To derive the weak problem of the ET model, we apply the methodology described in \cite{DG1} to  Equations \eqref{diffusion} and \eqref{diffusion1}. Let us multiply \Equaref{\eqref{diffusion}} by $u\in V_h$ and then integrate by parts on each interval $I_i$ to get
\begin{equation}\label{tot}
\int_{I_i} u\partial_{\tilde{t}}Cdz + \int_{I_i} D\partial_zC \partial_zu dz - D\partial_zC(z_i^-)u(z_i^-) + D\partial_zC(z_{i-1}^+)u(z_{i-1}^+)= 0,
\end{equation}
for all $i$ in $\{1,\cdots,n\}$. By summing  \Equaref{\eqref{tot}} over all values of $i$, we obtain
\begin{equation*}
\sum_{i=1}^{n} \int_{I_i} u\partial_{\tilde{t}}Cdz  + \sum_{i=1}^{n} \int_{I_i} D\partial_zC \partial_zu dz - \sum_{i=1}^{n-1}[D\partial_zC(z_i)u(z_i)]
+ D\partial_zC(z_{0}^+)u(z_{0}^+) - D\partial_zC(z_n^-)u(z_n^-)= 0.
\end{equation*}
When we combine the above equation with the boundary conditions shown in Equations \eqref{diff1} and \eqref{diff2}, we obtain
\begin{equation}
\sum_{i=1}^{n} \int_{I_i} u\partial_{\tilde{t}}Cdz + \sum_{i=1}^{n} \int_{I_i} D\partial_zC \partial_zu dz - \sum_{i=1}^{n-1}[D\partial_zC(z_i)u(z_i)] 
+(K_fC(z_{0}) - K_bC^+(z_{0}))u(z_{0}) = 0.\label{weak_P}
\end{equation}
For all $X_j,Y_j, j\in \{1,2\}$, we have
\begin{equation}\label{identi}
X_1Y_1 - X_2Y_2 = \frac{1}{2}(Y_1 + Y_2)(X_1 -X_2) + \frac{1}{2}(X_1 +X_2)(Y_1 - Y_2).
\end{equation}
Letting $X_1= D\partial_zC(z_i^-), X_2= D\partial_zC(z_i^+), Y_1 = u(z_i^-)$ and $Y_2 = u(z_i^+)$, we obtain
\begin{equation}\label{ident}
[D\partial_zC(z_i)u(z_i)] = \{D\partial_zC(z_i)\} [u(z_i)]  + \{u(z_i)\} [D\partial_zC(z_i)],
\end{equation}
for all interior nodes $z_i, \;i \in \{1,\cdots,n-1\}$. Substituting the \Equaref{\eqref{ident}} into \eqref{weak_P} and noting that the exact solution $C$ satisfies $[D\partial_zC(z_i)] = 0$ for all nodes $z_i$, $i\in \{1,\cdots,n-1\}$, we obtain
\begin{equation*}
\sum_{i=1}^{n} \int_{I_i} u\partial_{\tilde{t}}Cdz + \sum_{i=1}^{n} \int_{I_i} D\partial_zC \partial_zu dz - \sum_{i=1}^{n-1}\{D\partial_zC(z_i)\} [u(z_i)]
 +(K_fC(z_{0}) - K_bC^+(z_{0}))u(z_{0}) = 0.
\end{equation*}
We now note that the exact solution also satisfies $[C(z_i)] = 0$  for all interior nodes $z_i$, $i\in \{1,\cdots,n-1\}$. Therefore, $C$ satisfies 
\begin{align}
\sum_{i=1}^{n} &\int_{I_i} u\partial_{\tilde{t}}Cdz + \sum_{i=1}^{n} \int_{I_i} D\partial_zC \partial_zu dz+(K_fC(z_{0}) - K_bC^+(z_{0}))u(z_{0}) \nonumber\\
&+ \sum_{i=1}^{n-1}-\{D\partial_zC(z_i)\} [u(z_i)] +s \{D\partial_zu(z_i)\} [C(z_i)] + D\sigma_i[C(z_i)][u(z_i)]
  = 0,\label{last}
\end{align}
where $\sigma_i = \sigma_0/h_{i,i+1}, \;\;h_{i,i+1} = \max(h_i,h_{i+1})$.

By applying the same methodology to \Equaref{\eqref{diffusion1}}, the concentration $C^+$ satisfies 
\begin{align}
\sum_{i=1}^{n} &\int_{I_i} u\partial_{\tilde{t}}C^+dz + \sum_{i=1}^{n} \int_{I_i} D^+\partial_zC^+ \partial_zu dz - (K_fC(z_{0}) - K_bC^+(z_{0}))u(z_{0}) \nonumber\\
&+ \sum_{i=1}^{n-1}-\{D^+\partial_zC^+(z_i)\} [u(z_i)] +s \{D^+\partial_zu(z_i)\} [C^+(z_i)] + D^+\sigma_i[C^+(z_i)][u(z_i)]
  = 0.\label{last1}
\end{align}
Let us define the following bilinear functions
\begin{subequations}
\begin{align}
\mathcal{A}_0(v,u) =& \sum_{i=1}^{n} \int_{I_{i}} uv\;dz,\;\qquad
\mathcal{B}_D(v,u) = D\sum_{i=1}^{n} \int_{I_{i}} \partial_zv \partial_zu \;dz,\label{Mas_M}\\
\mathcal{C}_D^{s,\sigma_0}(v,u) =& D\sum_{i=1}^{n-1}\underbrace{-\{\partial_zv(z_i)\} [u(z_i)] +s \{\partial_zu(z_i)\} [v(z_i)] + \sigma_i[v(z_i)][u(z_i)]}_{=:\;\mathcal{C}_{i}^{s,\sigma_0}(v,u)},\label{bil_S}\\
\mathcal{D}_{K_f}(v,u) =& -K_fv(z_{0})u(z_{0}),\;\qquad
\mathcal{E}_{K_b}(v,u) = -K_bv(z_{0})u(z_{0}).
\end{align}
\end{subequations}
Then,  \Equaref{\eqref{last}} for $C$ and \Equaref{\eqref{last1}} for $C^+$, can be rewritten as follows
\begin{subequations}
\begin{align}
\mathcal{A}_0(\partial_{\tilde{t}} C,u) + \overbrace{\mathcal{B}_D(C,u) + \mathcal{C}_D^{s,\sigma_0}(C,u) - \mathcal{D}_{K_f}(C,u)}^{=:\; \mathcal{A}_{D,K_f}^{s,\sigma_0}(C,u)} + \mathcal{E}_{K_b}(C^+,u) = 0,\label{wf1}\\
\mathcal{A}_0(\partial_{\tilde{t}} C^+,u) + \underbrace{\mathcal{B}_{D^+}(C^+,u) + \mathcal{C}_{D^+}^{s,\sigma_0}(C^+,u)  - \mathcal{E}_{K_b}(C^+,u)}_{=:\; \mathcal{A}_{D^+,K_b}^{s,\sigma_0}(C^+,u)} + \mathcal{D}_{K_f}(C,u) = 0.\label{wf2}
\end{align}
\end{subequations}
Depending on the choices of the parameters $\sigma_0$ and $s$, we obtain different variations of DG methods that have appeared in the literature:
\begin{itemize}
\item If $s = -1$ and $\sigma_0$ is bounded below by a large enough constant, the resulting method is called the symmetric interior penalty Galerkin (SIPG) method, introduced in 1970 by Wheeler \cite{SP2} and Arnold \cite{SP1}. It has also been shown that in one space dimension, the SIPG  method for second order elliptic problems is stable for polynomial orders $k_i\geq 2$ using any stabilization parameter ($\sigma_0 = 0$) \cite{SP3}.
\item If $s = 1$ and $\sigma_0 = 1$, the resulting method is called the non symmetric interior penalty Galerkin (NIPG) method, introduced in 1999 by Riviere, Wheeler and Girault \cite{NSP}.
\end{itemize}
In order not to increase the number of parameters of our system, we use NIPG.

To derive the linear system from the weak formulation of ET model, we first define the local basis function of $\mathcal{P}^{k_j}(I_j)$ from the shifted Legendre polynomials of degree less or equal to $k_j$ on $I_j$. Since the application $z\longmapsto X(z) = 2(z-z_{j-1})/h_j - 1$ is a bijection from $I_j$ onto $I = [-1,1]$,  the shifted Legendre polynomials of degree $r$ on $I_j$ are defined as 
\begin{equation}\label{link_B}
\phi^j_r(z) = P_r\left( X(z)\right), \;\;\forall z\in I_j.
\end{equation}
Here, $P_r$ is the Legendre polynomial of degree $r$ on $I$. 
It is straightforward to see
\begin{equation}\label{mass}
\int_{I_j}\phi^j_r(z)\phi^j_l(z)dz = \frac{h_j}{2} \int_I P_r(x)P_l(x)dx\;\; \text{ and } \;\;\parallel \phi^j_r\parallel_{L^2(I_j)} = \frac{\sqrt{h_j}}{\sqrt{2r+1}}.
\end{equation}

For a given value of $k_j$, we define the set $\mathcal{P}^{k_j}(I_j)$ by
$$\mathcal{P}^{k_j}(I_j) = \text{ span}\left\lbrace  \frac{\phi^j_r}{\parallel \phi^j_r\parallel_{L^2(I_j)}},\;r\in\{0,\cdots,k_j\}\right\rbrace .$$

For efficiency, we do not compute all the basis functions, but instead we only use the Legendre polynomials on the reference interval $I$. The orthonormal global basis functions of the finite space
\begin{align}
V_h \approx \text{span}\left\lbrace  \Phi^j_r,j\in\{1,\cdots,n\},r\in\{0,\cdots,k_j\}\right\rbrace,\nonumber\\
\text{where }\Phi^j_r(z)=
\left\lbrace 
\begin{array}{ll}
\frac{\phi^j_r(z)}{\parallel \phi^j_r\parallel_{L^2(I_j)}},&z\in I_j\\
0,&z\notin I_j
\end{array}
\right. ,\label{extend}
\end{align}
are obtained from the orthonormal local basis functions by extending them with zero.
For all $z$ in $\Omega$, we can expand the DG solutions $C$ and $C^+$ and their time derivative as
\begin{align}
C(z,\tilde{t})  =& \sum_{j=1}^{n}\sum_{r=0}^{k_j}\alpha^j_r(\tilde{t})\Phi^j_r(z),\qquad 
\partial_{\tilde{t}} C(z,\tilde{t})  = \sum_{j=1}^{n}\sum_{r=0}^{k_j}d_{\tilde{t}}\alpha^j_r(\tilde{t})\Phi^j_r(z),\label{interp1}\\
C^+(z,\tilde{t})  =& \sum_{j=1}^{n}\sum_{r=0}^{k_j}\alpha^{+,j}_r(\tilde{t})\Phi^j_r(z),\qquad 
\partial_{\tilde{t}} C^+(z,\tilde{t})  = \sum_{j=1}^{n}\sum_{r=0}^{k_j}d_{\tilde{t}}\alpha^{+,j}_r(\tilde{t})\Phi^j_r(z),\label{interp2}
\end{align}
where the time-dependant coefficients $\alpha^j_r$ and $\alpha^{+,j}_r $ are the unknown variables to be solved for. Plugging Equations \eqref{interp1} and \ref{interp2} into the variational formulations shown in Equations \eqref{wf1} and \eqref{wf2} for all $i \in\{1,\cdots,n \}, \; l\in \{0,\cdots,k_i \}$, we obtain
\begin{align*}
\sum_{j=1}^{n}\sum_{r=0}^{k_j}d_{\tilde{t}}\alpha^j_r\mathcal{A}_0(\Phi^j_r,\Phi^i_l) +\sum_{j=1}^{n}\sum_{r=0}^{k_j}\alpha^j_r\mathcal{A}_{D,K_f}^{s,\sigma_0}(\Phi^j_r,\Phi^i_l) + \sum_{j=1}^{n}\sum_{r=0}^{k_j}\alpha^{+,j}_r\mathcal{E}_{K_b}(\Phi^j_r,\Phi^i_l) = 0,\\
\sum_{j=1}^{n}\sum_{r=0}^{k_j}d_{\tilde{t}}\alpha^{+,j}_r\mathcal{A}_0(\Phi^j_r,\Phi^i_l) + \sum_{j=1}^{n}\sum_{r=0}^{k_j}\alpha^{+,j}_r\mathcal{A}_{D^+,K_b}^{s,\sigma_0}(\Phi^j_r,\Phi^i_l) + \sum_{j=1}^{n}\sum_{r=0}^{k_j}\alpha^j_r\mathcal{D}_{K_f}(\Phi^j_r,\Phi^i_l) = 0.
\end{align*}

This is a linear system that can be rewritten as  
\begin{equation}\label{syst_L}
\underbrace{\left( 
\begin{array}{cc}
M & 0\\
0 & M
\end{array}
\right)}_{\mathcal{M}}
\underbrace{\left( 
\begin{array}{c}
d_t\alpha\\
d_t\alpha^+
\end{array}
\right)}_{d_t\chi}
 + 
\underbrace{\left( 
\begin{array}{cc}
A_{D,K_f}^{s,\sigma_0}&E_{K_b}\\
D_{K_f}&A_{D^+,K_b}^{s,\sigma_0}
\end{array}
\right) }_{\mathcal{L}^{s,\sigma_0}_{D,D^+,K_f,K_b}}
\underbrace{\left( 
\begin{array}{c}
\alpha\\
\alpha^+
\end{array}
\right)}_{\chi} =
\left( 
\begin{array}{c}
0\\
0
\end{array}
\right),
\end{equation}
where the vectors $\alpha,\;\alpha^+,\;d_t\alpha\;\text{and }d_t\alpha^+$ are given by 
$$\alpha = (\alpha^j_r), \quad \alpha^+ = (\alpha^{+,j}_r),\quad
d_t\alpha = (d_t\alpha^j_r),\quad d_t\alpha^+=(d_t\alpha^{+,j}_r),  $$
and the block matrices $M,\; A_{D,K_f}^{s,\sigma_0},\;A_{D^+,K_b}^{s,\sigma_0},\;E_{K_b},\;D_{K_f}$ are given by
\begin{align*}
M = (\mathcal{A}_0(\Phi^j_r,\Phi^i_l)),\quad  A_{D,K_f}^{s,\sigma_0} =(\mathcal{A}_{D,K_f}^{s,\sigma_0}(\Phi^j_r,\Phi^i_l)),\quad A_{D^+,K_b}^{s,\sigma_0} =(\mathcal{A}_{D^+,K_b}^{s,\sigma_0}(\Phi^j_r,\Phi^i_l)),\\
\;E_{K_b}=(\mathcal{E}_{K_b}(\Phi^j_r,\Phi^i_l)),\quad D_{K_f}=(\mathcal{D}_{K_f}(\Phi^j_r,\Phi^i_l)),
\end{align*}
where $i,j \in\{1,\cdots,n\} $,  $r \in\{0,\cdots,k_j\} $ and $l \in\{0,\cdots,k_i\} $.
The matrices $\mathcal{M}$ and $\mathcal{L}^{s,\sigma_0}_{D,D^+,K_f,K_b}$ are respectively called the mass matrix and the stiffness matrix. Due to the local definition of the basis functions, most entries of their block matrices are equal to zeros. More details on the assembling of the mass and stiffness matrices, in the case of the ET model, is given in Appendix \ref{sec:append}.

To obtain the DG formulation of the EC' model, we use the same methodology as for the DG formulation of the ET model. We apply the DG method  to each equation of the dimensionless governing equations of EC' model (from Equations \eqref{diffusion22} to  \eqref{diffusion52}) and we obtain the system of ODEs 
\begin{equation}\label{OSSM2}
\mathbb{M} d_{\tilde{t}}X + (\underbrace{\mathfrak{S}_0 + \mathfrak{S}_1(t)+ \mathfrak{S}_2}_{\mathfrak{S}(\tilde{t})})X = N(X,K_r),
\end{equation}
where the vector valued $X= (\alpha, \alpha^+, \beta^-, \beta)^T$ is the coupled vector of the component of the dimensionless concentration of the species $Q$($\alpha = (\alpha_r^j)$),  $Q^{+}$($\alpha^+ = (\alpha_r^{+,j})$),  $B$($\beta^-=(\beta^{-,j}_r)$) and $A$($\beta=(\beta^{j}_r)$) in the finite space $V_h$. The associated mass matrix $\mathbb{M}$ and  stiffness matrix $\mathfrak{S}(\tilde{t})$ can be assembled, for $\tilde{ D}_{Q} = 1, \tilde{ D}_{Q^{+}} = D^+,\tilde{ D}_{B} = D^-$ and $\tilde{ D}_{A} = D$   as follows  
\begin{align*}
\mathbb{M} &= \left( 
\begin{array}{c|c|c|c}
M & 0 & 0 & 0\\\hline
0 & M & 0 & 0\\\hline
0 & 0 & M & 0\\\hline
0 & 0 & 0 & M
\end{array}
\right),&
 \mathfrak{S}_0 &= \left( 
\begin{array}{c|c|c|c}
L_0^{s,\sigma_0} & 0 & 0 & 0\\\hline
0&D^+L_0^{s,\sigma_0} & 0 & 0\\\hline
0& 0 & D^-L_0^{s,\sigma_0} & 0\\\hline
0& 0& 0 & DL_0^{s,\sigma_0} 
\end{array}
\right),\\\\
\mathfrak{S}_2 &= K_{A}\left( 
\begin{array}{c|c|c|c}
0 & 0 & 0& 0\\\hline
0 & 0 & 0& 0\\\hline
0 & 0 &   M& - M\\\hline
0 & 0 & -  M&  M
\end{array}
\right),& 
\mathfrak{S}_1(t) &= \left( 
\begin{array}{c|c|c|c}
K_f(t)L1 &-K_b(t)L1 & 0& 0\\\hline
-K_f(t)L1& K_b(t)L1 & 0& 0\\\hline
0 & 0 & 0& 0\\\hline
0 & 0 & 0& 0
\end{array}
\right),
\end{align*}
where the block matrices $M,\;L_0^{s,\sigma_0} $ and $L1$ are defined while assembling the mass and stiffness matrices associated to the DG formulation of the dimensionless governing equations of ET model (see Appendix \ref{sec:append} for more details). So the mass matrix $\mathbb{M}$ is then the $4n_T\times 4 n_T$ identity matrix and $n_T$ the dimension of $V_h$. The source term is given by 
\begin{equation*}
N(X,K_r)^T = K_r(2N_1(X), -2N_1(X), -N_1(X), 0),
\end{equation*}
where $N_1$ is a $1\times n_T$ vector defined as follows
\begin{equation}\label{source}
N_1(X) = (N_{1,r}^j)_{j=1,\cdots,n}^{r = 0,\cdots, k_j}, \quad  N_{1,r}^j = \int_\Omega \tilde{ C}_{Q^{+}}\tilde{ C}_{B}\Phi^j_r(z)dz.
\end{equation}

To simulate the voltammogram, in the case of the ET or EC' model, we respectively have to solve the system of ODEs shown in \Equaref{\eqref{syst_L}} or \Equaref{\eqref{OSSM2}} subject to their initial condition. Since the block matrix $M$ is the identity matrix, both systems of ODEs can be rewritten as follows
\begin{equation}\label{ode}
\dfrac{d}{d\tilde{t}} X =\overbrace{L(\tilde{t})X + N(X)}^{F(X)},
\end{equation}
which subject to the corresponding initial condition can be solved by the linearly Implicit Euler (Impl) method. The main idea behind the linearly Implicit Euler method is to first linearise \Equaref{\eqref{ode}} using the Taylor expansion described in  \cite{Tayl1} for the nonlinear function $N$ in each step at $X^n$. This leads to 
\begin{equation}\label{ode1}
\dfrac{\partial X}{\partial t } = \overbrace{\left(  L^{n+1} + \dfrac{\partial N}{\partial X }(X^n) \right)}^{J_n} X + \overbrace{\left( N(X_n) - \dfrac{\partial N}{\partial X }(X^n)X_n  \right)}^{N^n_n},
\end{equation}
where the matrix $\dfrac{\partial N}{\partial X }(X^n)$ is the Jacobian of the function $N$ at $X^n$. Finally, to complete the Impl scheme for \Equaref{\eqref{ode}}, the backward Euler scheme is applied to the linear part of \Equaref{\eqref{ode1}} to get
\begin{equation}\label{simpl1}
(I-\Delta tJ_n)X^{n+1} = X^{n} + \Delta tN^n_n.
\end{equation}
Then, starting from an initial condition $X^0$, the function $X^n$ can be computed by repeatedly solving \Equaref{\eqref{simpl1}}. In the particular case of $N(X) = 0$, i.e. ET model, \Equaref{\eqref{simpl1}} is reduced to (I-$\Delta tL^{n+1})X^{n+1} = X^{n}$. For more details on the computation of the source term $N$ and its Jacobian in the case of EC' model, see Appendix \ref{sec:append1}. Throughout the resolution of the system of  ODEs \eqref{ode1}, the dimensionless current is simultaneously computed with 
\begin{equation}\label{ccur1}
G(\tilde{t}) = K_f(\tilde{t})\sum_{r=0}^{k_1}(-1)^r\alpha_r^1(\tilde{t})\Gamma_{\{ r\}}^1 - K_b(\tilde{t})\sum_{r=0}^{k_1}(-1)^r\alpha^{+,1}_r(\tilde{t})\Gamma_{\{ r\}}^1,
\end{equation}
where the coefficient $\Gamma_{\{ S\}}^j$ computed by \Equaref{\eqref{norm_inv1}} for given set $S$ and index $j$. Note from \Equaref{\eqref{ccur1}} that the dimensionless current depends only on the concentration of the species $Q$ and $Q^+$ at the boundary $z=0$.


%
%

\section{Results and Discussion\label{sec:experiments1}} In this section, we discuss the numerical method  applied to ET and EC' models in \subsectref{\ref{sec:experiments}} and \subsectref{\ref{NumEC}}, respectively. We compare the accuracy and computational efficiency of our method against the one obtained using an implementation based on MATLAB in-built solver $pdepe$. 
More details on the MATLAB in-built solver $pdepe$ can be found in \cite{pdeFE, PDEPEMATLAB}. The accuracy is measured with respect to a finely resolved solution with respect to the time and space step obtained using $pdepe$ (the finest time and  space steps considered is such that the relative error between the solution obtained using both solvers is less than $10^{-3}$). We measure computational efficiency by CPU time required for a fixed error relative to the finely resolved solution with respect to the time and space step.

\subsection{Numerical experiments for ET model}
\label{sec:experiments}
We start by validating our solver $\text{DG}_{\text{Impl}}$ against some known properties of the ET model, namely \remarkref{\ref{Aoki_approx}} and \proref{\ref{pro_sum}}. We then compare our solver $\text{DG}_{\text{Impl}}$ against MATLAB's in-built solver $pdepe$. 


Unless stated otherwise, we used \begin{equation}\label{ETO_defo_set}
\tilde{D}_{Q} = \tilde{D}_{Q^+} = 1,\; \alpha = 0.5,\; P_1 = -10,\; P_2 = 10,\; K_0 = 20,\; \Delta t = 0.1, k_j=2
\end{equation}
in \Equaref{\eqref{simpl1}} and from  Equations \eqref{eq:RDE-ETO} to \eqref{eq:P-ETO}  during our numerical simulations. 
\subsubsection{Numerical validation of $\text{DG}_{\text{Impl}}$ for ET model \label{ETO_valid}}
\begin{figure}[H]
\centering
\subfloat[$ \tilde{C}_{Q}$ for $ \tilde{D}_{Q^+} = 1$]{\includegraphics[height=2.5cm, width = 3.5cm]{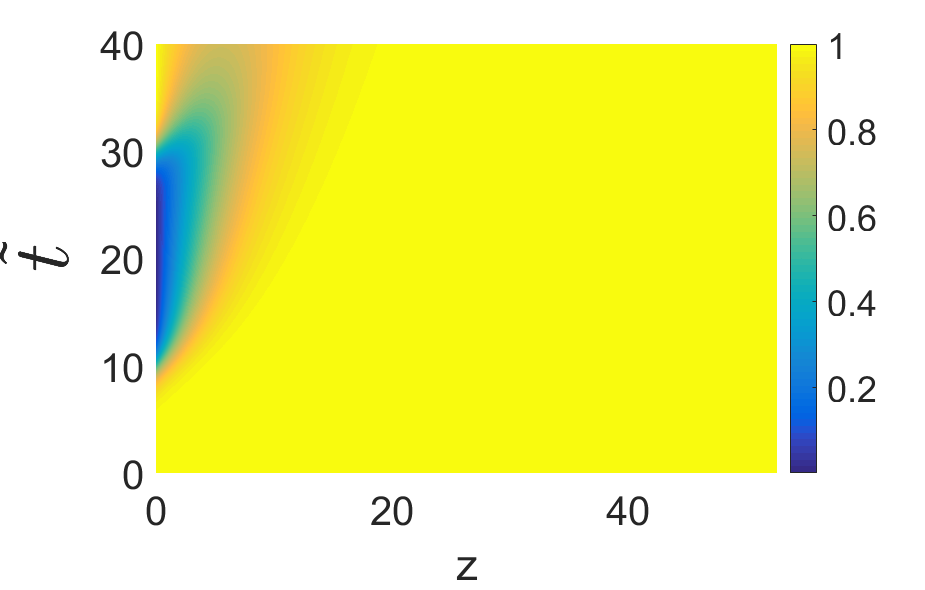}} 
\hspace{0.5cm} 
\subfloat[$ \tilde{C}_{Q^+}$ for $ \tilde{D}_{Q^+} = 1$]{\includegraphics[height=2.5cm, width = 3.5cm]{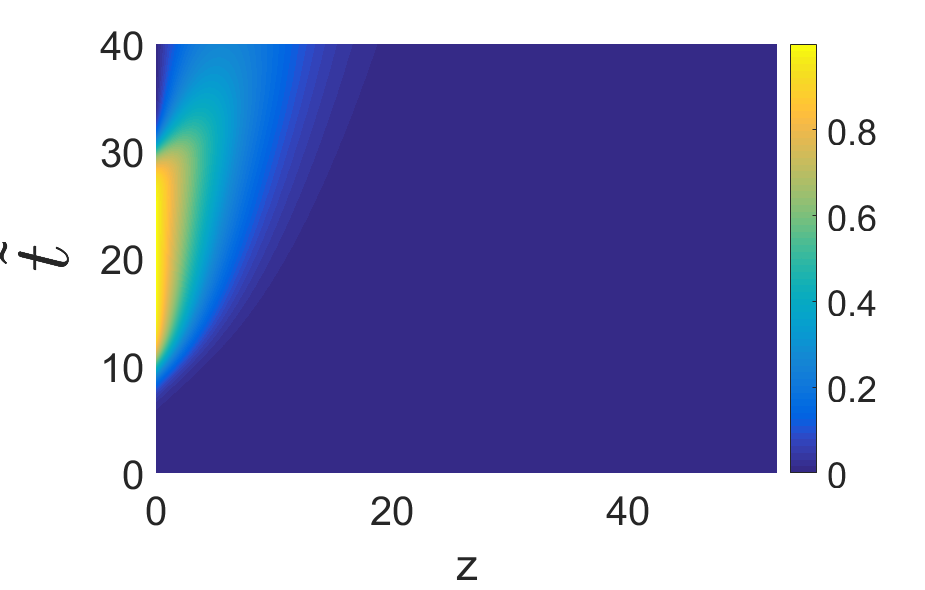}}
\hspace{0.5cm} 
\subfloat[$\tilde{S}$ for $ \tilde{D}_{Q^+} = 1$]{\includegraphics[height=2.5cm, width = 3.5cm]{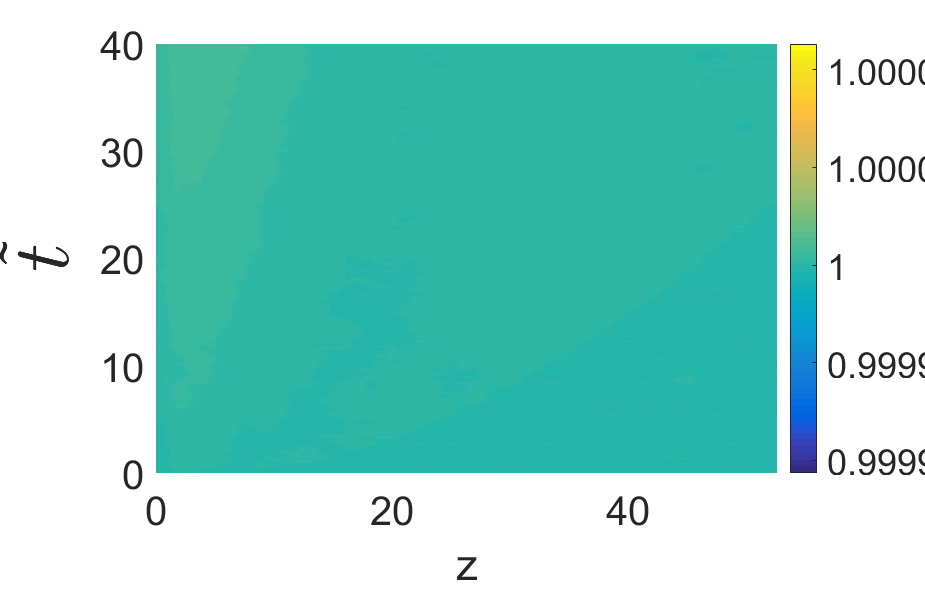}}
\hspace{0.5cm}
\subfloat[$R$ for $ \tilde{D}_{Q^+} = 1$]{\includegraphics[height=2.5cm, width = 3.5cm]{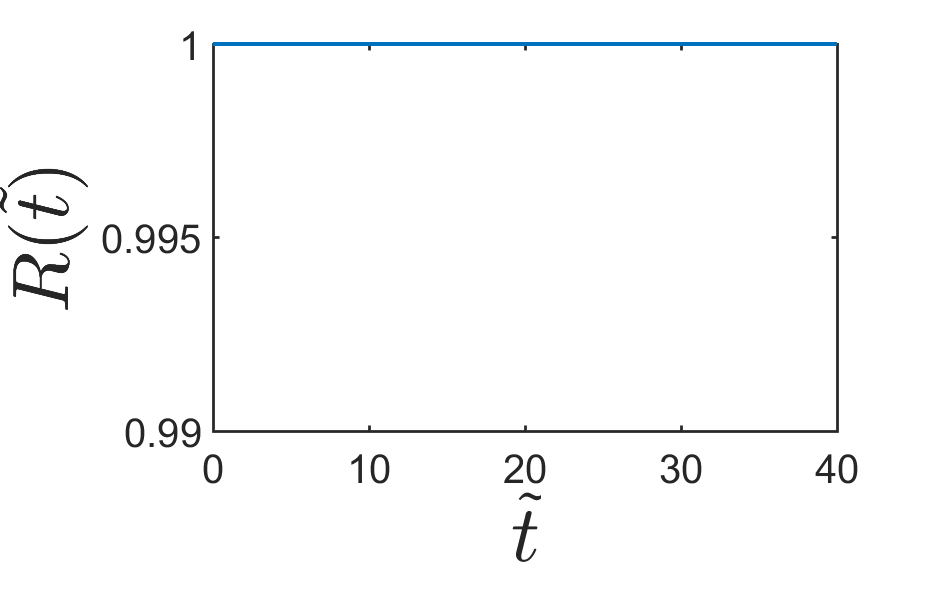}}

\subfloat[$ \tilde{C}_{Q}$ for $ \tilde{D}_{Q^+} = 5$]{\includegraphics[height=2.5cm, width = 3.5cm]{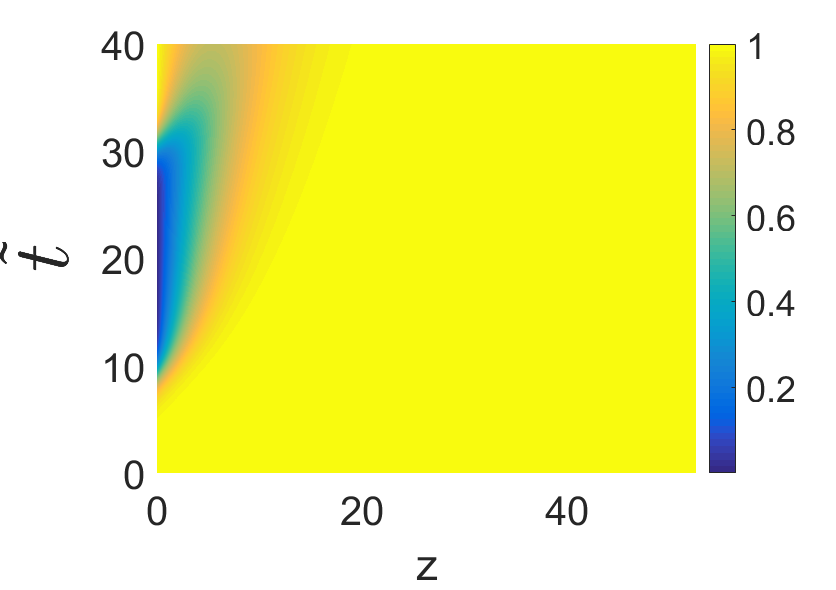}}                 
\hspace{0.5cm}     
\subfloat[$ \tilde{C}_{Q^+}$ for $ \tilde{D}_{Q^+} = 5$]{\includegraphics[height=2.5cm, width = 3.5cm]{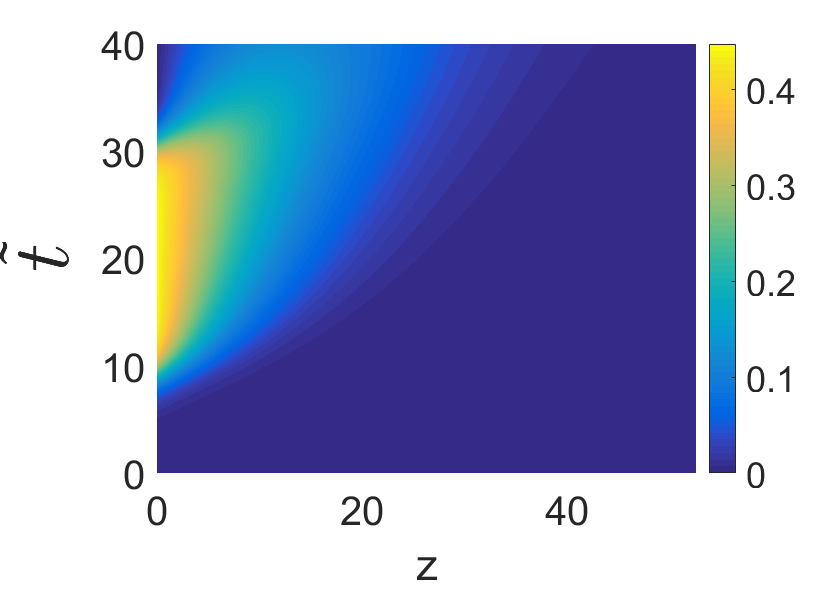}}          
\hspace{0.5cm}
\subfloat[$\tilde{S}$ for $ \tilde{D}_{Q^+} = 5$]{\includegraphics[height=2.5cm, width = 3.5cm]{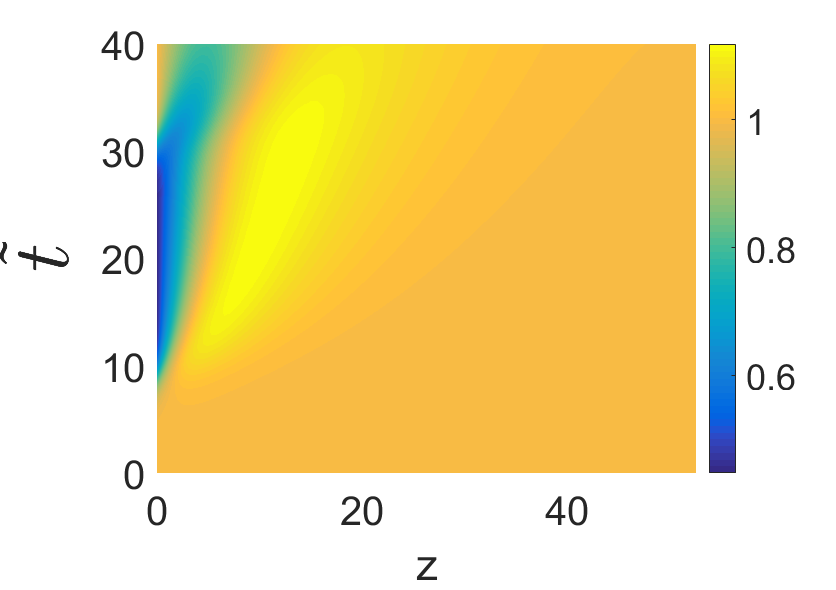}}
\hspace{0.5cm}
\subfloat[$R$ for $ \tilde{D}_{Q^+} = 5$]{\includegraphics[height=2.5cm, width = 3.5cm]{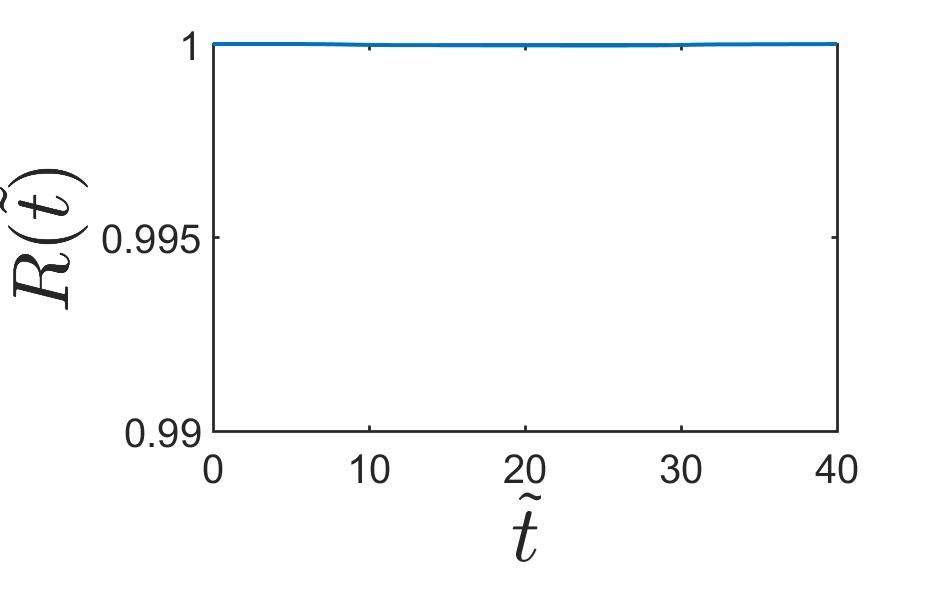}}
\caption{\textbf{Concentration profile.} Simulated concentrations of the species $Q$, $Q^+$,  their sum  $\tilde{S}$, the function $R$ in (a), (b), (c), (d) and  (e), (f), (g), (h), respectively, for $ \tilde{D}_{Q^+} = 1$ (first row) and $ \tilde{D}_{Q^+} = 5$ (second row). Note from (g) and (h) that the mass conservation is satisfied both cases, $\tilde{D}_{Q^+} = 1$  and $\tilde{D}_{Q^+} = 5$ . (c) shows that \Equaref{\eqref{Eto_equa1}} is verified for $ \tilde{D}_{Q^+} = 1$.
\label{ETO_fig1}}
\end{figure}

One way to validate our solver $\text{DG}_{\text{Impl}}$ is to check if the simulated concentrations of the species $Q$ and $Q^+$ in the ET model satisfy  Equations \eqref{Eto_equa0} and \eqref{Eto_equa1} in \proref{\ref{pro_sum}}. For 
$ \tilde{D}_{Q^+} = 1$ and $ \tilde{D}_{Q^+} = 5$,  simulated concentrations of the species $Q$, $Q^+$,  their sum  $\tilde{S}$, the function $R$ are respectively presented in Figure \ref{ETO_fig1} (a), (b), (c), (d) and Figure \ref{ETO_fig1} (e), (f), (g), (h). From Figure \ref{ETO_fig1} (d) and Figure \ref{ETO_fig1} (h) we see that \Equaref{\eqref{Eto_equa0}} is verified for $ \tilde{D}_{Q^+} = 1$ and $ \tilde{D}_{Q^+} = 5$. Similarly, 
Figure \ref{ETO_fig1} (c) and Figure \ref{ETO_fig1} (g) shows that \Equaref{\eqref{Eto_equa1}} is satisfied only for $ \tilde{D}_{Q^+} = 1$ but not for $ \tilde{D}_{Q^+} = 5$, as expected.

\begin{figure}[H]
\centering
\includegraphics[height=6cm, width = 12cm]{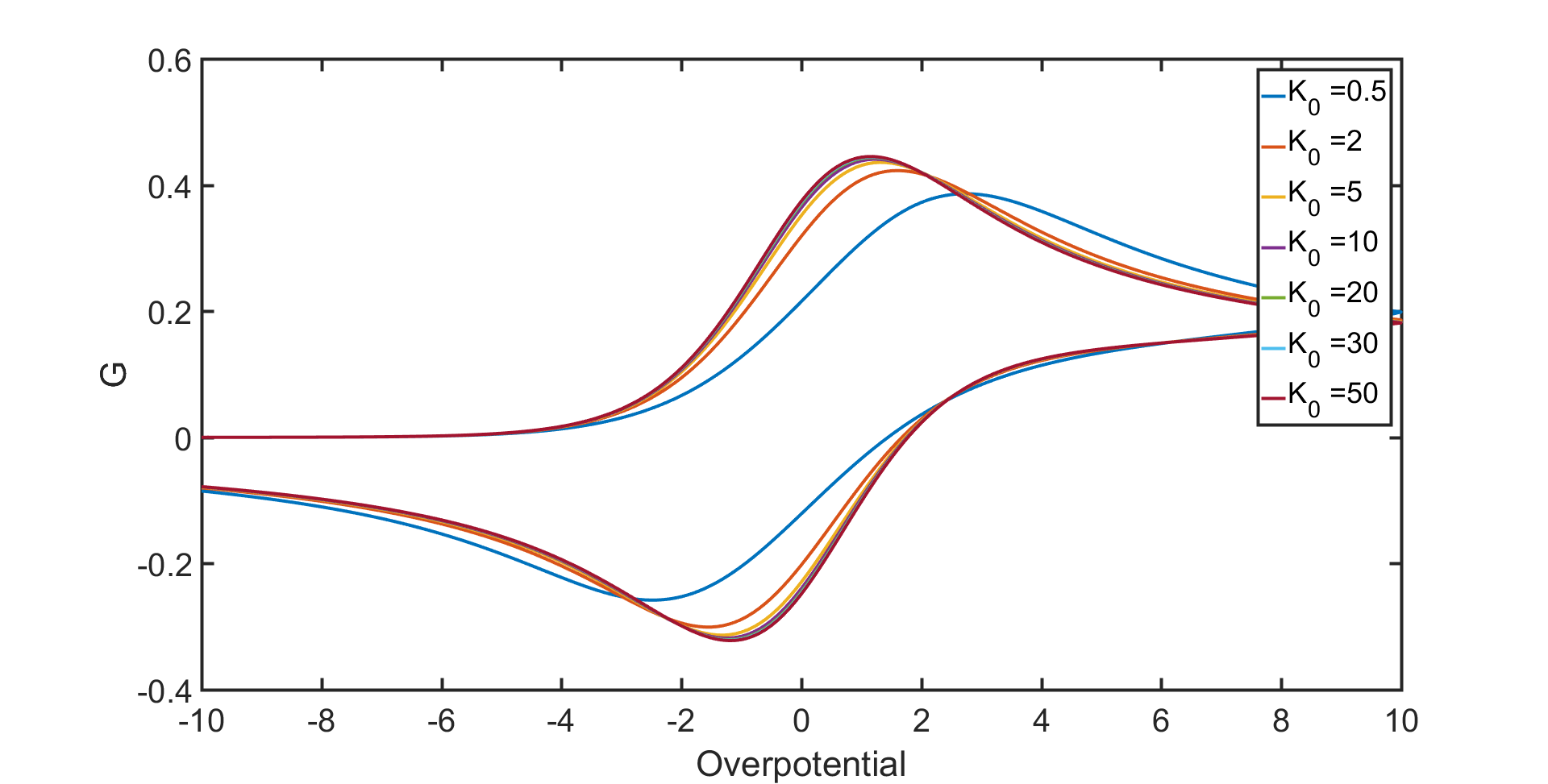}
\caption{Variation of the voltammogram with respect to the dimensionless electron transfer rate constant. This shows that the dimensionless cathodic peak current increases to a certain value as we increase the dimensionless electron transfer rate constant.\label{Var_K0}}
\end{figure}

\begin{figure}[H]
\centering
\subfloat[]{\includegraphics[height=6cm, width = 7cm]{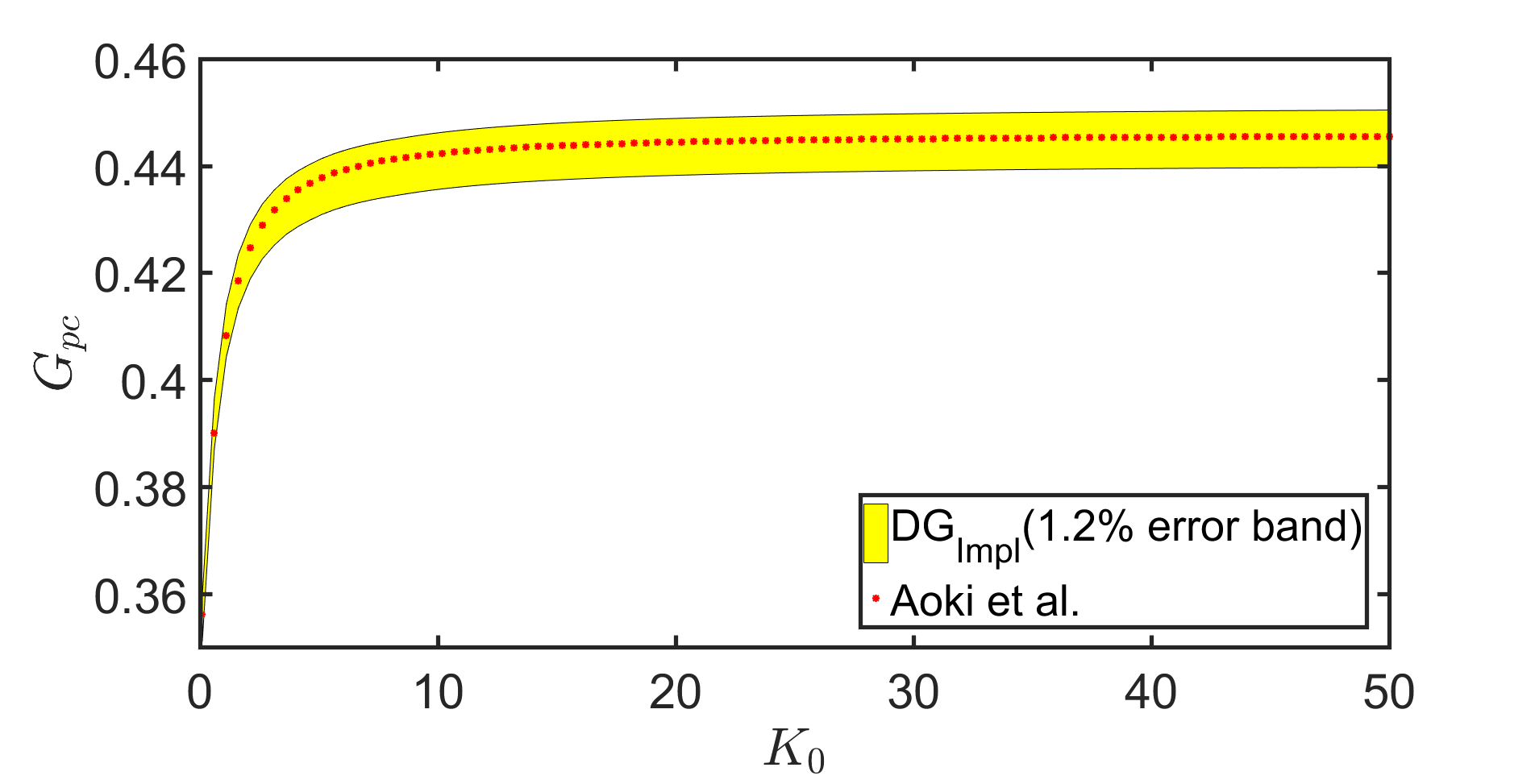}}\hspace{0.5cm}
\subfloat[]{\includegraphics[height=6cm, width = 7cm]{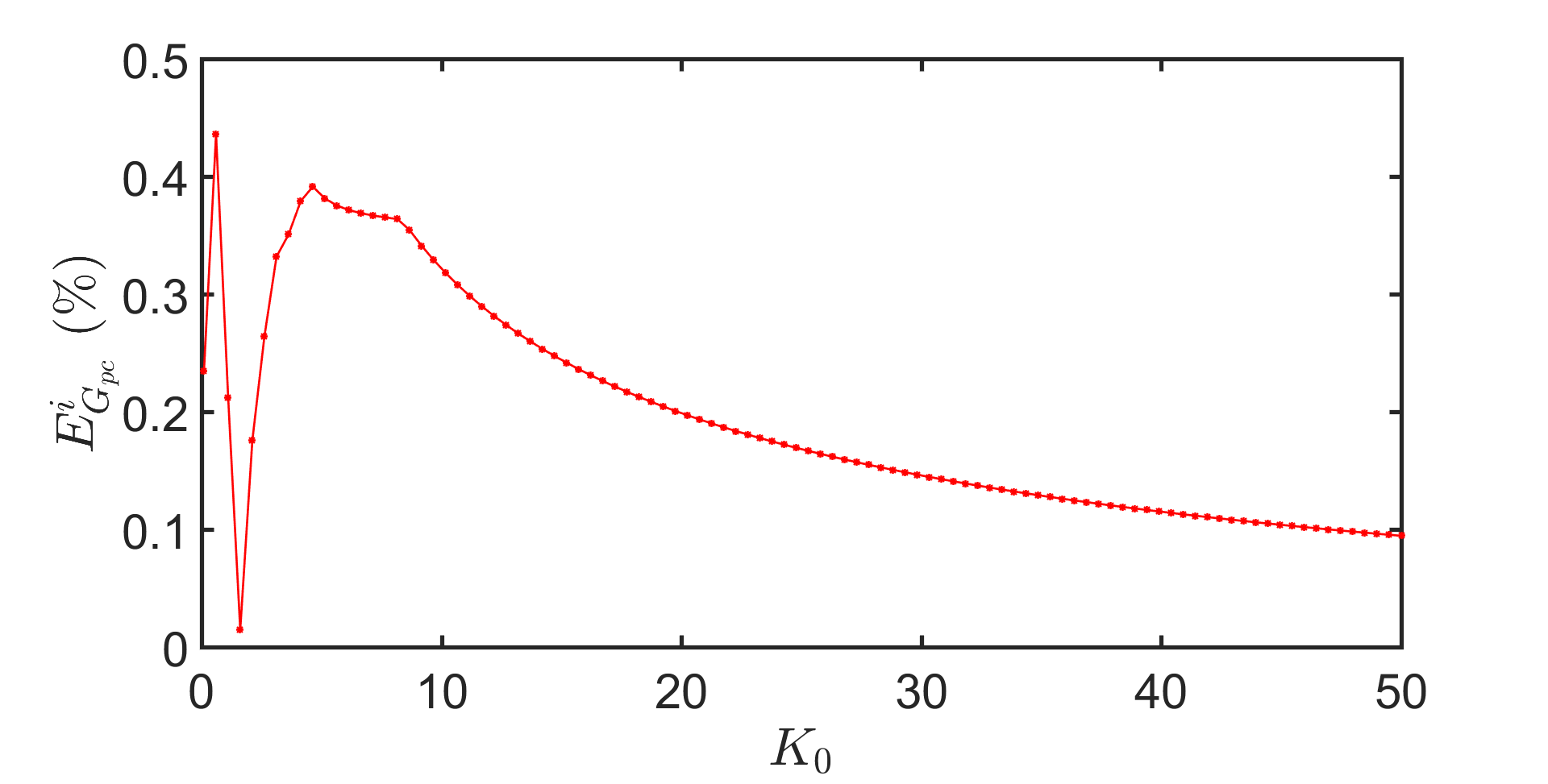}}
\caption{\textbf{Comparison of the simulated and the analytic approximated cathodic peak current.} (a) The dimensionless cathodic peak obtained with the Aoki et al. \cite{aoki84} approximation formula and corresponding $1.2\%$ relative error band around the
$\text{DG}_{\text{Impl}}$ simulated dimensionless cathodic peak current with respect to the dimensionless electron transfer rate constant $K_0$, for $\alpha = 0.5$. (b) The relative error between the simulated and the approximated dimensionless cathodic peak current with respect to $K_0$.\label{Aoki_K0}}
\end{figure}

Another way to validate our solver $\text{DG}_{\text{Impl}}$ is to show that its solutions are in agreement with \remarkref{\ref{Aoki_approx}}. 
To that end, we first simulate the dimensionless current, $G_{K_0}$, with our solver $\text{DG}_{\text{Impl}}$ for the default dimensionless diffusion coefficients $D_{Q^+}$ and various dimensionless electron transfer rate constants $K_0 \in \{0.5, 2,5,10,20,30,50\}$. Figure \ref{Var_K0} presents the dimensionless currents $G_{K_0}$, showing that the dimensionless cathodic peak current increases with $K_0$ and approaches a limiting value, in agreement with \remarkref{\ref{Aoki_approx}}.
Finally, for several dimensionless electron transfer rate constant, $K_0^i\in[0.1, 50]$ and transfer coefficient $\alpha = 0.5$, we simulate the dimensionless cathodic peak currents $\tilde{G}_{pc}^i$  
and $G_{pc}^i$ using respectively our solver $\text{DG}_{\text{Impl}}$ and the approximation \Equaref{\eqref{AokiEq}}. We compute the relative error $E^i_{G_{pc}}$ as follows
\begin{equation*}
E^i_{G_{pc}} = 100\frac{\mid \tilde{G}_{pc}^i - G_{pc}^i\mid}{\tilde{G}_{pc}^i},\quad\text{ for all } K_0^i = 0.1 + i\times\frac{50-0.1}{100}, \;\; \text{ and } i\in \{0,1,\cdots,100\}.
\end{equation*}
As presented in Figure \ref{Aoki_K0}(a), the dimensionless cathodic peak current,  $G_{pc}^i$, remains within the prescribed $1.2\%$ relative error band for all considered values of the dimensionless electron transfer rate constant $K_0^i$. Figure \ref{Aoki_K0}(b) further shows that the corresponding relative error $E^i_{G_{pc}}$ is less than $0.5\%$. These results are consistent with the theoretical behavior established in \remarkref{\ref{Aoki_approx}} and provide additional validation of the $\text{DG}_{\text{Impl}}$ solver.

\subsubsection{Comparison of $\text{DG}_{\text{Impl}}$ and MATLAB's $pdepe$ for ET model}

We compare the performances of $\text{DG}_{\text{Impl}}$ and MATLAB in-built solver $pdepe$ with respect to their accuracy and computational efficiency for simulating the dimensionless current from ET model.

\begin{figure}[H]
\centering
\subfloat[]{\includegraphics[height=5cm, width = 7cm]{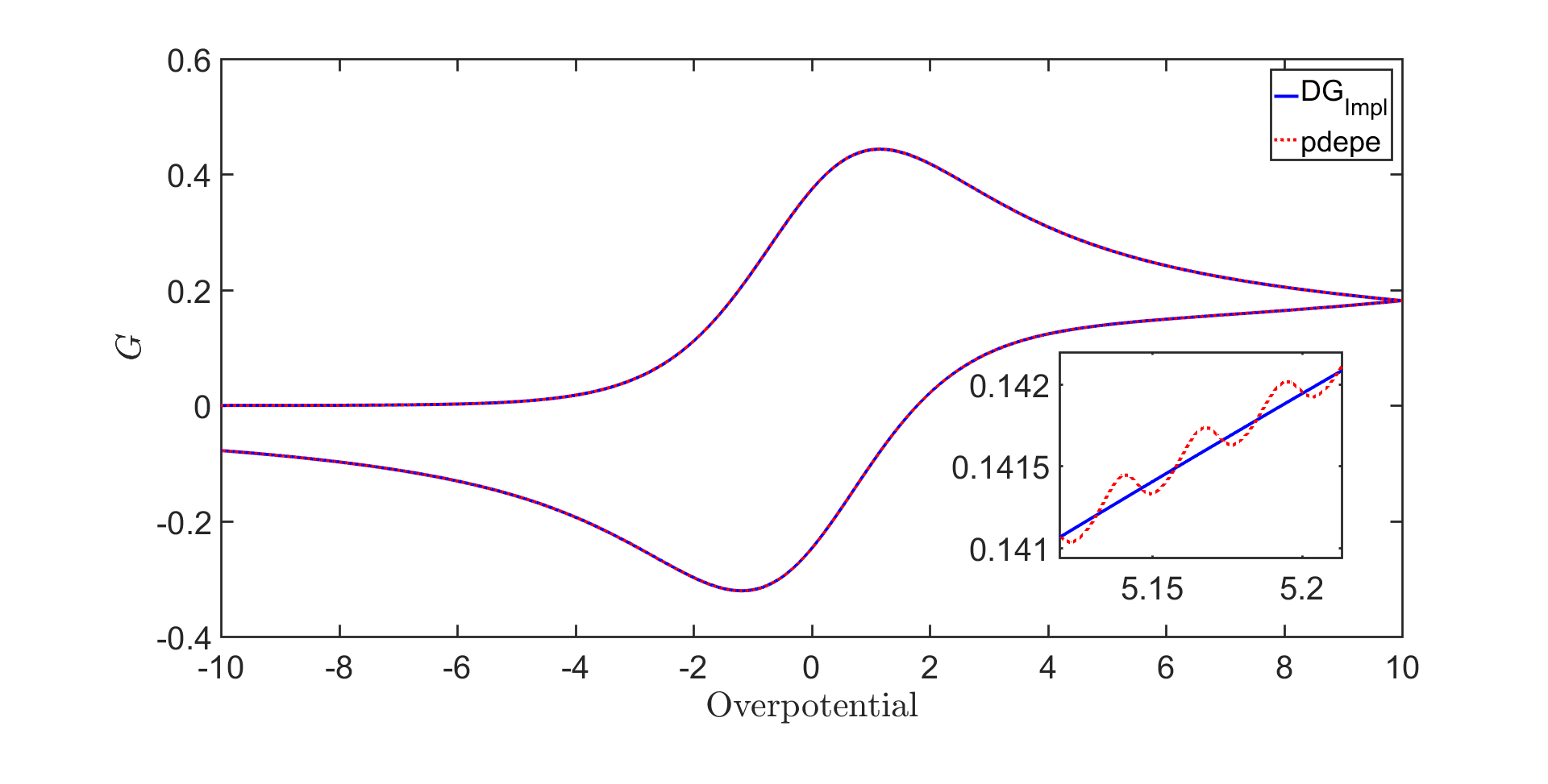}}\hspace{0.5cm}
\subfloat[]{\includegraphics[height=5cm, width = 7cm]{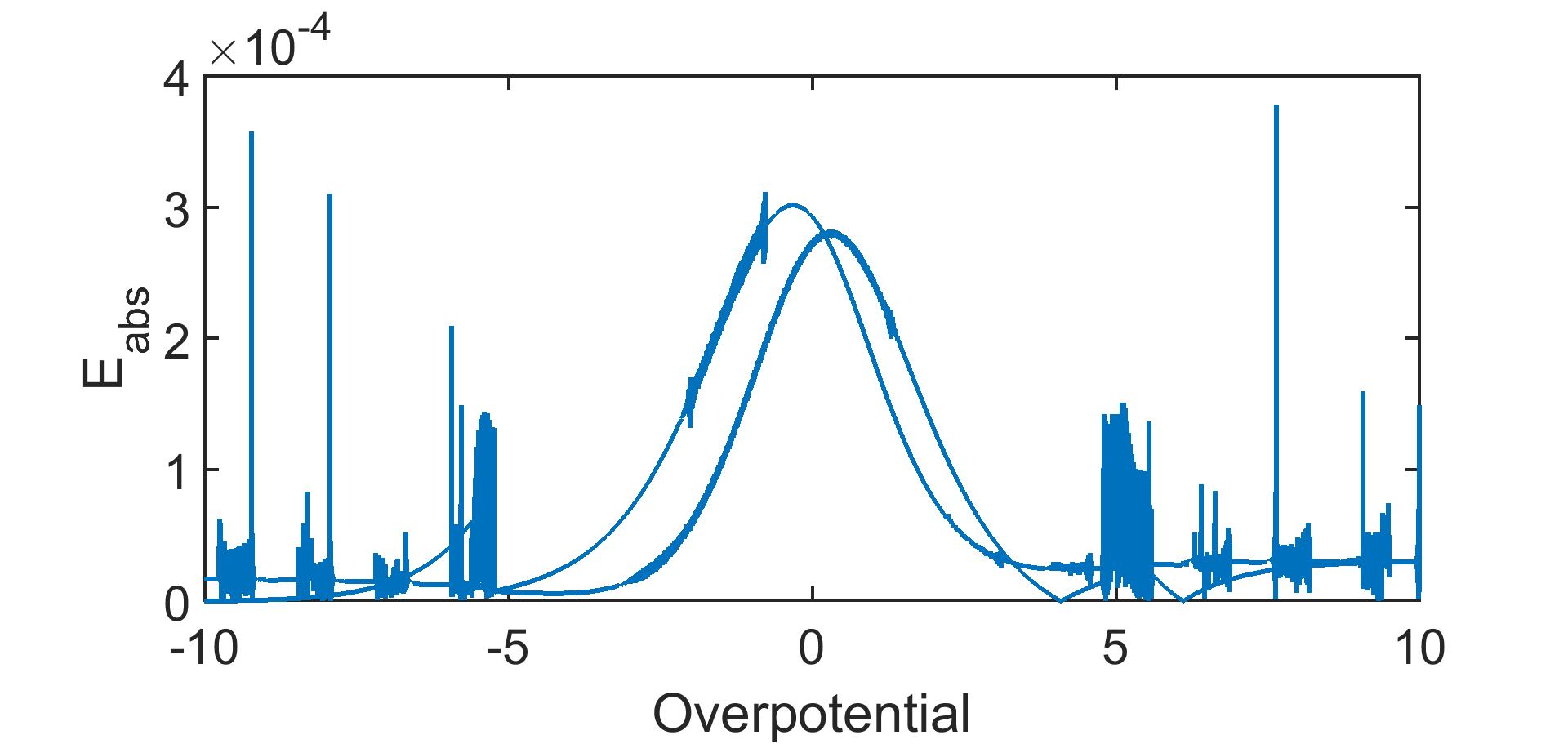}}

\subfloat[]{\includegraphics[height=5cm, width = 7cm]{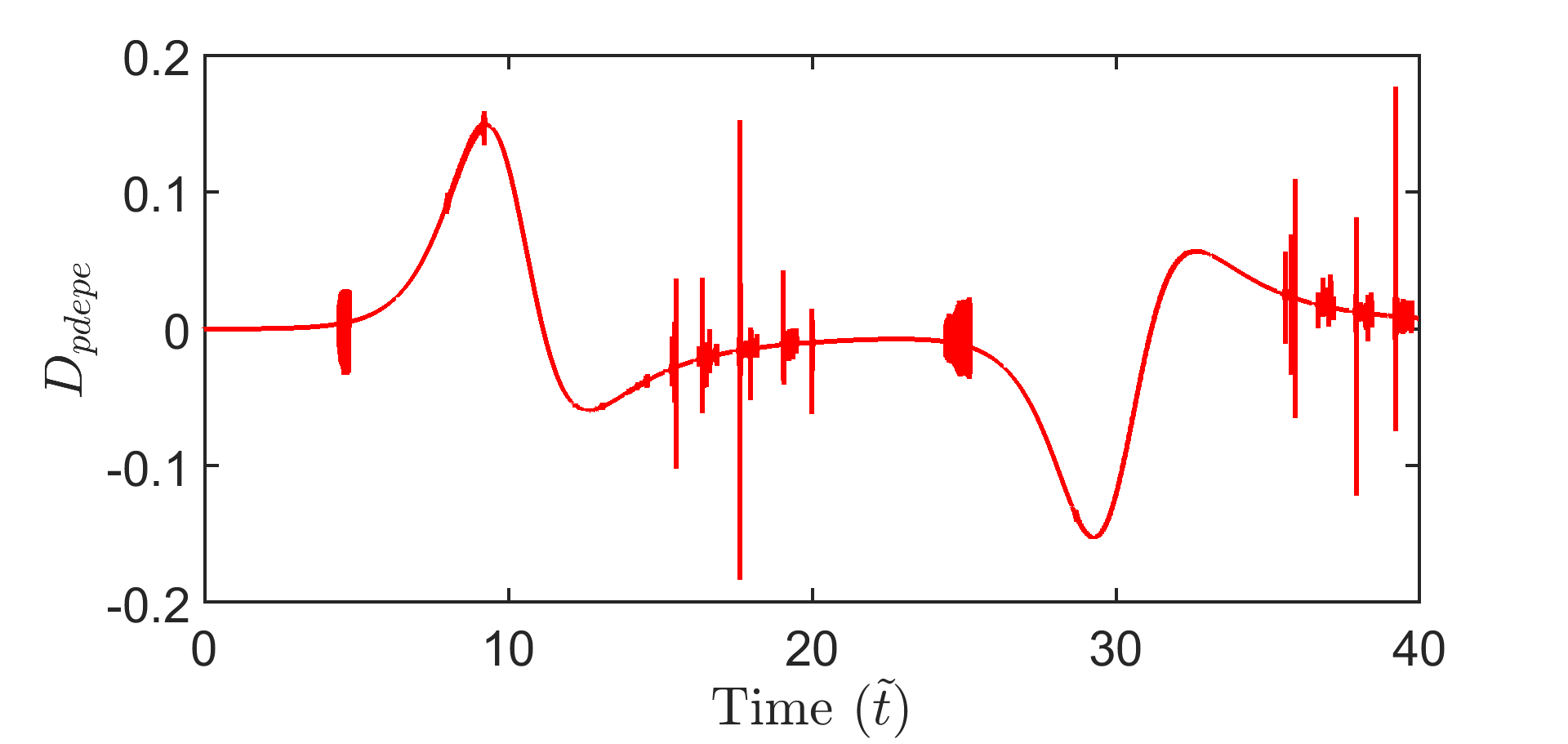}}\hspace{0.5cm}
\subfloat[]{\includegraphics[height=5cm, width = 7cm]{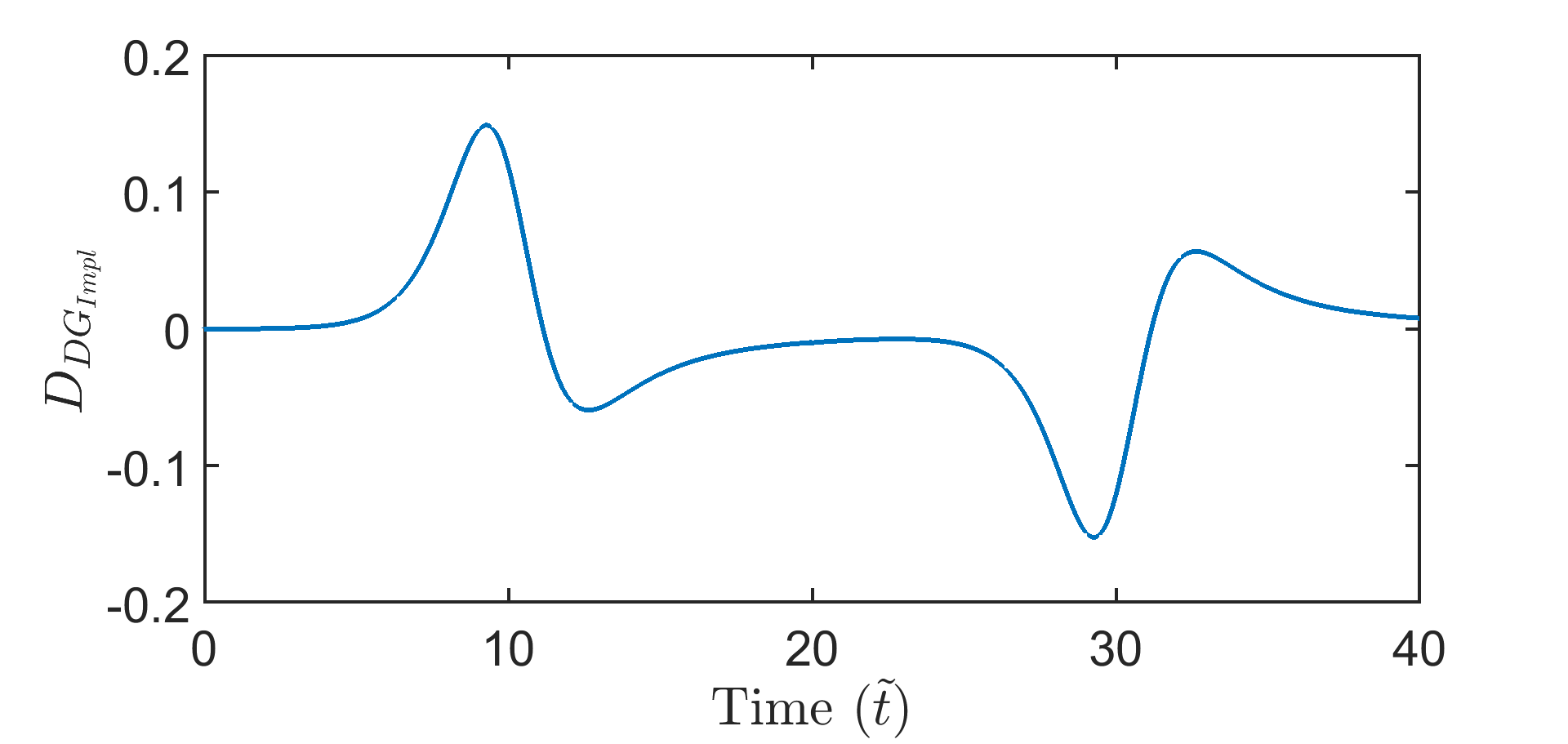}}
\caption{\textbf{Comparison of the simulated dimensionless current from ET model.} 
(a) Dimensionless currents computed using $pdepe$ and $\text{DG}_{\text{Impl}}$ as functions of the overpotential; the inset highlights small oscillations in the $pdepe$ dimensionless current.
(b) Absolute difference between the two simulated dimensionless currents, which remains below $4\times10^{-4}$ over the range considered.
(c)–(d) Approximate time derivatives of the currents computed using $pdepe$ and $\text{DG}_{\text{Impl}}$, respectively, as functions of the dimensionless time $\tilde{t}$. The derivative obtained from $pdepe$ exhibits pronounced oscillations, whereas the $\text{DG}_{\text{Impl}}$ derivative remains comparatively smooth.
\label{compi}}           
\end{figure}

First, we compare the dimensionless current simulated using the solver $r\in \{ pdepe,  \text{DG}_{\text{Impl}}\}$ under the same spatial and temporal discretization settings. Let $G^0_r$  denote the dimensionless current simulated using the solver $r$ for time step $\Delta t^0 = 2.44 \times 10^{-3}$. Simulated currents $G^0_r$ are shown in Figure \ref{compi} (a), while their absolute difference
\begin{equation}\label{abs_diff}
E_{abs} = \mid G^0_{pdepe} - G^0_{\text{DG}_{\text{Impl}}} \mid,
\end{equation}
is presented in Figure \ref{compi} (b). The results show that the absolute difference $E_{abs}$ is less than $4\times 10^{-4}$, indicating close agreement between the two solvers for the ET model. Nevertheless, a closer inspection of Figure \ref{compi} (a) reveals small oscillations in the current obtained with $pdepe$. To further examine these oscillations, we approximate the time derivative of the simulated dimensionless current using a forward finite-difference formula

\begin{equation}\label{varia_comp}
D_{r}(\tilde{t}_n^0) = \dfrac{G^0_r(\tilde{t}_{n+1}^0) - G^0_r(\tilde{t}_{n}^0)}{\tilde{t}_{n+1}^0 - \tilde{t}_n^0},\text{ with } \tilde{t}_n^0 = n\Delta t^0, \;\;\forall n\in  \left\lbrace 0,\cdots,\dfrac{2 t_{\lambda}}{\Delta t^0}\right\rbrace.
\end{equation}
Figure \ref{compi} (c) and (d), present the derivatives $D_{pdepe}$ and $D_{\text{DG}_{\text{Impl}}}$, respectively, as a function of the dimensionless time $\tilde{t}$. The derivative obtained from $pdepe$ exhibits noticeable oscillations, whereas the derivative associated with $\text{DG}_{\text{Impl}}$ remains comparatively smooth. These results indicate that, although the two solvers produce very similar dimensionless currents, $\text{DG}_{\text{Impl}}$ yields a smoother numerical response for the discretization considered.


\begin{figure}[H]
\centering
\subfloat[]{\includegraphics[height=6cm, width = 7cm]{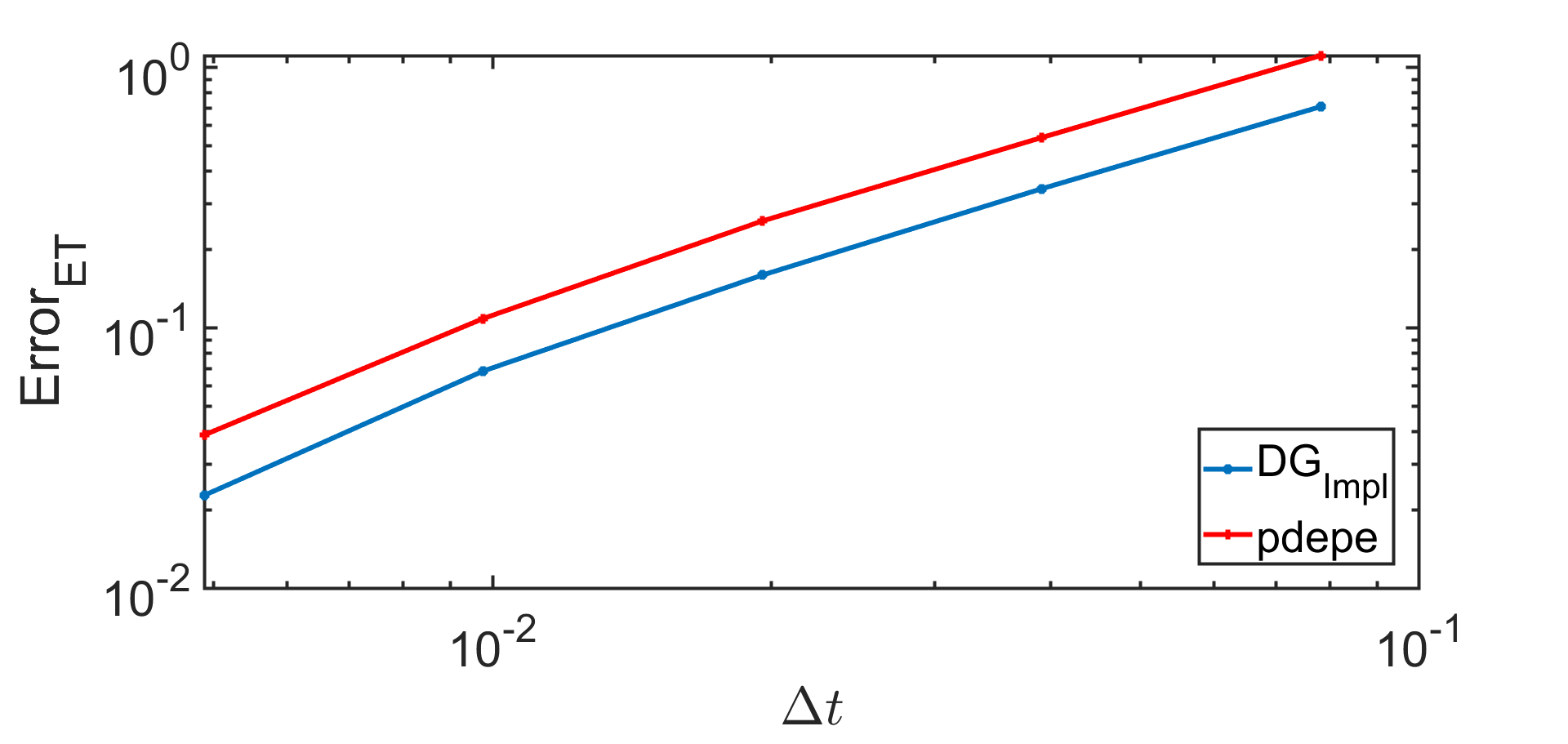}}\hspace{0.5cm}
\subfloat[]{\includegraphics[height=6cm, width = 7cm]{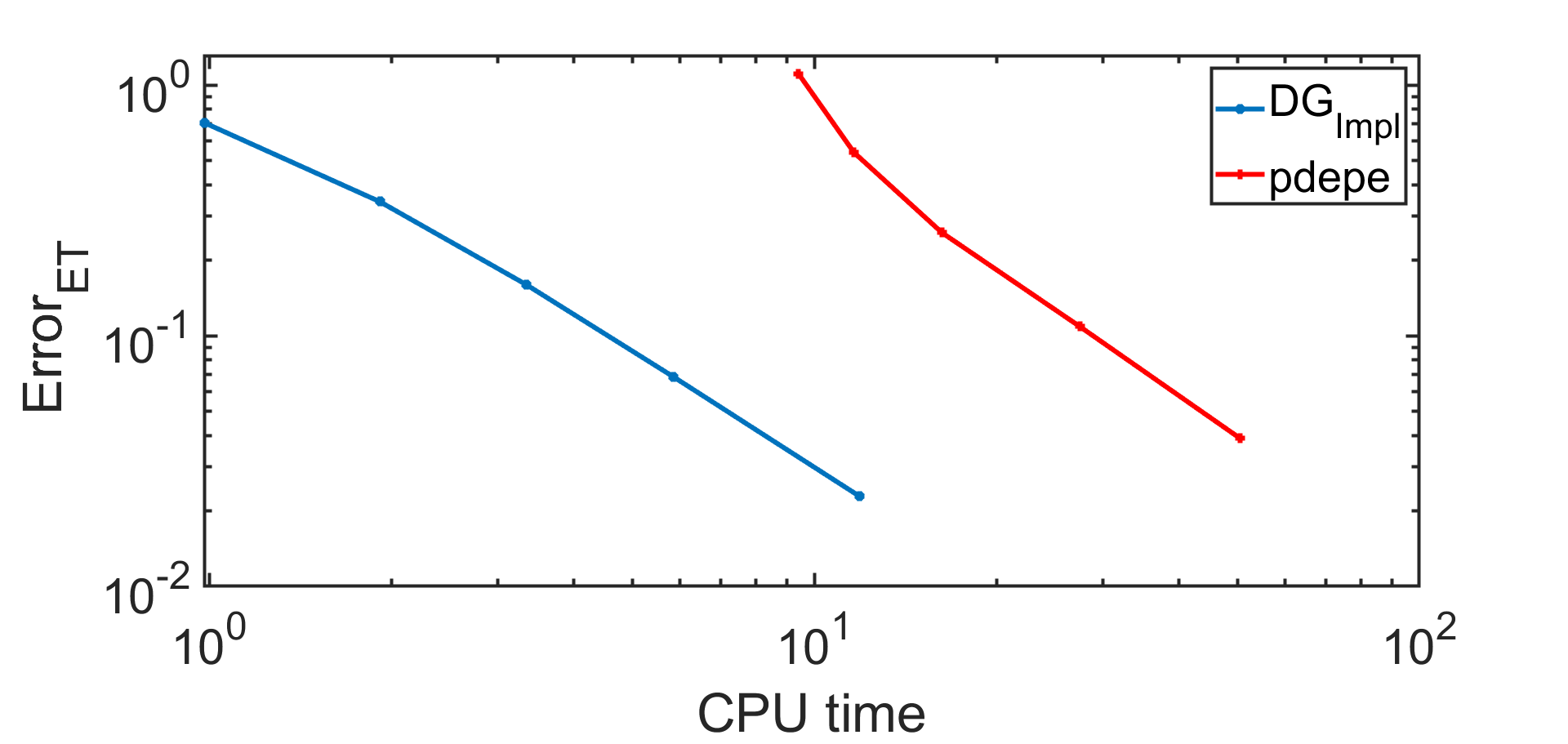}}
\caption{\textbf{Temporal convergence and computational efficiency of $\text{DG}_{\text{Impl}}$ and $pdepe$ for the ET model.}  The relative error $Error_{ETO,r}^i$ against the time step $\Delta t^i$ (a) and the computational time $\text{CPU}^i_r$ (b). $Error_{ETO,r}^i$ decreases with the time step $\Delta t^i$ implying convergence of both solvers. $\text{DG}_{\text{Impl}}$, compared to $pdepe$, takes less time to simulate the dimensionless current for a given  relative error.\label{compi1}}   
\end{figure}

The accuracy and computational efficiency with respect to the temporal discretization for ET model is examined next. For $r\in \{ pdepe,  \text{DG}_{\text{Impl}}\}$, the dimensionless current $G^i_r$ is computed using the time-step sequence
$$\Delta t^i = 2^i \times \Delta t^0, \;\; i\in \{0,\cdots,4\},$$
with $\Delta t^0=6.1\times 10^{-4}$. The solution $G_r^0$, obtained using the finest prescribed time step, $\Delta t^0$, is taken as the reference solution for solver $r$. The corresponding relative error is defined by
\begin{equation}\label{current_relER}
Error_{ET,r}^i = 100\dfrac{\Bigr| G^{i}_{r} -  G^{0}_{r}\Bigr|_{L^2(P)}}{\Bigr|  G^{0}_{r}\Bigr|_{L^2(P)}},\;\; \Bigr|  G\Bigr|_{L^2(P)}^2 = \int_P G^2 dP.
\end{equation}
Figures \ref{compi1} (a) and (b) present $Error_{ET,r}^i$ as a function of the time step $\Delta t^i$ and the corresponding computational time $\text{CPU}^i_r$, respectively. Figure \ref{compi1} (a) shows that the relative error decreases as the prescribed time step is refined, indicating convergence of both solvers toward their respective reference solutions. Figure \ref{compi1} (b) compares the computational cost required to achieve a given level of accuracy. Over the range of errors considered, $\text{DG}_{\text{Impl}}$ requires less computational time than $pdepe$. 

It should be noted that, for a prescribed $\Delta t^i$, $\text{DG}_{\text{Impl}}$ advances the solution using a uniform time step, whereas $pdepe$ employs adaptive time integration with internal time steps not exceeding the prescribed maximum step. Thus, the two solvers use different temporal integration strategies. Nevertheless, the results in Figures \ref{compi} and \ref{compi1} indicate that, for the ET model considered here, $\text{DG}_{\text{Impl}}$ produces a smoother current response and achieves the prescribed accuracy at a lower computational cost than $pdepe$.

\begin{figure}[H]
\centering
\subfloat[]{\includegraphics[height=6cm, width = 7cm]{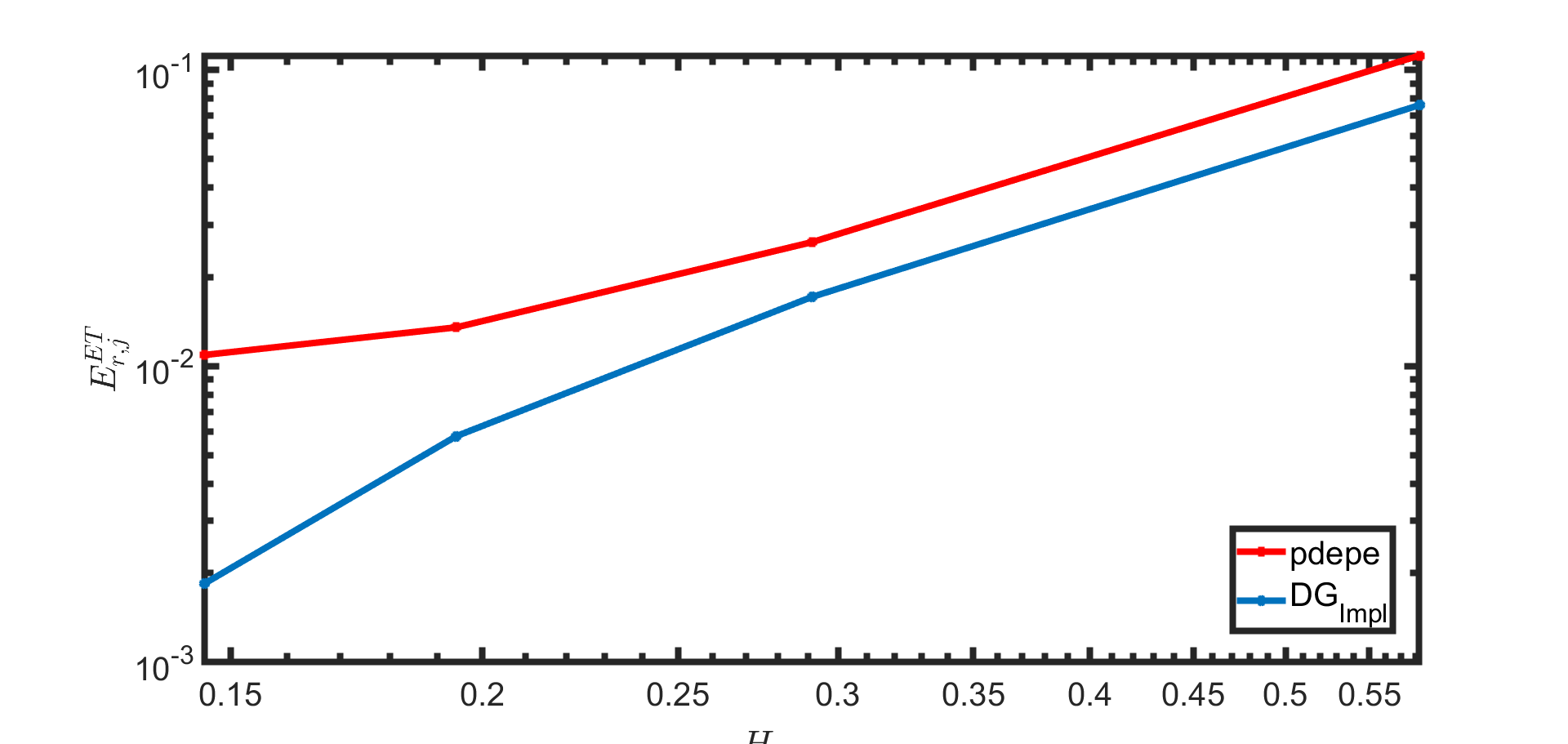}}
\subfloat[]{\includegraphics[height=6cm, width = 7cm]{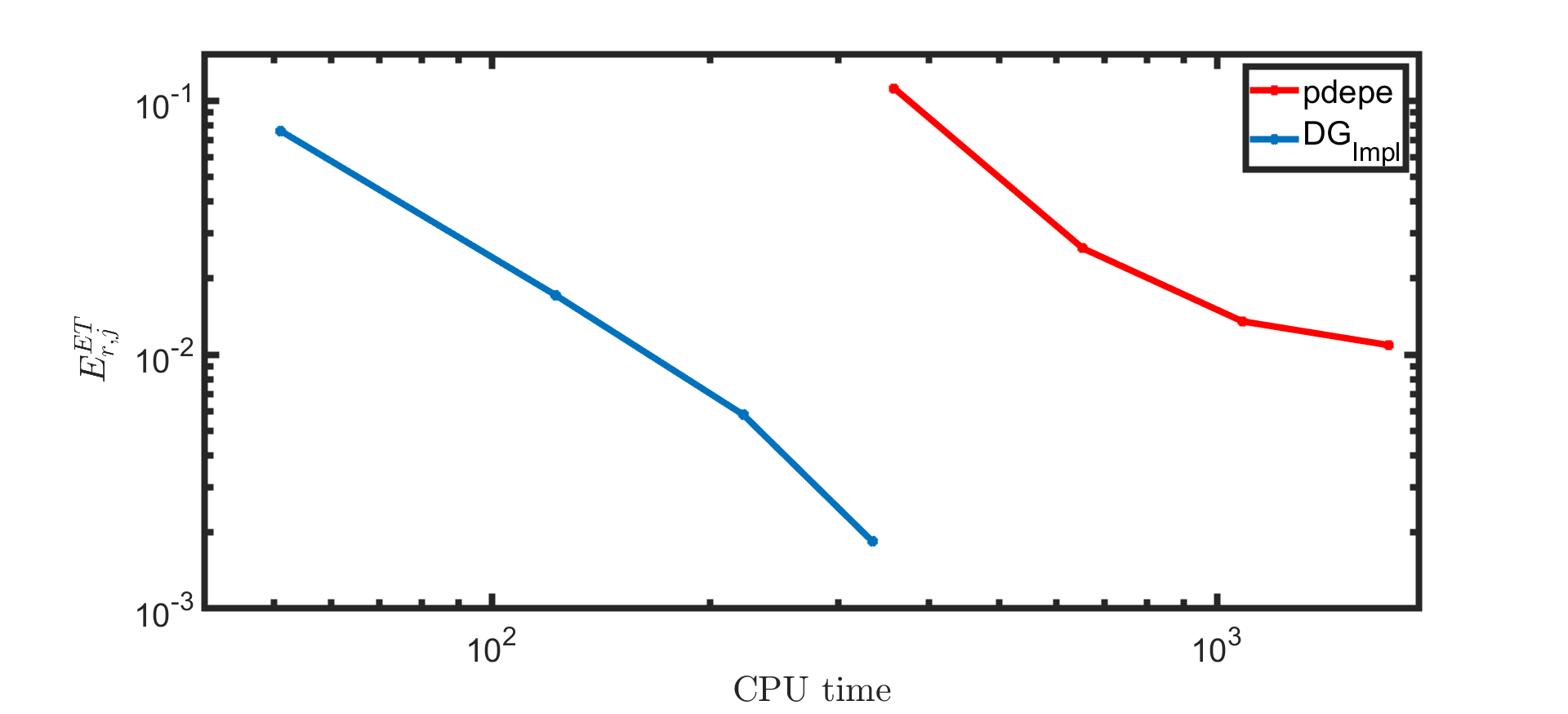}} 
\caption{\textbf{ Spatial convergence and computational efficiency of $\text{DG}_{\text{Impl}}$ and $pdepe$ for the ET model.} The relative error $E^{ET}_{r,j}$ against the mean space step  $H_j$ (a) and the computation time ($\text{CPU time}$) (b) for all $j=1,\cdots,4$ and $r\in \{ pdepe,  \text{DG}_{\text{Impl}}\}$. The relative error $E^{ET}_{r,j}$ tend to zeros with the mean space step, indicating the convergence of both solvers. Compared to $pdepe$,  $\text{DG}_{\text{Impl}}$ spends less time to simulate the dimensionless current with a relative error equal $E_0$. \label{compi11}}           
\end{figure}

We finally investigate the accuracy and computational efficiency with respect to the spatial discretization. For the solver $r\in \{ pdepe,  \text{DG}_{\text{Impl}}\}$, let \(G_r^j\) denote the dimensionless current computed on the partition \(\Omega_j\), \(j=1,\ldots,5\). Starting from the initial nonuniform partition \(\Omega_1\), each subsequent partition is obtained by subdividing every element of \(\Omega_1\) into \(j\) equal subelements. Since the resulting meshes remain nonuniform, we characterize the spatial resolution of each partition \(\Omega_j\) by its mean element size, denoted by \(H_j\). The dimensionless current \(G_r^5\), computed on the finest partition \(\Omega_5\) with mean element size \(H_5\), is considered as the reference solution for solver \(r\). The relative error associated with the partition \(\Omega_j\) is then defined as
\begin{equation}\label{cuu_error4}
E^{ET}_{r,j} = 100\dfrac{ \Bigr| G^r_j -  G^r_5 \Bigr|_{L^2(P)}}{ \Bigr| G^r_5 \Bigr|_{L^2(P)}}.
\end{equation}
For these simulations, the values of $H_5 = 0.1165$ and $\Delta t = 10^{-3}$ were used. Figures \ref{compi11} (a) and (b) show \(E_{r,j}^{ET}\) as a function of the mean element size \(H_j\) and the corresponding computational time, respectively. Figure \ref{compi11} (a) shows that the relative error decreases as \(H_j\) decreases, indicating convergence toward the respective reference solutions under successive mesh refinement for both solvers. Figure \ref{compi11} (b) compares the computational cost required by the two solvers to attain a prescribed level of accuracy. For a relative error $E_0 = 0.05\%$, $\text{DG}_{\text{Impl}}$ requires less computational time than $pdepe$. Hence, for the ET model and the family of refined nonuniform meshes considered here, $\text{DG}_{\text{Impl}}$ achieves the prescribed accuracy at a lower computational cost.

\subsection{Numerical experiments for EC' model}\label{NumEC}
The numerical assessment of $\text{DG}_{\text{Impl}}$ for EC' model follows the same procedure as for the ET model in Section \ref{sec:experiments}. $\text{DG}_{\text{Impl}}$ is first validated against known analytical and physical properties of the EC' model by examining the simulated concentration profiles of the reacting species. The comparison with MATLAB's $pdepe$ is then carried out using the same accuracy and computational efficiency measures introduced for the ET model (i.e. convergence and CPU time comparisons at prescribed relative error levels with respect to temporal/spatial discretizations).


Unless otherwise stated, the ET numerical settings in \Equaref{\eqref{ETO_defo_set}} are retained, together with the following EC' specific parameters
\begin{equation}\label{SSM2_defo_set}
K_{A}=10,\; \tilde{D}_{A} = 1,\; K_r=8,\; \tilde{D}_{B} = 1,\; \tilde{ C}_{A^{total}}^0 = 1,
\end{equation}
during the numerical simulations.

\subsubsection{Numerical validation of $\text{DG}_{\text{Impl}}$ for EC' model}
\begin{figure}[H]
  \centering
  \subfloat[$ \tilde{C}_{Q}$ for $ \tilde{D}_{Q^+} = 1$]{\includegraphics[height=2.5cm, width = 3.5cm]{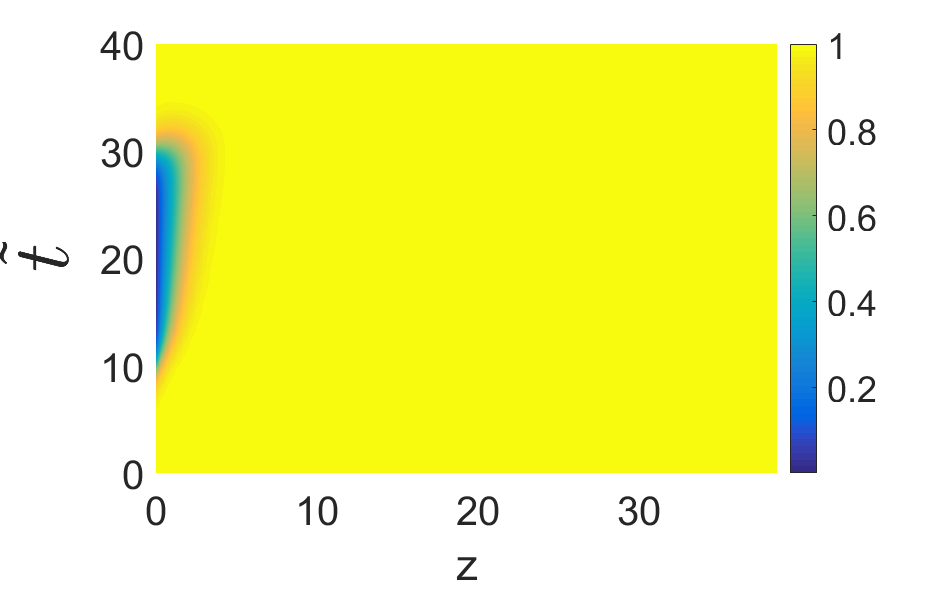}} \hspace{0.5cm}
   \subfloat[$ \tilde{C}_{Q^+}$ for $ \tilde{D}_{Q^+} = 1$]{\includegraphics[height=2.5cm, width = 3.5cm]{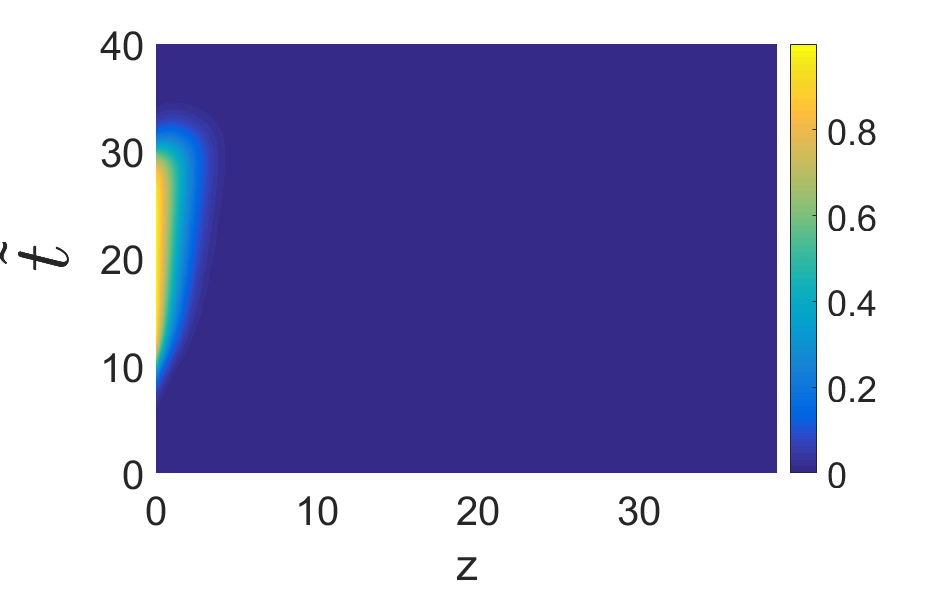}}\hspace{0.5cm}
    \subfloat[$ \tilde{C}_{B}$ for $ \tilde{D}_{Q^+} = 1$]{\includegraphics[height=2.5cm, width = 3.5cm]{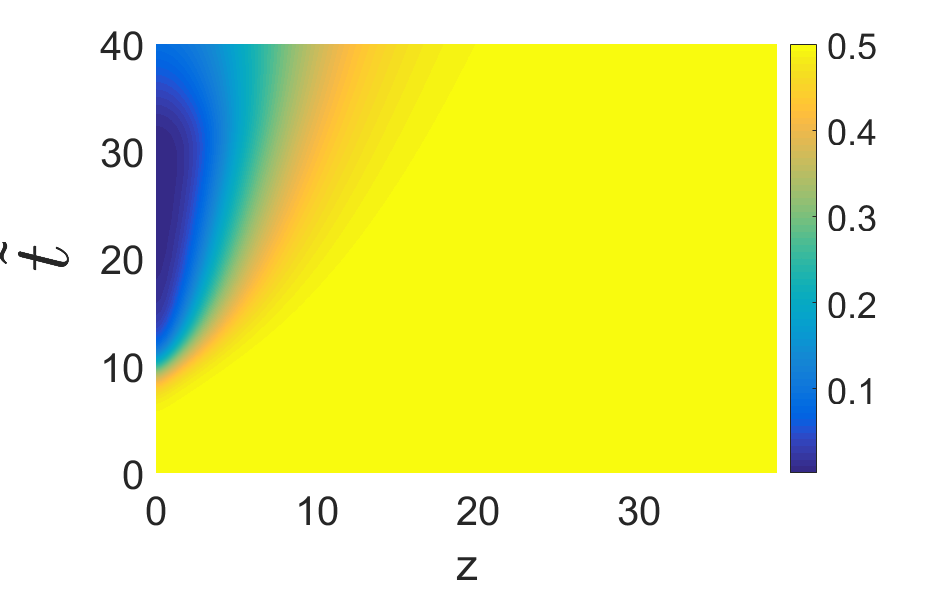}}\hspace{0.5cm}
    \subfloat[$ \tilde{C}_{A}$ for $ \tilde{D}_{Q^+} = 1$]{\includegraphics[height=2.5cm, width = 3.5cm]{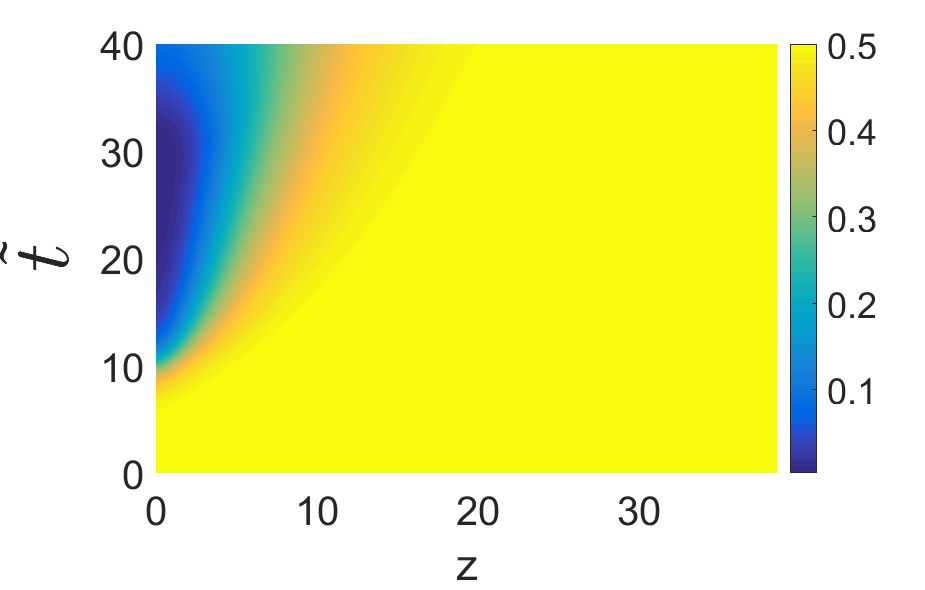}}

  \subfloat[$ \tilde{C}_{Q}$ for $ \tilde{D}_{Q^+} = 5$]{\includegraphics[height=2.5cm, width = 3.5cm]{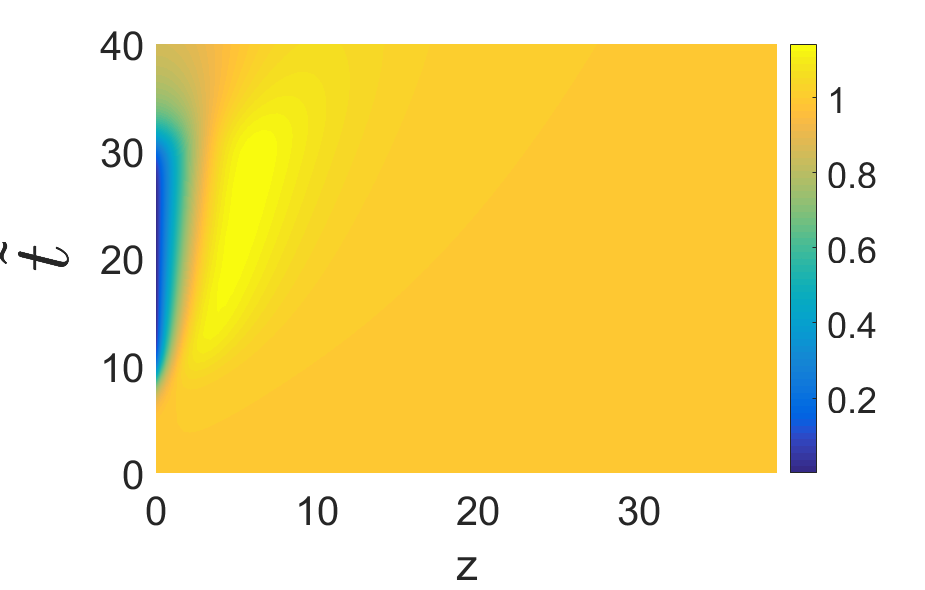}} \hspace{0.5cm}                  
  \subfloat[$ \tilde{C}_{Q^+}$ for $ \tilde{D}_{Q^+} = 5$]{\includegraphics[height=2.5cm, width = 3.5cm]{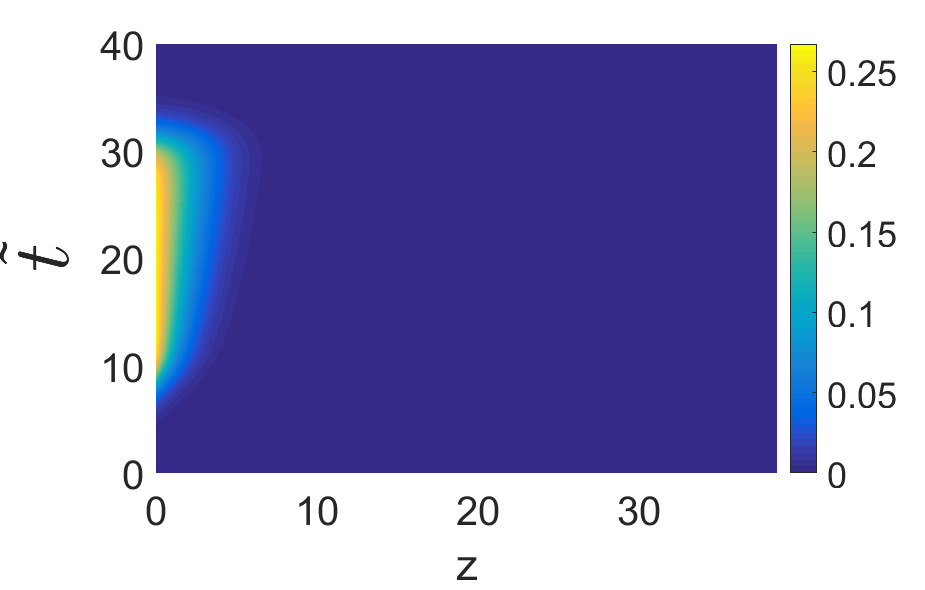}}\hspace{0.5cm}
  \subfloat[$ \tilde{C}_{B}$ for $ \tilde{D}_{Q^+} = 5$]{\includegraphics[height=2.5cm, width = 3.5cm]{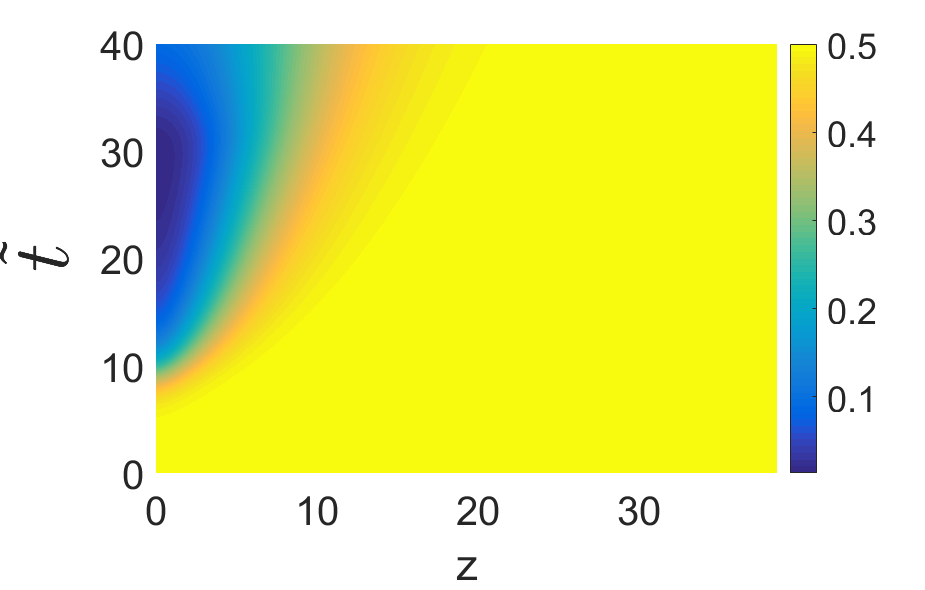}} \hspace{0.5cm}
  \subfloat[$ \tilde{C}_{A}$ for $ \tilde{D}_{Q^+} = 5$]{\includegraphics[height=2.5cm, width = 3.5cm]{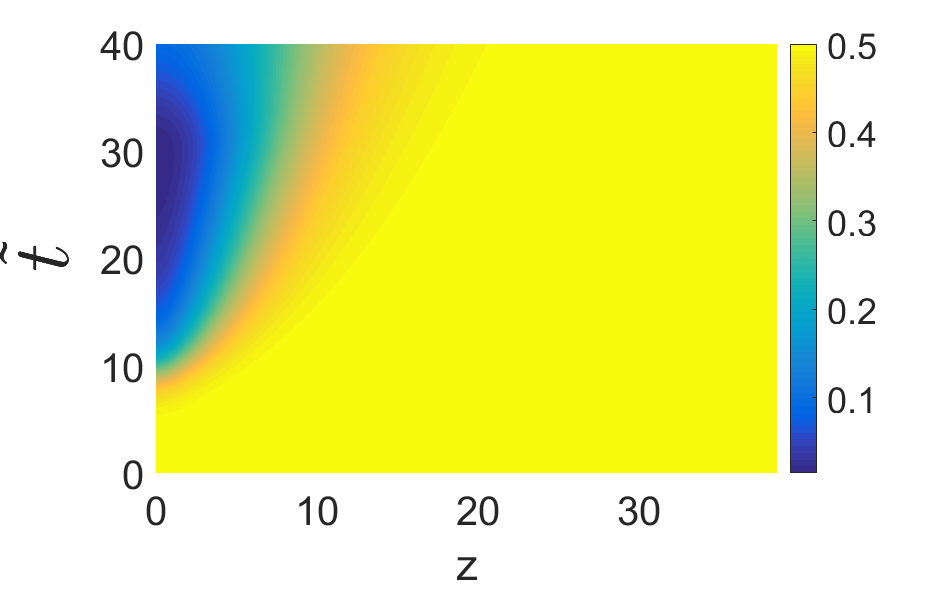}}          
  
  \subfloat[$\tilde{S}$ for $ \tilde{D}_{Q^+} = 1$]{\includegraphics[height=2.5cm, width = 3.5cm]{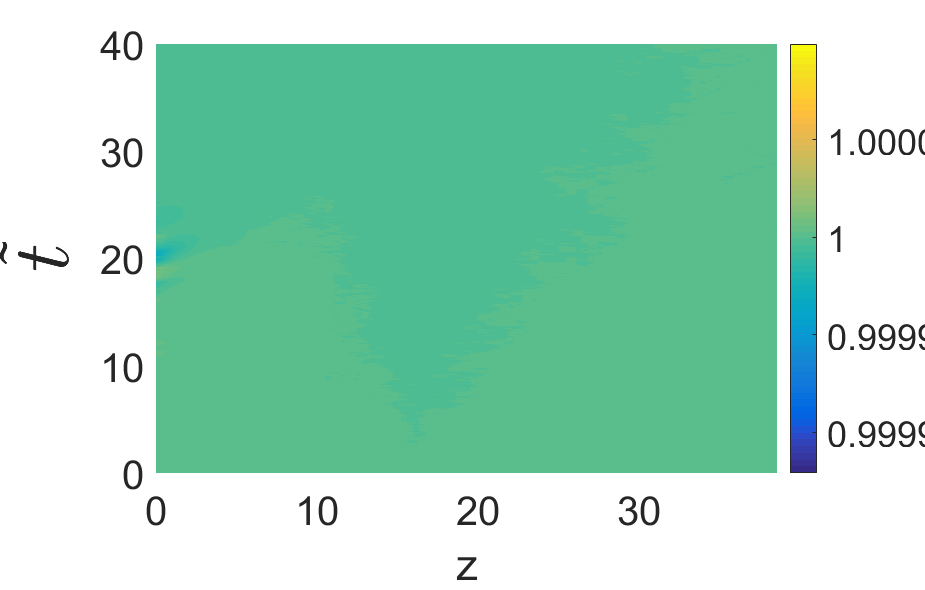}}\hspace{0.5cm}
   \subfloat[$\tilde{S}$ for $ \tilde{D}_{Q^+} = 5$]{\includegraphics[height=2.5cm, width = 3.5cm]{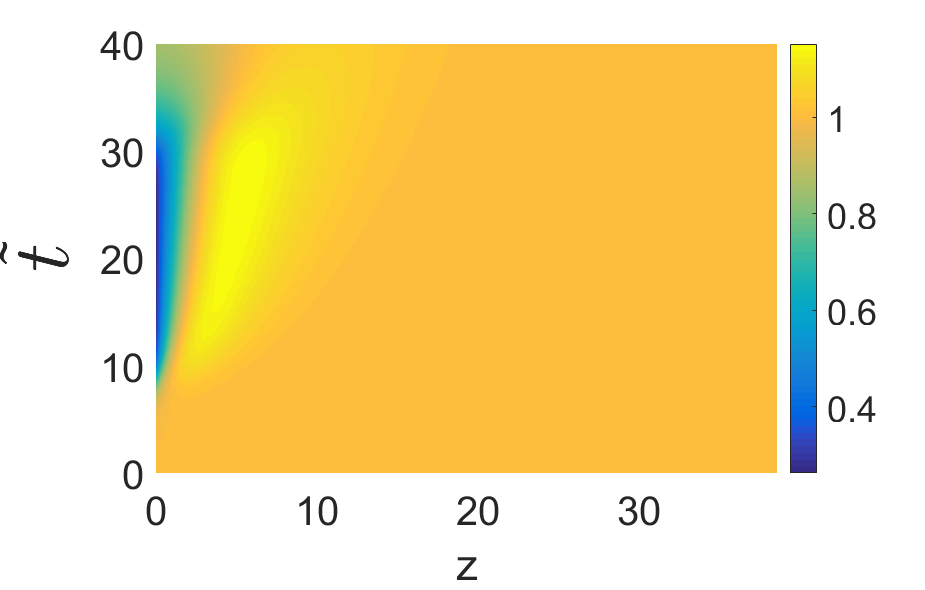}}  
  \subfloat[$R$ for $ \tilde{D}_{Q^+} = 1$]{\includegraphics[height=2.5cm, width = 3.5cm]{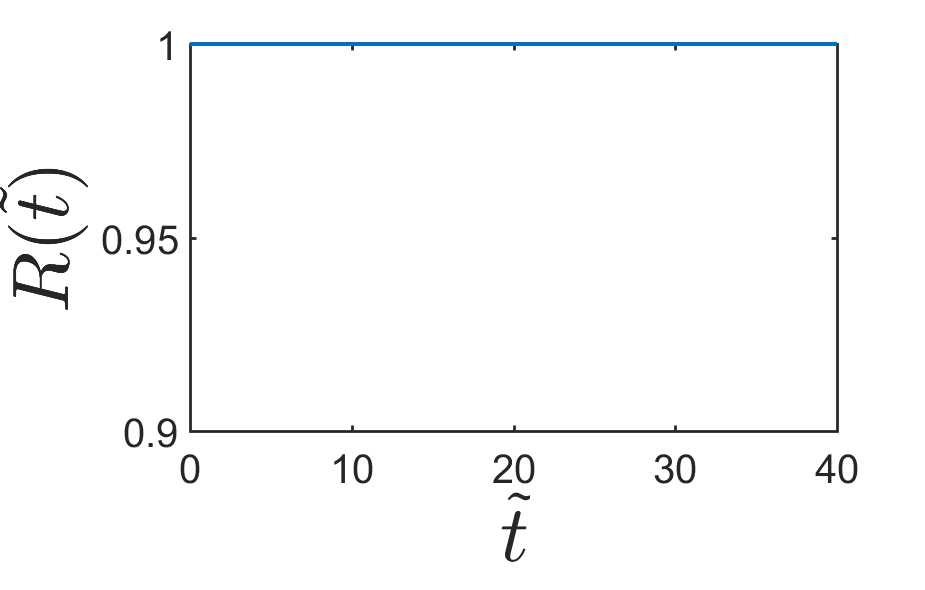}}\hspace{0.5cm}
  \subfloat[$R$ for $ \tilde{D}_{Q^+} = 5$]{\includegraphics[height=2.5cm, width = 3.5cm]{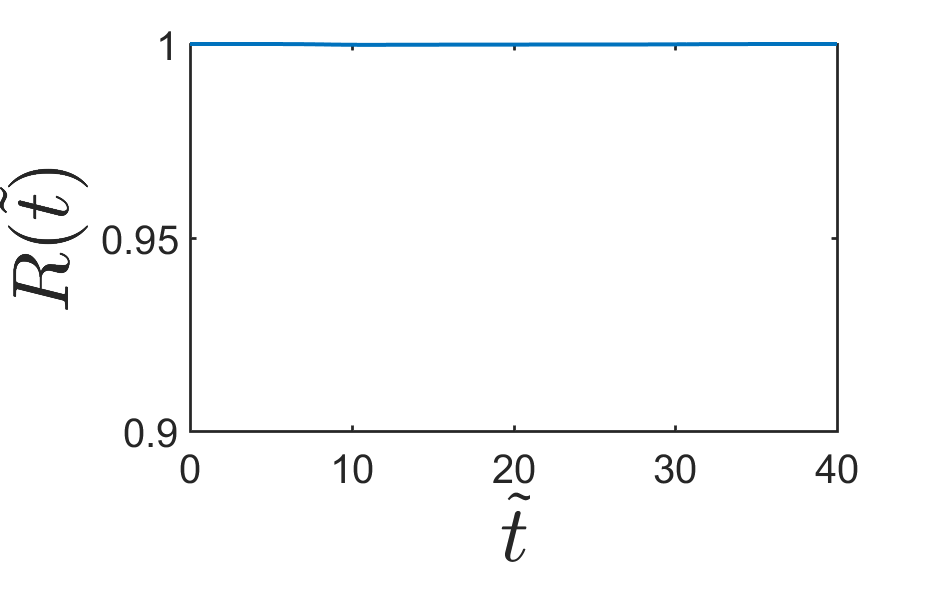}}
  \caption{\textbf{Concentration profile.} Simulated concentrations of the species $Q$, $Q^+$, $B$, $A$,  the sum  $\tilde{S}$, the function $R$ in (a), (b), (c), (d), (i), (k) and  (e), (f), (g), (h), (j), (l), for $ \tilde{D}_{Q^+} = 1$  and $ \tilde{D}_{Q^+} = 5$, respectively. Note from (k) and (l) that the mass conservation is satisfied in both cases, $\tilde{D}_{Q^+} = 1$  and $\tilde{D}_{Q^+} = 5$. (i) shows that Equation \ref{Eto_equa1} is verified for $ \tilde{D}_{Q^+} = 1$. \label{SSM2}}
\end{figure}

We validate $\text{DG}_{\text{Impl}}$ for EC' model by checking the consistency of the simulated concentrations of $Q$ and $Q^+$ with  \proref{\ref{pro_sum}} and the effect of the catalytic reaction. Figure \ref{SSM2} presents the simulated concentrations of species $Q$, $Q^+$, $B$, $A$, the sum  $\tilde{S}$ of Equation \ref{Eto_equa1}, the function $R$ of Equation \ref{Eto_equa0}, for $ \tilde{D}_{Q^+} = 1$ and $ \tilde{D}_{Q^+} = 5$.

As expected, the concentrations of the species $Q$ and $Q^+$ decrease  faster compared to ET model. This is due to the catalytic reaction shown in \Equaref{\eqref{EC}}, which involves the consumption of the species $Q^+$ and the regeneration of the species $Q$. Similary to the results for ET model, Figure \ref{SSM2} (i) and Figure \ref{SSM2} (j) show that Equation \ref{Eto_equa1} is verified for $ \tilde{D}_{Q^+} = 1$ but not for $ \tilde{D}_{Q^+} = 5$. Note from Figure \ref{SSM2} (k) and (l) that the mass conservation is satisfied in both cases, $\tilde{D}_{Q^+} = 1$  and $\tilde{D}_{Q^+} = 5$.

\subsubsection{Comparing $\text{DG}_{\text{Impl}}$ and MATLAB's $pdepe$ for EC' model}
We compare the performance of $\text{DG}_{\text{Impl}}$ and MATLAB's $pdepe$, by investigating their accuracy and computational efficiency for simulating the dimensionless current from EC' model. 

\begin{figure}[H]
\centering
\subfloat[]{\includegraphics[height=5cm, width = 7cm]{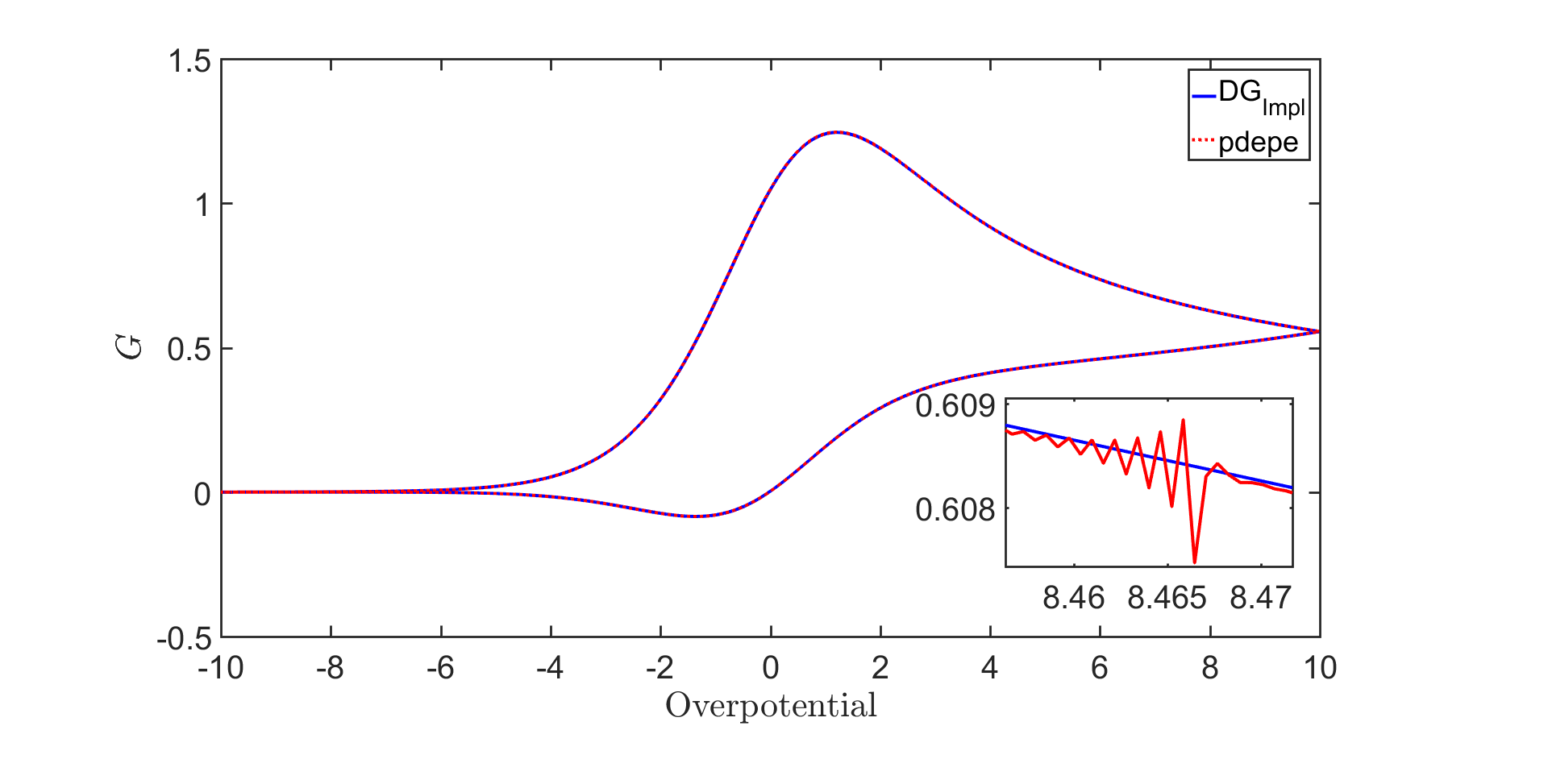}}\hspace{0.5cm}
\subfloat[]{\includegraphics[height=5cm, width = 7cm]{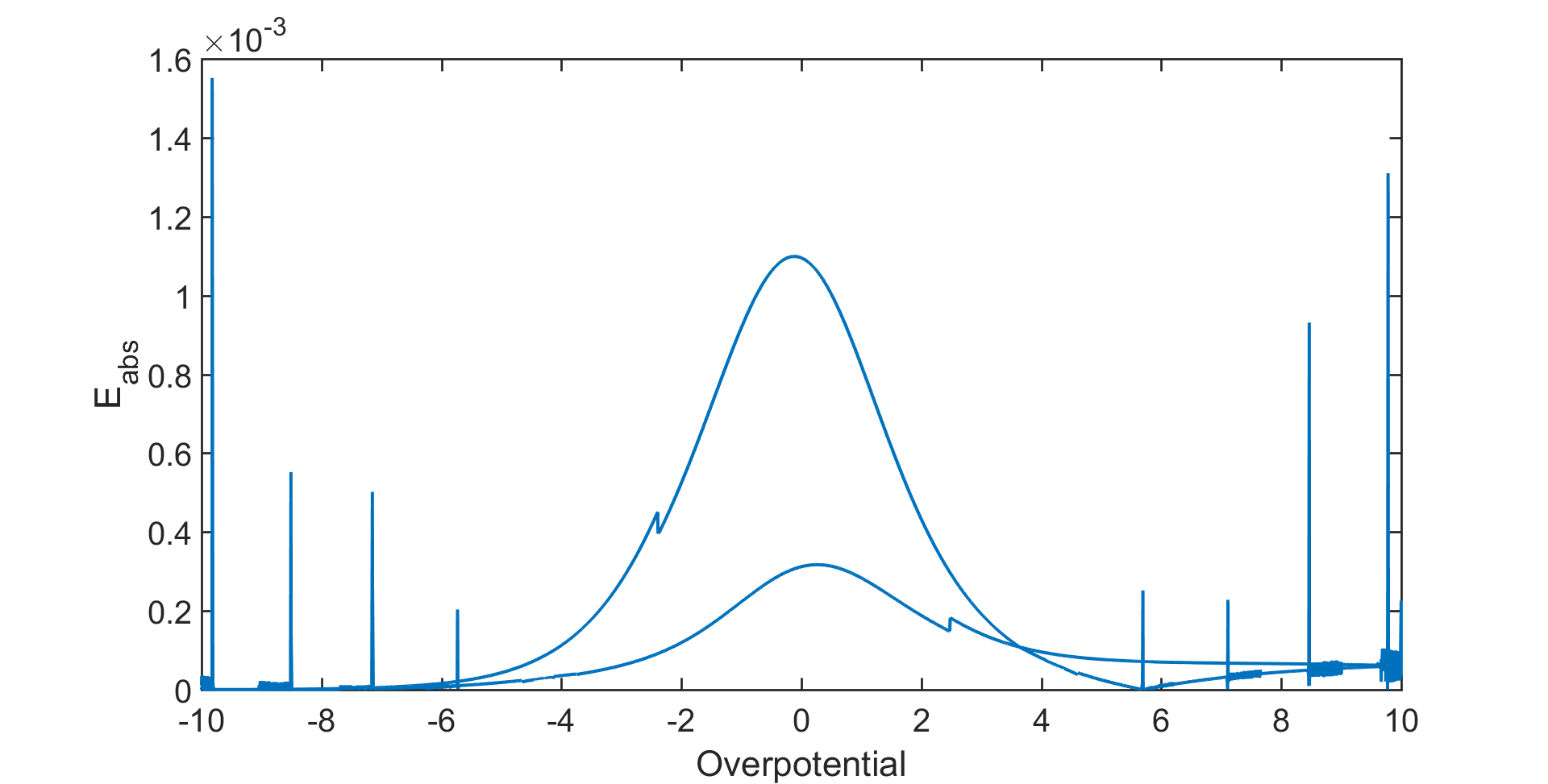}}

\subfloat[]{\includegraphics[height=5cm, width = 7cm]{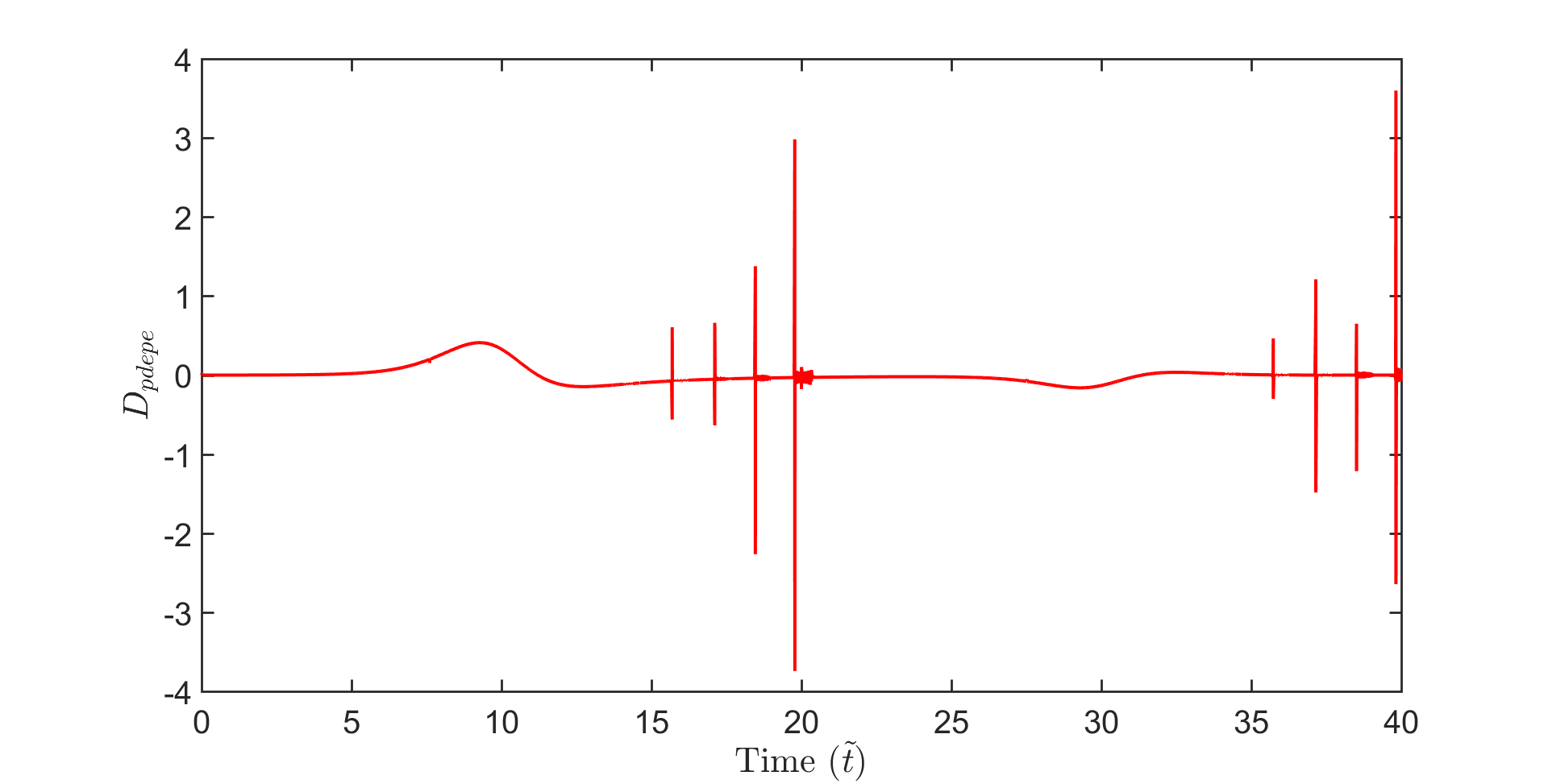}}\hspace{0.5cm}
\subfloat[]{\includegraphics[height=5cm, width = 7cm]{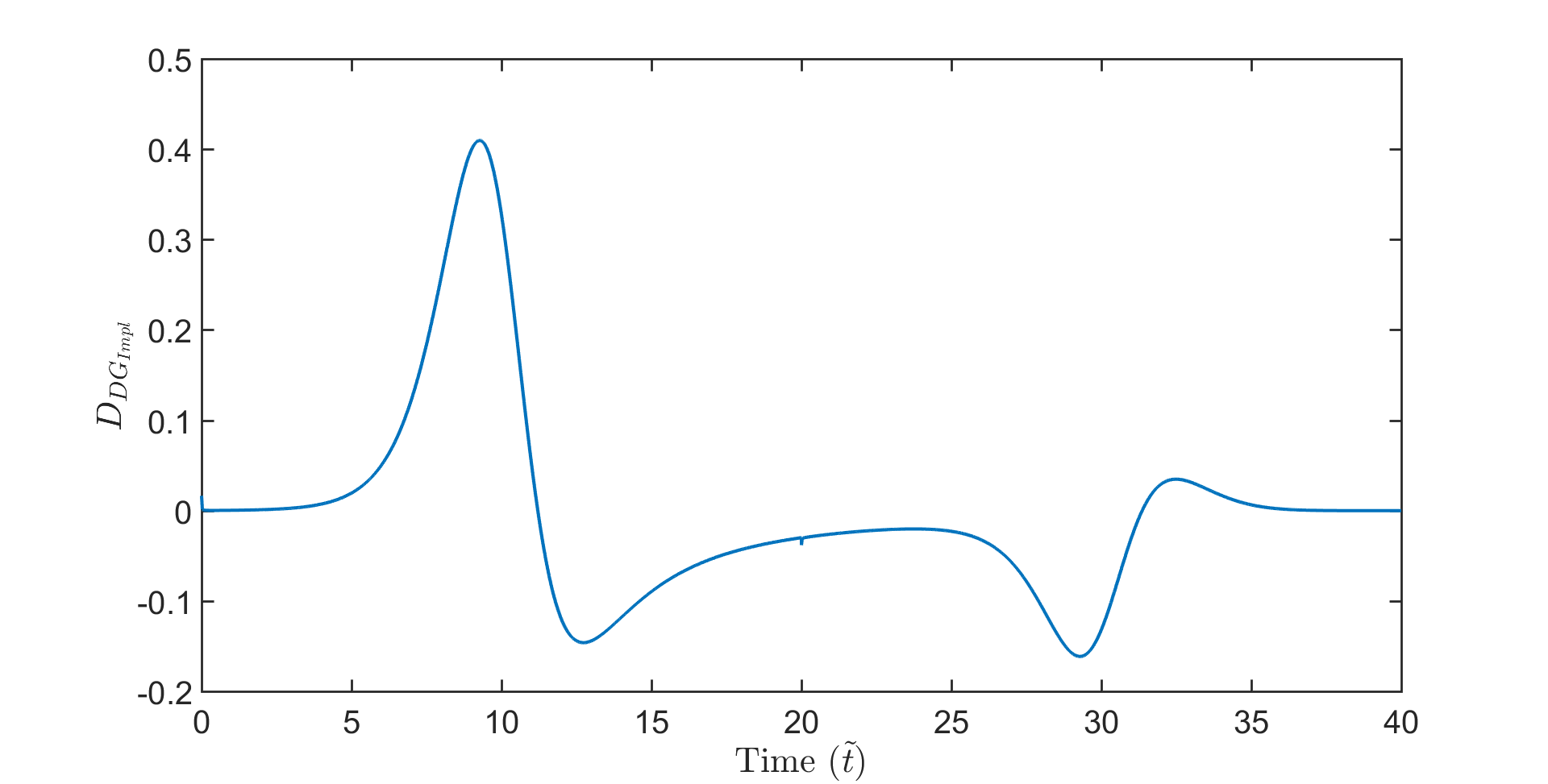}}
\caption{\textbf{Comparison of the simulated dimensionless current from EC' model.}
(a) Dimensionless currents computed using $pdepe$ and $\text{DG}_{\text{Impl}}$ as functions of the overpotential; the inset highlights small oscillations in the $pdepe$ dimensionless current.
(b) Absolute difference between the two simulated dimensionless currents, which remains below $1.6\times10^{-3}$ over the range considered.
(c)–(d) Approximate time derivatives of the currents computed using $pdepe$ and $\text{DG}_{\text{Impl}}$, respectively, as functions of the dimensionless time $\tilde{t}$. The derivative obtained from $pdepe$ exhibits pronounced oscillations, whereas the $\text{DG}_{\text{Impl}}$ derivative remains comparatively smooth.
\label{SSM2_comparison}}           
\end{figure}

First, we compare the dimensionless current simulated using the solver $r\in \{ pdepe,  \text{DG}_{\text{Impl}}\}$ under the same spatial and temporal discretization settings. Let $G^0_r$ denote the dimensionless current simulated using the solver $r\in \{ pdepe$, $  \text{DG}_{\text{Impl}}\}$ for time step $\Delta t^0 = 6.1 \times 10^{-4}$. Simulated dimensionless currents $G^0_r$ are shown in Figure \ref{SSM2_comparison} (a) while their absolute difference $E_{abs}$, given by  \Equaref{\eqref{abs_diff}}, is presented in Figure \ref{SSM2_comparison} (b).  The results show that the absolute difference $E_{abs}$ is less than $1.6\times 10^{-3}$, indicating close agreement between the two solvers for EC' model. Nevertheless, the inset in Figure \ref{SSM2_comparison} (a) reveals small oscillations in the dimensionless current obtained with $pdepe$. 

To further examine these oscillations, we approximate the time derivatives $D_{pdepe}$ and $D_{\text{DG}_{\text{Impl}}}$, of the simulated dimensionless currents using \Equaref{\eqref{varia_comp}} and present the results in Figures \ref{SSM2_comparison} (c) and (d), respectively. The derivative obtained from $pdepe$ shows noticeable oscillations, whereas the derivative associated with $\text{DG}_{\text{Impl}}$ remains comparatively smooth. These results indicate that, although the two solvers produce very similar dimensionless currents, $\text{DG}_{\text{Impl}}$ yields a smoother dimensionless current for the discretization considered.

\newpage

\begin{figure}[H]
\centering
\subfloat[]{\includegraphics[height=5cm, width = 7cm]{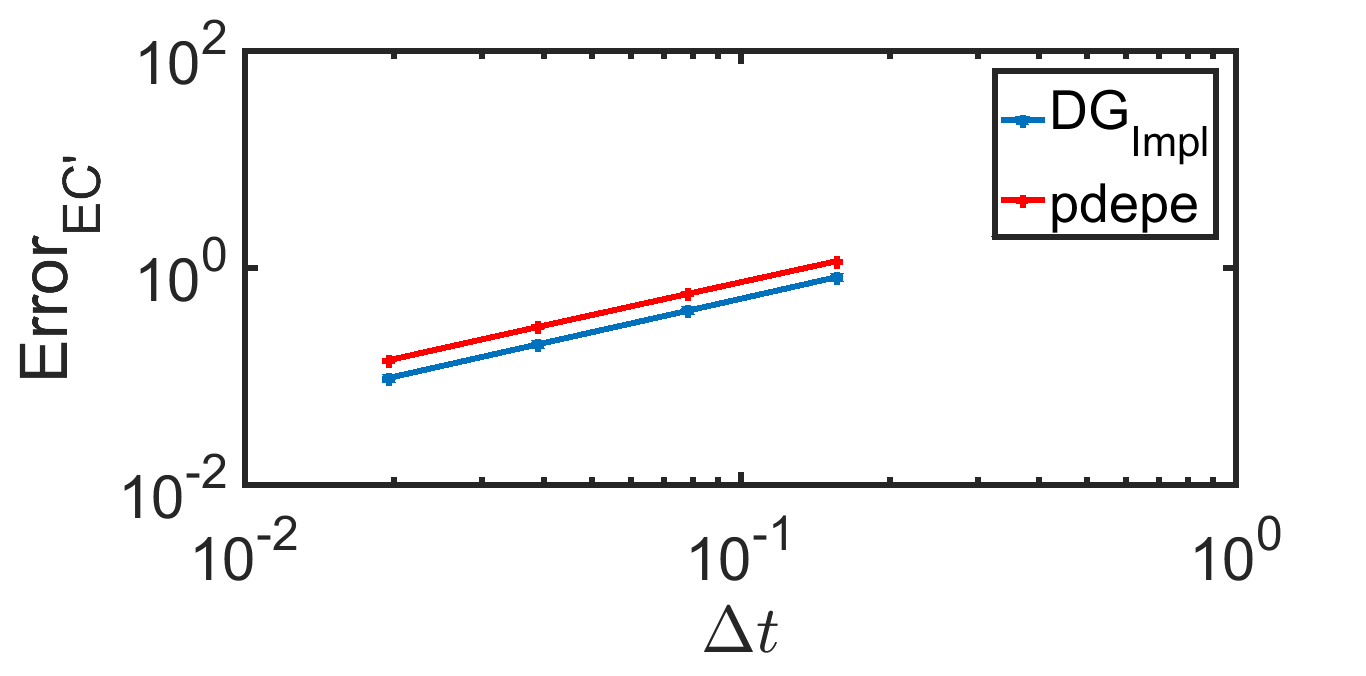}}\hspace{0.5cm}
\subfloat[]{\includegraphics[height=5cm, width = 7cm]{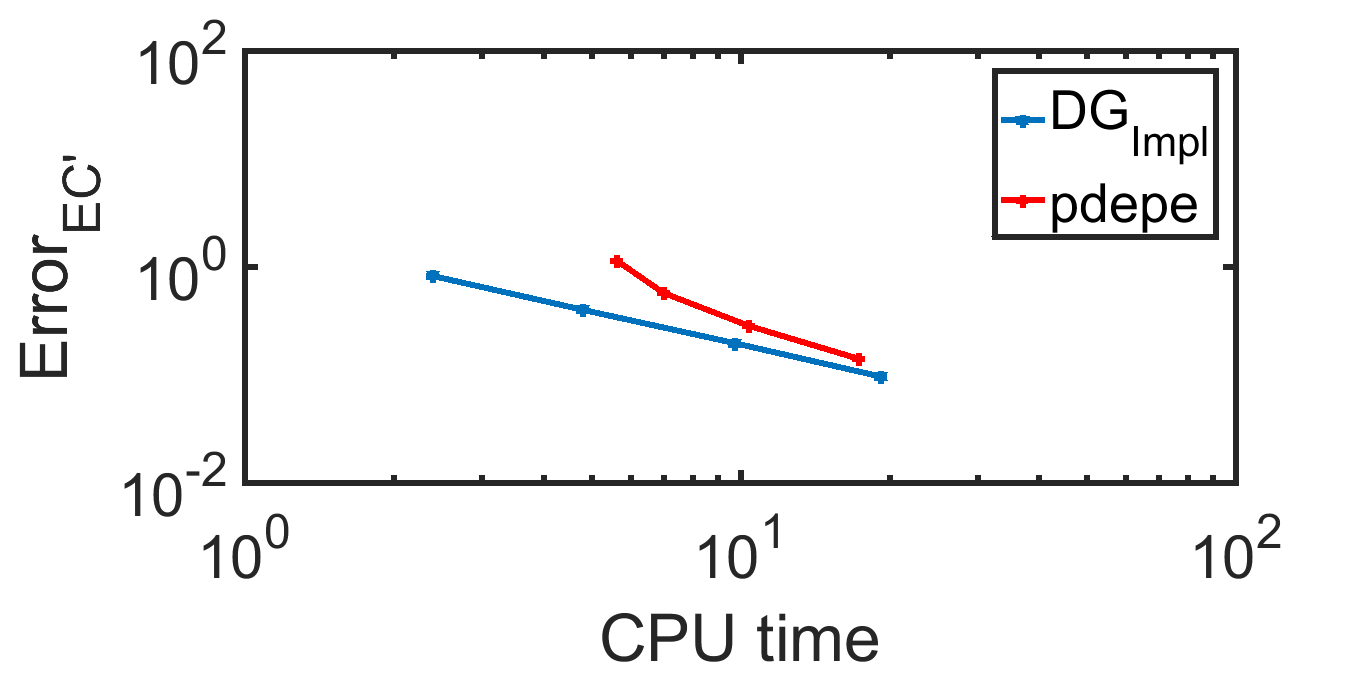}}
\caption{\textbf{Temporal convergence and computational efficiency of $\text{DG}_{\text{Impl}}$ and $pdepe$ for the EC' model.} 
(a) Relative error, $Error_{EC',r}^i$, as a function of the prescribed time step $\Delta t^i$ for both solvers. The decrease in error as $\Delta t^i$ is reduced indicates convergence of both solvers toward their respective reference currents.
(b) Relative error as a function of the corresponding CPU time. For a prescribed error level, $\text{DG}_{\text{Impl}}$ requires less computational time than $pdepe$ over the range considered.
\label{compi4}}           
\end{figure}

The accuracy and computational efficiency with respect to the temporal discretization for EC' model is examined next. For $r\in \{ pdepe,  \text{DG}_{\text{Impl}}\}$, the dimensionless current $G^i_r$ is computed using the time step sequence
$$\Delta t^i = 2^{i+4} \times \Delta t^0, \;\; i\in \{1,\cdots,4\},$$
with $\Delta t^0=6.1\times 10^{-4}$. The dimensional current, $G_r^0$, obtained using the finest prescribed time step, $\Delta t^0$, is taken as the reference current for solver $r$. We compute the relative error, $Error_{EC',r}^i$, using \Equaref{\eqref{current_relER}}.  
Figure \ref{compi4}(a) and (b) show the relative error $Error_{EC',r}^i$ against the time step $\Delta t^i$ and the computation time $CPU^i_r$, respectively. Figure \ref{compi4} (a) shows that the relative error decreases as the prescribed time step is refined, indicating convergence of both solvers toward their respective reference currents. Figure \ref{compi4} (b) compares the computational cost required to achieve a given level of accuracy. Over the range of errors considered, $\text{DG}_{\text{Impl}}$ requires less computational time than $pdepe$.

\newpage 
Finally, we investigate the accuracy and the computational efficiency of MATLAB's $pdepe$ and $\text{DG}_{\text{Impl}}$ with respect to the space discretization for EC' model. We simulate the dimensionless currents for all  partition $\Omega_j,\;j=1,\cdots, 5$  obtained by splitting each element of an initial domain $\Omega_0$ into $j=2(6-j)$ equidistant elements. By considering the dimensionless current associated to the finest mean space step $H_1$ of the partition $\Omega_1$ as the reference current, we compute the relative error $E^{EC'}_{r,j}$ associated to the mean space step $H_j$ of the partition $\Omega_j,\;j\in\{2,\cdots,5\}$ using \Equaref{\eqref{cuu_error4}}. 

\begin{figure}[H]
  \centering
  \subfloat[]{\includegraphics[height=5cm, width = 7cm]{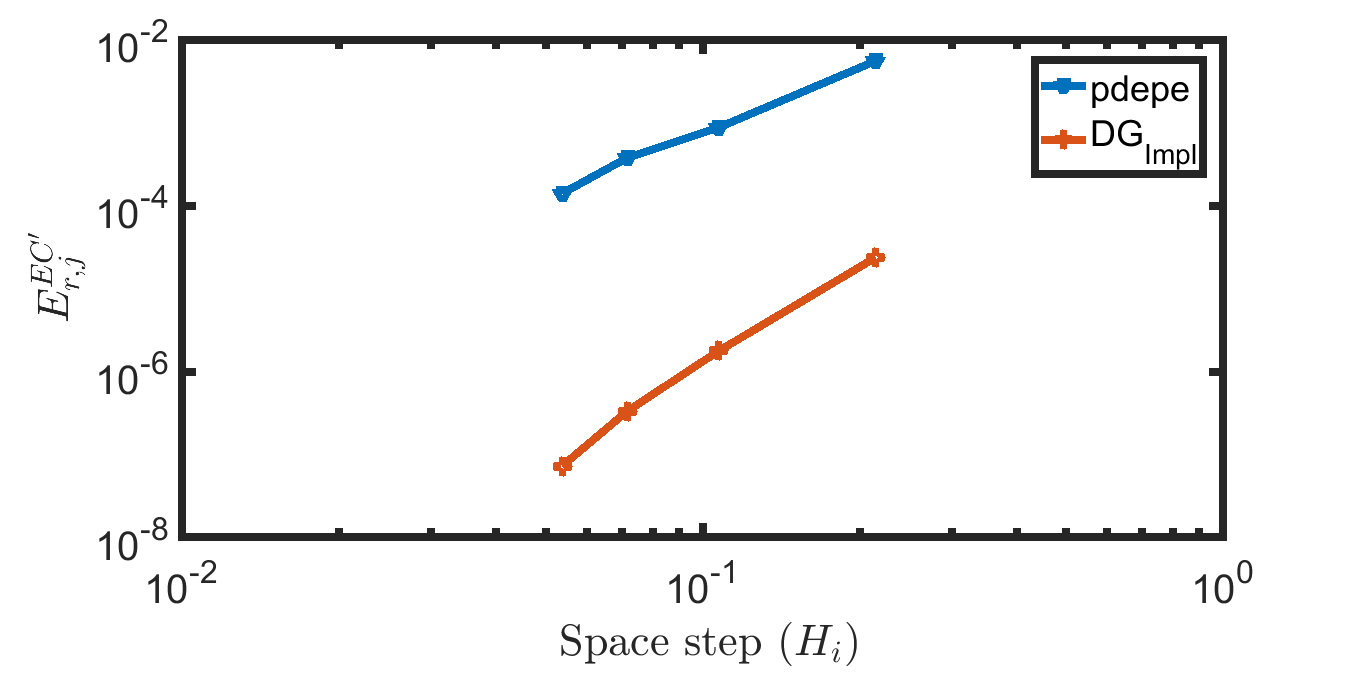}} 
\subfloat[]{\includegraphics[height=5cm, width = 7cm]{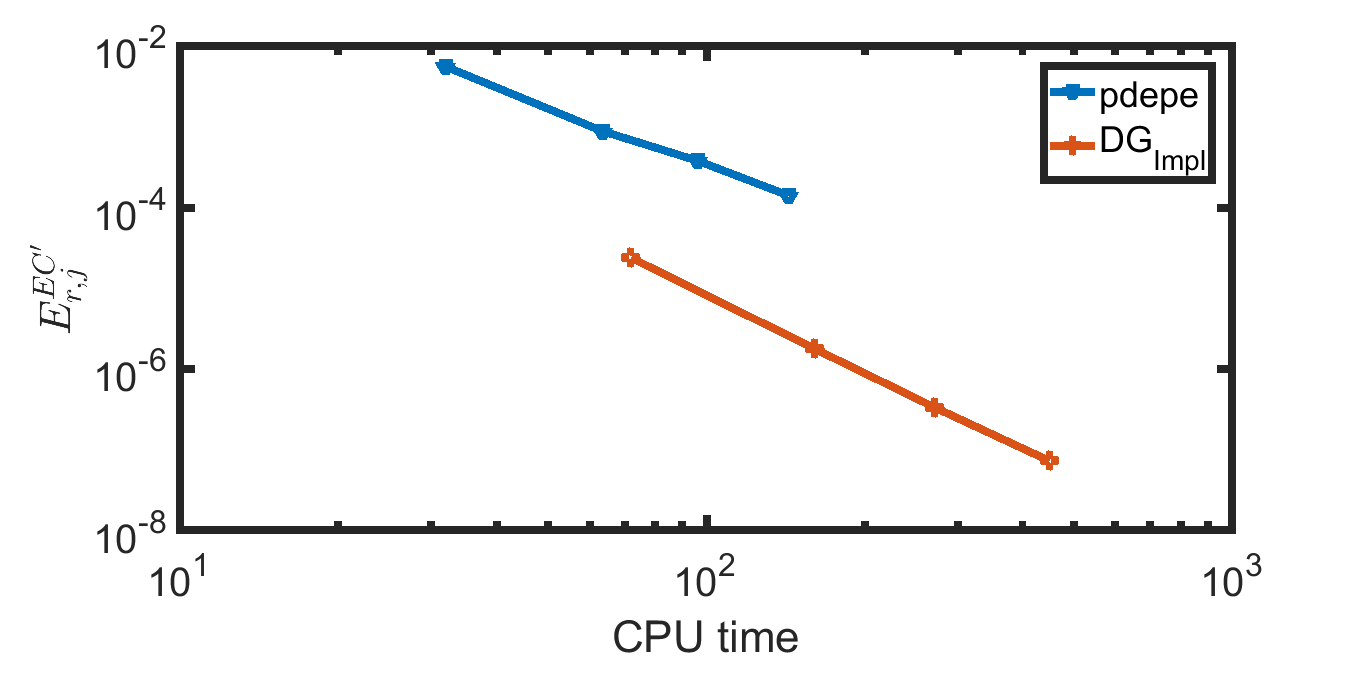}}
\caption{\textbf{Spatial convergence and computational efficiency of $\text{DG}_{\text{Impl}}$ and $pdepe$ for the EC' model}. 
(a) Relative error $E^{EC'}_{r,j}$ versus the mean element size $H_j$ for $pdepe$ and $\text{DG}_{\text{Impl}}$. Both solvers converge toward their respective reference currents as the nonuniform mesh is refined.
(b) Relative error versus CPU time, showing that $\text{DG}_{\text{Impl}}$ attains a prescribed accuracy at a lower computational cost than $pdepe$ over the refinements considered.
\label{compi_EC11}}           
\end{figure}

For these simulations, values of $H_1 = 0.0431$ and $\Delta t = 10^{-3}$ were used. Figures \ref{compi_EC11} (a) and (b) show $E^{EC'}_{r,j}$ as a function of the mean element size $H_j$ and the corresponding computational time, respectively. Figure \ref{compi_EC11} (a) shows that the relative error decreases as $H_j$ decreases, indicating convergence toward the respective reference current under successive mesh refinement for both solvers. Figure \ref{compi_EC11} (b) compares the computational cost required by the two solvers to attain a prescribed level of accuracy. For a relative error $E_0 = 0.05\%$, $\text{DG}_{\text{Impl}}$ requires less computational time than $pdepe$. Hence, for the EC' model and the family of refined nonuniform meshes considered here, $\text{DG}_{\text{Impl}}$ achieves the prescribed accuracy at a lower computational cost.

\newpage
\section{Conclusions}
\label{sec:conclusions}
We developed a new PDE solver, $\text{DG}_{\text{Impl}}$, for the efficient simulation of cyclic voltammetry (CV) responses governed by heterogeneous electron-transfer (ET) and electron-transfer coupled with a catalytic chemical reaction (EC') mechanisms. $\text{DG}_{\text{Impl}}$ employs a discontinuous Galerkin (DG) spatial discretization to transform the governing reaction diffusion PDEs into a system of ordinary differential equations (ODEs), which is subsequently advanced in time using a linearly implicit Euler method. The DG spatial discretization uses locally defined orthonormal Legendre polynomials of degree at most $k_i$ as basis functions on each mesh element $I_i$ of the nonuniform partition of the spatial domain $\Omega$.

Numerical experiments demonstrate that $\text{DG}_{\text{Impl}}$ accurately resolves transient CV reaction diffusion models, preserves relevant electrochemical properties, and reproduces the expected kinetic behavior of the ET and EC' mechanisms. For the numerical tests considered, $\text{DG}_{\text{Impl}}$ produces smooth voltammetric responses and achieves accuracy comparable to MATLAB's $pdepe$ solver at a lower computational cost. These properties make the proposed framework particularly attractive for applications requiring repeated forward simulations, such as parameter estimation, sensitivity analysis and inverse problems.

Although the present study focuses on ET and EC' mechanisms, the $\text{DG}_{\text{Impl}}$ framework is not restricted to these reaction schemes. Its modular formulation can be extended to more complex electrochemical mechanisms by incorporating additional species, homogeneous reaction terms, and appropriate electrode boundary conditions. This provides a route for treating mechanisms such as electron-transfer followed or preceded by chemical reactions (EC and CE), multiple electron-transfer and chemical steps (e.g. ECE), proton-coupled electron transfer (PCET) and models involving adsorption or desorption of electroactive species at the electrode surface. More generally, the framework can accommodate coupled reaction diffusion systems representing increasingly complex homogeneous and heterogeneous electrochemical kinetics, while retaining the same DG spatial discretization strategy.

In a follow up study, having established an efficient and extensible forward solver for CV models, we will investigate PDE-based parameter estimation and inverse problems using adjoint methods. The adjoint formulation enables efficient computation of the gradient of an objective function with respect to model parameters, making it particularly suitable for gradient-based optimization when repeated solutions of the forward problem would otherwise be computationally expensive. This extension will provide a framework for estimating kinetic and transport parameters directly from measured CV responses.

\section*{Acknowledgments}
The authors gratefully acknowledge Schlumberger Cambridge Research, where this work was carried out during the research internship of the first author, i.e. Assionvi Hove Kouevi. The work formed part of the first author's PhD program in the Department of Mathematics within the School of Mathematical and Computer Sciences (MACS) at Heriot-Watt University, whose support is also gratefully acknowledged.

\appendix
\section{Assembly of the mass and stiffness matrices}
\label{sec:append} 
In order to assemble the mass matrix, $\mathcal{M}$, we only need to assemble the block matrix $M = (M^{I_jI_i}) $, where $M^{I_jI_i}$ is a $(k_i+1)\times (k_j+1)$ matrix, for all $i,j \in\{1,\cdots,n\} $, with entries $\mathcal{A}_0(\Phi^j_r,\Phi^i_l)$, $r \in\{0,\cdots,k_j\} $ and $l \in\{0,\cdots,k_i\} $.
Since the basis function $\Phi^j_r(z) = 0$ for all $z\notin I_j$ and are orthonormal, then from \Equaref{\eqref{Mas_M}} we obtain 
\begin{equation}
\mathcal{A}_0(\Phi^j_r,\Phi^i_l) = \left\lbrace 
\begin{array}{ll}
0\;\; &\text{if } i\neq j,\;r \in\{0,\cdots,k_j\},\;l \in\{0,\cdots,k_i\} \\
\int_{I_j} \Phi^j_r\Phi^j_l= \delta_r^l \;\; &\text{if }i= j,\;r \in\{0,\cdots,k_j\},\;l \in\{0,\cdots,k_i\}
\end{array}
\right..
\end{equation}
Therefore, the block matrix $M$  is a $n_T\times n_T$ identity matrix and $\mathcal{M}$ is $(2n_T)\times (2n_T)$ identity matrix, where $n_T$ is the dimension of the finite space $V_h$, given by
\begin{equation}\label{dimen}
n_T=\sum_{j=1}^n(k_j + 1).
\end{equation}


%
%
%

In order to assemble the stiffness matrix, $\mathcal{L}^{s,\sigma_0}_{D,D^+,K_f,K_b}$, we need to assemble the block matrices $E_{K_b}, D_{K_f},A_{D,K_f}^{s,\sigma_0}$ and $A_{D^+,K_b}^{s,\sigma_0}$. 
First we investigate the assembly of $E_{K_b}$ and $D_{K_f}$. Since the entries of the first two block matrices depend only on the boundary $z_0 = 0$, then their entries will be different to zero only for the basis functions defined locally on the element $I_1$. We have for $r,l \in \{0,\cdots,k_1\}$
\begin{equation}
\mathcal{D}_{K_f}(\Phi^1_r,\Phi^1_l) = -K_f \Phi^1_r(0)\Phi^1_l(0)
,\qquad \mathcal{E}_{K_b}(\Phi^1_r,\Phi^1_l)  = -K_b \Phi^1_r(0)\Phi^1_l(0). 
\end{equation}
According to the definition of the global basis functions (see Equations \eqref{link_B} and \eqref{extend}), the above bilinear functions are equivalent to 
$$\mathcal{D}_{K_f}(\Phi^j_r,\Phi^i_l)  = -K_f P_r(-1)P_l(-1)\Gamma_{\{ r,l\}}^1
,\quad \mathcal{E}_{K_b}(\Phi^1_r,\Phi^1_l)   = -K_b P_r(-1)P_l(-1)\Gamma_{\{ r,l\}}^1,$$
where $\Gamma_S^j$ is given by  
\begin{equation}\label{norm_inv1}
\Gamma_S^j = \prod_{r\in S}\frac{1}{\parallel \phi^j_r\parallel_{L^2(I_j)}} = \prod_{r\in S}\frac{\sqrt{2r+1}}{\sqrt{h_j}}= \sqrt{h_j^{-N_S}}\underbrace{\prod_{r\in S}\sqrt{2r+1}}_{\Gamma_S}.
\end{equation}
for a given $j$ in $\{ 1,\cdots,n \}$ and set $S$ with $N_S$ elements. More generally, the non zero entries of the block matrix $E_{K_b}$ and $D_{K_f}$ are given by
\begin{equation*}
\mathcal{E}_{K_b}(\Phi^1_r,\Phi^1_l) = -K_b (-1)^{r+l}\Gamma_{\{ r,l\}}^1,\quad
\mathcal{D}_{K_f}(\Phi^1_r,\Phi^1_l) = -K_f (-1)^{r+l}\Gamma_{\{ r,l\}}^1.
\end{equation*}
Therefore, the matrices $E_{K_b}$ and $D_{K_f}$ are defined as follows
\begin{equation}\label{mat_L1}
E_{K_b} = -K_b(t)
\underbrace{ \begin{pmatrix}
  F_1 & 0        & \cdots & 0 \\
  0        & 0 & \cdots & 0 \\
  \vdots   & \vdots   & \ddots & \vdots  \\
  0       & 0        & \cdots & 0
 \end{pmatrix}}_{L1},\quad
 D_{K_f} = -K_f(t)
 \underbrace{ \begin{pmatrix}
  F_1 & 0        & \cdots & 0 \\
  0        & 0 & \cdots & 0 \\
  \vdots   & \vdots   & \ddots & \vdots  \\
  0       & 0        & \cdots & 0
 \end{pmatrix}}_{L1},
\end{equation}
where $F_1$ is $(k_1+1)\times(k_1+1)$ squared matrix with entry 
$ F_1(r,l) = (-1)^{r+l}\Gamma_{\{ r,l\}}^1$.
%
%

Finally, we investigate the assembly of $A_{D,K_f}^{s,\sigma_0}$ and $A_{D^+,K_b}^{s,\sigma_0}$. Let us introduce two bilinear function $\mathcal{B}$ and  $\mathcal{C}^{s,\sigma_0}$ such that for any function $u$, $v$ in $V_h$
\begin{equation}
\mathcal{B}_D(u,v) = D\mathcal{B}(u,v),\qquad \mathcal{C}_D^{s,\sigma_0}(u,v) = D\mathcal{C}^{s,\sigma_0}(u,v).
\end{equation}
The entries of the the matrices $A_{D,K_f}^{s,\sigma_0}$ and $A_{D^+,K_b}^{s,\sigma_0}$ are given by
\begin{align}
\mathcal{A}_{D,K_f}^{s,\sigma_0}(\Phi^j_r,\Phi^i_l) = D(\mathcal{B}(\Phi^j_r,\Phi^i_l) + \mathcal{C}^{s,\sigma_0}(\Phi^j_r,\Phi^i_l)) - \mathcal{D}_{K_f}(\Phi^j_r,\Phi^i_l),\label{mat_A}\\
\mathcal{A}_{D^+,K_b}^{s,\sigma_0}(\Phi^j_r,\Phi^i_l) = D^+(\mathcal{B}(\Phi^j_r,\Phi^i_l) + \mathcal{C}^{s,\sigma_0}(\Phi^j_r,\Phi^i_l)) - \mathcal{E}_{K_b}(\Phi^j_r,\Phi^i_l),\label{mat_A1}
\end{align}
for all $i,j \in\{1,\cdots,n \}, \; r\in\{0,\cdots,k_j \} , \; l\in\{0,\cdots,k_i \}$. Since the matrices $E_{K_b}$ and $D_{K_f}$ have previously been investigated, all that remains to do is to assemble the spatial contribution matrix, $B$, with entries $\mathcal{B}(\Phi^j_r,\Phi^i_l)$ and the interior nodal contribution matrix, $C^{s,\sigma_0}$, with entries $\mathcal{C}^{s,\sigma_0}(\Phi^j_r,\Phi^i_l)$. Thus, we can evaluate the matrices $A_{D,K_f}^{s,\sigma_0}$ and $A_{D^+,K_b}^{s,\sigma_0}$ as follows
\begin{equation*}
A_{D,K_f}^{s,\sigma_0} = D(B + C^{s,\sigma_0}) - D_{K_f}, \qquad 
A_{D^+,K_b}^{s,\sigma_0}  = D^+(B + C^{s,\sigma_0}) - E_{K_b}.
\end{equation*}
\begin{itemize}
\item Spatial contribution matrix $B$: Due to the extension of the local basis function to  global, \Equaref{\eqref{extend}}, we have
\begin{equation*}
\mathcal{B}(\Phi^j_r,\Phi^i_l) = \int_{I_j\cap I_i}d_z\Phi^j_r(z)d_z\Phi^i_l(z) \;dz.
\end{equation*}
Therefore, $\mathcal{B}(\Phi^j_r,\Phi^i_l)$ is equal to zero if $i\neq j$, since $I_i, i= 1,\cdots,n$ represent a partition. Using the set of  reference basis functions, $\{P_r,r\geq 0\}$, on the interval $I = [-1,1]$, we obtain
\begin{equation*}
\mathcal{B}(\Phi^i_r,\Phi^i_l) = \Gamma_{\{ r,l\}}^i\int_{I_i}d_z\phi^i_r(z)\phi^i_l(z)\;dz = \frac{2}{h_i^2}\underbrace{\left( \Gamma_{\{ r,l\}} \int_{I}d_xP_r(x)d_xP_l(x)\;dx\right)}_{\mathcal{B}_1(P_r,P_l)}.
\end{equation*}
The spatial contribution matrix then takes the form
\begin{equation*}
B =
 \begin{pmatrix}
  \frac{2}{h_1^2} B^{I_1} & 0        & \cdots & 0 \\
  0        & \frac{2}{h_2^2} B^{I_2} & \ddots & \vdots \\
  \vdots   & \ddots   & \ddots       & 0  \\
  0        & \cdots   & 0            & \frac{2}{h_n^2} B^{I_n}
 \end{pmatrix},
\end{equation*}
where $B^{I_i}$ is a $(k_i+1)\times(k_i+1)$ symmetric  matrix with entries $\mathcal{B}_1(P_r,P_l)$ for all $r,l=0,\cdots,k_i$. All block matrices $B^{I_i}$ can be extracted from $B^{I_{max}}$ associated to $k_{max} = \max \{ k_i, i=1,\cdots,n\}$. To assemble $B^{I_{max}}$ efficiently, we only need to compute its entries for all $l=0,\cdots,k_{max} \text{ and } l\leq r$. Using the Legendre polynomials properties, we have for $l\leq r, \;l=0,\cdots,k_{max}$ 
\begin{equation}\label{toto2}
B^{I_{max}}(r,l) = \left\lbrace 
\begin{array}{ll}
0 \;\; \qquad\quad &\text{if }\;  r+l \text{ is odd }\\
 l(l+1)\Gamma_{\{ r,l\}},\;\; &\text{otherwise}
\end{array}
\right..
\end{equation}
\item Interior nodal contribution matrix $C^{s,\sigma_0}$: Due to the extension of the local basis function to global, \Equaref{\eqref{extend}}, we have
\begin{equation*}
\mathcal{C}^{s,\sigma_0}(\Phi^j_r,\Phi^i_l) \neq 0 \quad\text{ if and only if  }\quad I_i\cap I_j \neq \emptyset.
\end{equation*}
Hence, the interior nodal contribution matrix takes the form
\begin{equation*}
C^{s,\sigma_0} =
 \begin{pmatrix}
 C_{I_1I_1}^{s,\sigma_0} & C_{I_2I_1}^{s,\sigma_0}  & 0       & \cdots & 0 \\\\
 C_{I_1I_2}^{s,\sigma_0} & C_{I_2I_2}^{s,\sigma_0}  & C_{I_3I_2}^{s,\sigma_0}  & \ddots & \vdots \\\\
  0                       &  \ddots                   & \ddots  & \ddots & 0  \\\\
 \vdots                   & \ddots                    & C_{I_{n-2}I_{n-1}}^{s,\sigma_0}  & C_{I_{n-1}I_{n-1}}^{s,\sigma_0} & C_{I_nI_{n-1}}^{s,\sigma_0}\\\\
  0                       & \cdots                    & 0  & C_{I_{n-1}I_{n}}^{s,\sigma_0} & C_{I_nI_{n}}^{s,\sigma_0} 
 \end{pmatrix},
\end{equation*}
where the block matrix $C_{I_jI_i}^{s,\sigma_0}$ is a $(k_i+1)\times(k_j+1)$  matrix, 
which according to  \Equaref{\eqref{bil_S}}, is defined as follows
\begin{equation*}
C_{I_jI_i}^{s,\sigma_0}  = \left(\mathcal{C}^{s,\sigma_0}(\Phi^j_r,\Phi^i_l) \right)_{r=0,\cdots,k_j}^{l=0,\cdots,k_i}, \quad \mathcal{C}^{s,\sigma_0}(\Phi^j_r,\Phi^i_l) = \sum_{m=1}^{n-1}\mathcal{C}_m^{s,\sigma_0}(\Phi^j_r,\Phi^i_l).
\end{equation*}
The bilinear function $ \mathcal{C}_i^{s,\sigma_0}$ is defined at the interior node $z_i$ and each interior node $z_i$ contributes to four blocks of the global matrix,  $C^{s,\sigma_0}$, which are
\begin{equation*}
A^i = \begin{pmatrix}
 C_{I_iI_i}^{s,\sigma_0} & C_{I_{i+1}I_i}^{s,\sigma_0}\\
C_{I_iI_{i+1}}^{s,\sigma_0} & C_{I_{i+1}I_{i+1}}^{s,\sigma_0}
\end{pmatrix}.
\end{equation*}
The local stiffness matrix stemming from the interior node can be block-partitioned in the form
\begin{equation}
C^i = \begin{pmatrix}
C^i_{I_iI_i} & C^i_{I_{i+1}I_i}\\
C^i_{I_iI_{i+1}} & C^i_{I_{i+1}I_{i+1}}
\end{pmatrix},\;\;C^i_{I_{r1}I_{r2}} = \left(\mathcal{C}_i^{s,\sigma_0}(\Phi^{r1}_r,\Phi^{r2}_l) \right)_{r=0,\cdots,k_{r1}}^{l=0,\cdots,k_{r2}},
\end{equation}
with $r1, r2\in \{i,i+1\}$, such that the matrix $A^i$ is updated by
$$ A^i \longrightarrow A^i + C^i,\;\;\text{ for }\;i=1,\cdots,n-1 .$$
Combining the formula of the bilinear function $\mathcal{C}_i^{s,\sigma_0}$, in  \Equaref{\eqref{bil_S}}, the definition of our global basis function in Equations \eqref{link_B}, \eqref{extend} and the properties of Legendre polynomials; we obtain 
\begin{equation}\label{stem_mat}
\left. 
\begin{array}{llll}
\mathcal{C}_{i}^{s,\sigma_0}(\Phi^i_r,\Phi^i_l) &= \left(  - \dfrac{r(r+1)}{2h_i}  +s \dfrac{l(l+1)}{2h_i} + \sigma_i\right) \Gamma_{\{ r,l\}}^i, \\
\mathcal{C}_{i}^{s,\sigma_0}(\Phi^i_r,\Phi^{i+1}_l) &= \left( \dfrac{r(r+1)}{2h_i} - s \dfrac{l(l+1)}{2h_{i+1}} - \sigma_i\right) (-1)^l\Gamma_{\{ r\}}^{i}\Gamma_{\{l\}}^{i+1},\\
\mathcal{C}_{i}^{s,\sigma_0}(\Phi^{i+1}_r,\Phi^i_l) &= \left( \dfrac{r(r+1)}{2h_{i+1}} - s \dfrac{l(l+1)}{2h_i} - \sigma_i\right) (-1)^r\Gamma_{\{ r\}}^{i+1}\Gamma_{\{ l\}}^{i},\\
\mathcal{C}_{i}^{s,\sigma_0}(\Phi^{i+1}_r,\Phi^{i+1}_l) &= \left( - \dfrac{r(r+1)}{2h_{i+1}}+ s \dfrac{l(l+1)}{2h_{i+1}} + \sigma_i\right) (-1)^{l+r}\Gamma_{\{ r,l\}}^{i+1},
\end{array}
\right. 
\end{equation}
\end{itemize}
By taking $L_0^{s,\sigma_0} = B + C^{s,\sigma_0}$, the stiffness matrix is then defined as follows
\begin{equation}\label{ETO_matr}
\mathcal{L}^{s,\sigma_0}_{D,D^+,K_f,K_b}
=\underbrace{\left( 
\begin{array}{c|c}
DL_0^{s,\sigma_0} + K_f(\tilde{t})L1&-K_b(\tilde{t})L1\\\hline
-K_f(\tilde{t})L1&D^+L_0^{s,\sigma_0}+K_b(\tilde{t})L1
\end{array}
\right)}_{\mathcal{L}_1(\tilde{t})}.
\end{equation}


\section{Computation of the source term and its Jacobian\label{sec:append1}}
The source term $N(X,K_r)$ can easily be formed once the components of the vector $N_1(X)$ are computed. In order to compute these components, we first rewrite \Equaref{\eqref{source}} as a linear combination of the integral of the triple product of Legendre polynomials \Equaref{\eqref{source1}}, where the coefficients are in terms of the components of $X$
\begin{align}
N_{1,l}^j =& \int_\Omega \tilde{ C}_{Q^{+}}\tilde{ C}_{B}\Phi^j_l(z)dz
=\int_{I_j} \tilde{ C}_{Q^{+}}\tilde{ C}_{B}\Phi^j_l(z)dz\nonumber,\\
=& \int_{I_j} \left( \sum_{r=0}^{k_j}\alpha^{+,j}_r\Phi^j_r(z)\right) \left( \sum_{m=0}^{k_j}\beta^{-,j}_m\Phi^j_m(z)\right) \Phi^j_l(z)dz\nonumber,\\
=& \sum_{r=0}^{k_j}\sum_{m=0}^{k_j} \alpha^{+,j}_r\left(  \int_{I_j}\Phi^j_r\Phi^j_m\Phi^j_l\right) \beta^{-,j}_m
\nonumber,\\
=& \sum_{r=0}^{k_j}\sum_{m=0}^{k_j} \alpha^{+,j}_r\left( \frac{h_j\Gamma_{\{ r,m,l\}}^j }{2} \int_{I}P_rP_mP_l\right) \beta^{-,j}_m\nonumber,\\
N_{1,l}^j=& \frac{1}{2\sqrt{h_j}}\sum_{r=0}^{k_j}\sum_{m=0}^{k_j} \alpha^{+,j}_r\left( \Gamma_{\{ r,m,l\}} \int_{I}P_rP_mP_l\right) \beta^{-,j}_m.\label{source2}
\end{align}
By introducing the $(k_j+1)\times (k_j+1)$ matrix, denoted $\mathbb{P}^j_l$, with the entries
$$ \mathbb{P}^j_l(r,m) = \Gamma_{\{ r,m,l\}} \int_{I}P_rP_mP_l,\qquad \forall r,m\in \{0,\cdots,k_j\},$$
then \Equaref{\eqref{source2}} can be rewritten as follows
\begin{equation}\label{source1}
N_{1,l}^j= \frac{1}{2\sqrt{h_j}} \alpha^{+,j^T}\mathbb{P}^j_l \beta^{-,j},\; \forall j\in\{1,\cdots,n\},\;\forall l\in\{1,\cdots,k_j\}.
\end{equation}

It is well known,  see \cite{LEG, LEG1} for more details, that the product of two Legendre polynomials can be expanded, in a series using the special Wigner $3j$ symbols, as follows
\begin{equation}\label{Wig}
P_iP_m = \sum_{q =\mid i-m\mid}^{i+m}\left( 
\begin{array}{ccc}
i & m & q\\
0 & 0 & 0
\end{array}
\right)^2(2q+1)P_q,
\end{equation}
where the special Wigner $3j$ symbols are computed by the following formula
\begin{equation*}
\left( 
\begin{array}{ccc}
i & m & q\\
0 & 0 & 0
\end{array}
\right) = \left\lbrace 
\begin{array}{ll}
\sqrt{\frac{(2g - 2i)!(2g - 2m)!(2g - 2q)!}{(2g+1)!}}\frac{g!(-1)^g}{(g - i)!(g - m)!(g - q)!} \;\; &\text{if }\;  J=2g\\
0\;\; &\text{if }\;  J=2g +1
\end{array}
\right.,
\end{equation*}
for $J = i +m +q$. Therefore, the integral of the triple product of Legendre polynomials is given by 
$$ \int_{I}P_iP_mP_l = 2\left( 
\begin{array}{ccc}
i & m & l\\
0 & 0 & 0
\end{array}
\right)^2, \qquad \text{if }\;\; \mid m-i\mid\leq l \leq m+i.$$
Using the symmetry of the triple product and the  orthogonality of the Legendre polynomials, the vector $N_1(X)$ can be efficiently computed.

Since the entries of the vector, $N_1(X)$ defined by \Equaref{\eqref{source1}}, depend only on the concentration of the species $Q^{+}$ and $B$ then the Jacobian of the source term $N(X,K_r)$ is given by
\begin{equation}
 \partial_XN(X,K_r) = K_r\left( 
\begin{array}{cccc}
0 & 2J_1 & 2J_2& 0\\
0 & -2J_1& -2J_2& 0\\
0 & -J_1 & -J_2& 0\\
0&0&0&0
\end{array}
\right),
\end{equation}
where the block diagonal matrices  $J_1$ and $J_2$ respectively represent the Jacobian of the vector $N_1(X)$ with respect to the entries $\alpha^{+}$ and $\beta^{-}$ of the vector $X$. According to the local definition of the basis function \Equaref{\eqref{extend}}, the block diagonal matrices  $J_1$ and $J_2$ are given by
\begin{equation}
J_1  =
 \begin{pmatrix}
  J_1^1 & 0        & \cdots & 0 \\
  0        & J_1^2 & \ddots & \vdots \\
  \vdots   & \ddots   & \ddots & 0  \\
  0       & \cdots        & 0 & J_1^n
 \end{pmatrix}, \qquad
 J_2 = \begin{pmatrix}
  J_2^1 & 0        & \cdots & 0 \\
  0        & J_2^2 & \ddots & \vdots \\
  \vdots   & \ddots   & \ddots & 0  \\
  0       & \cdots        & 0 & J_2^n
 \end{pmatrix},
\end{equation}
with the $(k_j+1)\times(k_j+1)$ matrices  $J_1^j$ and $J_2^j$ entries given by
\begin{equation}\label{source3}
J_1^j(m,l) = \sum_{i=0}^{k_j}\frac{\beta^{-,j}_i}{2\sqrt{h_j}}\Gamma_{\{ i,m,l\}}\int_IP_iP_lP_m,\quad J_2^j(m,l) = \sum_{i=0}^{k_j}\frac{\alpha^{+,j}_i}{2\sqrt{h_j}}\Gamma_{\{ i,m,l\}}\int_IP_iP_lP_m.
\end{equation}

\bibliographystyle{abbrv}  
\bibliography{references}  

%
%
%
%

\end{document}